\documentclass[10pt,reqno]{amsart}
\usepackage{amssymb,mathrsfs,graphicx,extpfeil,amsmath}
\usepackage{ifthen,latexsym,float,colortbl}

\usepackage{enumerate}

\usepackage[margin=1in]{geometry}
\usepackage{caption}
\usepackage{sidecap}
\usepackage[bookmarks=false]{hyperref}
\usepackage{rotating,dsfont,stackengine}
\usepackage{colortbl}
\definecolor{black}{rgb}{0.0, 0.0, 0.0}
\definecolor{red}{rgb}{1.0, 0.5, 0.5}
\provideboolean{shownotes} % define a boolean variable
\setboolean{shownotes}{true} % defined as ``true'' or ``false''
\newcommand{\margnote}[1]{
\ifthenelse{\boolean{shownotes}}%
{\marginpar{\raggedright\tiny\texttt{#1}}}%
{}%
}
\newcommand{\hole}[1]{
\ifthenelse{\boolean{shownotes}}%
{\begin{center} \fbox{ \rule {.25cm}{0cm} \rule[-.1cm]{0cm}{.4cm}
\parbox{.85\textwidth}{\begin{center} \texttt{#1}\end{center}} \rule
{.25cm}{0cm}}\end{center}} {} }

\graphicspath{{pics/}}
	\title{Shock-type singularity of the  hyperbolic-parabolic chemotaxis system}

\author{Woojae Lee}
\address[Woojae Lee]{\newline Department of Mathematics \newline
Yonsei University, 50 Yonsei-Ro, Seodaemun-Gu, Seoul 03722, Republic of Korea}
\email{woori0108@yonsei.ac.kr}

\numberwithin{equation}{section}

\newtheorem{theorem}{Theorem}[section]
\newtheorem{lemma}{Lemma}[section]
\newtheorem{corollary}{Corollary}[section]
\newtheorem{proposition}{Proposition}[section]
\newtheorem{remark}{Remark}[section]

\newcommand{\R}{\mathbb R}

\newcommand{\ls}{\lesssim}

\newcommand{\bq}{\begin{equation}}
\newcommand{\eq}{\end{equation}}

\makeatletter
\def\moverlay{\mathpalette\mov@rlay}
\def\mov@rlay#1#2{\leavevmode\vtop{%
   \baselineskip\z@skip \lineskiplimit-\maxdimen
   \ialign{\hfil$\m@th#1##$\hfil\cr#2\crcr}}}
\newcommand{\charfusion}[3][\mathord]{
    #1{\ifx#1\mathop\vphantom{#2}\fi
        \mathpalette\mov@rlay{#2\cr#3}
      }
    \ifx#1\mathop\expandafter\displaylimits\fi}
\makeatother

\begin{document}
\allowdisplaybreaks

\date{\today}

\subjclass[2020]{}
\keywords{Hyperbolic-parabolic chemotaxis, Finite-time blow-up, Blow-up profile, Shock-type singularity.}
	\maketitle
	\begin{abstract}
	This paper deals with the hyperbolic-parabolic chemotaxis (HPC) model, which is a hydrodynamic model describing vascular network formation at the early stage of the vasculature. We  study analytically the singularity formation associated with the shock-type structure, which was numerically observed by Filbet, Lauren{\c{c}}ot, and Perthame \cite{filbet2005derivation} and Filbet and Shu \cite{filbet2005approximation}. We construct the blow-up profile in a 1D HPC system on $\mathbb{R}$  as follows: The blow-up profile is stable in the sense of $H^m$ topology  ($m\geq 5$) prior to the occurrence of the singularity. For the first singularity, while the density and velocity $(\rho, u)$ of endothelial cells themselves remain bounded, the gradients of the density and velocity blow up. The chemoattractant concentration $\phi$ has $C^2$ regularity. However, the density and velocity  with $C^ {\frac{1}{3}}$ regularity exhibit a cusp singularity at a unique blow-up point, the location and time of which are explicitly estimated. Furthermore, the HPC system is $C^1$ differentiable except in any neighborhood of the blow-up point.
 \end{abstract}          \section{Introduction}\label{intro}
          
          Vasculogenesis, the process by which new blood vessels are formed from endothelial cells rather than from 
pre-existing vessels, is a significant biological phenomenon driven by chemotaxis, which describes the change of cells movements in response to environmental chemical substances.
 For the beginning of a vasculature, endothelial cells that are randomly distributed  can spontaneously  assemble  a vascular network, which is a pivotal factor in the tumor growth. This early stage of vasculogenesis has a difficulty to be explained by the Keller--Segel systems but is observed by the following hyperbolic-parabolic chemotaxis (HPC) system numerically.
 %         This paper deals with the hyperbolic-parabolic chemotaxis (HPC) model, which is a hydrodynamic model describing vascular network formation during the early stage of vascular development. We  study analytically the singularity formation associated with the shock-type structure, which was numerically observed by [Filbet-Lauren{\c{c}}ot-Perthame, J. Math. Biol, 50(2):189–207, 2005] and [Filbet-Shu, SIAM J. Sci. Comput, 27(3):850–872, 2005]. In the one-dimensional HPC system on the real line, we construct the blow-up profile in a HPC system on the real line  as follows: The blow-up profile is stable in Sobolev spaces prior to the occurrence of the singularity. For the first singularity, while the density and velocity of endothelial cells themselves remain bounded, the gradients of the density and velocity blow up. The chemoattractant concentration remains twice continuously differentialbe. However, the density and velocity exhibit a cusp singularity at a unique blow-up point, the location and time of which are explicitly estimated. Furthermore, the HPC system is continuously differentiable except in any neighborhood of the blow-up point.  
	  	\begin{align}\label{original}
		\partial_t\rho+\nabla\cdot (\rho u )=&0, \quad x\in\mathbb{R}^d, \quad t>0,\nonumber\\
		\partial_t(\rho u)+\nabla\cdot\left(\rho u\otimes u\right)+\nabla P(\rho)=&\mu\rho\nabla \phi-\beta\rho u,\nonumber\\
		\partial_t \phi-D\Delta \phi=&a\rho-b\phi,
	\end{align}
where $\rho(x,t)\geq 0$ and $u(x,t)\in\mathbb{R}^d$ are the density and the velocity of endothelial cells, respectively. $\phi$ indicates the concentration of the chemoattractant known as the vascular endothelial growth factor(VEGF).  The parameter $\mu > 0$ quantifies the intensity of cell response to the chemoattractant concentration gradient and the damping term $\beta \rho u$ accounts for friction force with coefficient
$\beta>0$ arising from the interaction between cells and the fixed substratum. The diffusion coefficient $D>0$ characterizes how the chemoattractant spreads. The positive constants $a$ and $b$ represent the growth and death rates of the chemoattractant, respectively. The smooth pressure function P is dependent only on the density reflecting that closely packed cells resist to compression because cellular matter cannot be penetrated.

Gamba et al. \cite{gamba2003percolation} proposed the HPC system. Filbet, Lauren{\c{c}}ot, and Perthame \cite{filbet2005derivation} provided a formal derivation of the HPC system in the case $\beta = 0$ and $P(\rho) =A \rho$ with $A>0$, starting from a kinetic transport framework modeling the run-and-tumble behavior of self-propelled particles.  Chavanis and Sire \cite{chavanis2007kinetic} formally derived the HPC system from a nonlinear mean-field Fokker-Planck equation by invoking a local thermodynamic equilibrium assumption.

Both the Keller Segel system and the HPC system are liable to exhibit singularity formation. However,
it is noteworthy that, in distinction to the Keller-Segel system, which typically results in the formation of Dirac peaks, the blow-up profile of the HPC system appears to align with the aggregation patterns of endothelial cells observed in experiments \cite{gamba2003percolation, serini2003modeling}.  In the case of the Keller-Segel model, the density of cells is observed to accumulate in the region surrounding an isolated point. In contrast, the profile  numerically illustrated by Filbet, and Shu \cite{filbet2005approximation} and Filbet, Lauren{\c{c}}ot, and Perthame  \cite{filbet2005derivation} in the HPC system is considered to show that shock-type structures emerge prior to the implosion of density. It implies that the density accumulates in the vicinity of the edges of a network.

The main objective of this paper is to investigate the singularity formation of shock-type structures in the following 1D HPC system on $\mathbb{R}$:
\begin{align}\label{hpc}
	\rho_t+(\rho u )_x&=0, \quad x\in\mathbb{R}, \quad t>0,\nonumber \\
	(\rho u)_t+(\rho u^2)_x +A (\rho^\gamma)_x&=\mu \rho\phi_x-\beta\rho u,\nonumber \\
	\phi_t-D \phi_{xx}&=a\rho-b\phi.
\end{align}
where the pressure is given by
\begin{equation*}
P(\rho)=A\rho^\gamma, \quad A>0, \ \gamma>1,
\end{equation*}
corresponding to the isentropic case in fluid dynamics. This formulation was introduced to describe in vitro experiments by Serini et al. \cite{serini2003modeling}, capturing the volume-filling effect. Numerical simulations have been conducted in \cite{natalini2014well, natalini2015numerical, di2011analysis}. Moreover, stationary solutions with vacuum in 1D were constructed in \cite{berthelin2016stationary}. 

 Next, we study the scaling properties of the HPC system by introducing the characteristic scales $L, T, R, U,$ and $\Phi$, and defining the corresponding dimensionless variables by
\begin{equation*}
x=L \tilde{x}, \quad t=T \tilde{t}, \quad \rho=R \tilde{\rho}, \quad u=U \tilde{u}, \quad \phi=\Phi \tilde{\phi}.
\end{equation*}
To preserve the form of the continuity equation, we set
\begin{equation*}
U=\frac{L}{T}.
\end{equation*}
Substituting these scalings into the HPC system \eqref{hpc} yields
\begin{align}\label{hpc1}
&\tilde{\rho}_{\tilde{t}}+(\tilde{\rho}\tilde{u})_{\tilde{x}}=0, \nonumber\\
&\tilde{u}_{\tilde{t}}+\tilde{u} \tilde{u}_{\tilde{x}}+A R^{\gamma-1} \frac{T^2}{L^2} \, \gamma \, \tilde{\rho}^{\gamma-2} \tilde{\rho}_{\tilde{x}}=\frac{\mu \Phi T^2}{L^2} \, \tilde\phi_{\tilde{x}}-\beta T \, \tilde{u}, \nonumber\\
&\tilde{\phi}_{\tilde{t}}-\frac{D T}{L^{2}} \tilde{\phi}_{\tilde{x} \tilde{x}}+b T \, \tilde{\phi}=\frac{a R T}{\Phi} \, \tilde{\rho}.
\end{align}
To fix the coefficients of the terms $\tilde{\phi}_{\tilde{x} \tilde{x}}$  and $\tilde{\phi}$, we choose the scales so that
\begin{equation*}
\frac{DT}{L^2}=1, \qquad bT=1,
\end{equation*}
which yields
\begin{equation*}
T=\frac{1}{b},  \qquad L=\sqrt{\frac{D}{b}}.
\end{equation*}
Next, we fix the coefficients of the pressure term and the term $\tilde{\phi}_{\tilde{x}}$ by requiring
\begin{equation*}
A R^{\gamma-1} \frac{T^{2}}{L^2} \gamma=1, \qquad \frac{\mu \Phi T^{2}}{L^{2}}=1.
\end{equation*}
With these choices, the system \eqref{hpc1} reduces to 
\begin{align}\label{hpc2}
&\rho_t+(\rho u)_x=0,\nonumber \\
&u_t+u \, u_x+ \rho^{\gamma-2} \rho_x=\phi_x-\frac{\beta}{b} \, u,\nonumber \\
&\phi_t-\phi_{xx}+ \phi=a\frac{\mu}{b^2D}\left(\frac{b D}{\gamma A}\right)^{\frac{1}{\gamma-1}} \, \rho.
\end{align}
Here, we have dropped the tildes and denote the dimensionless variables by
\begin{equation*}
x=\tilde{x},\qquad
t=\tilde{t},\qquad
\rho=\tilde{\rho},\qquad
u=\tilde{u},\qquad
\phi=\tilde{\phi}.
\end{equation*} For simplicity, we select the remaining parameters as
\begin{equation*}
a\frac{\mu}{D}\left(\frac{D}{\gamma A}\right)^{\frac{1}{\gamma-1}}=b=1.
\end{equation*}
By introducing the variables $q=\frac{1}{\alpha}\rho^\alpha$ where $\alpha=\frac{\gamma-1}{2}$, the HPC system \eqref{hpc2}  is transformed into  
\begin{align}\label{simple model}
		q_t+u q_x +\alpha q  u_x=&0,\nonumber \\
		u_t+u u_x +\alpha q  q_x=& \phi_x-\beta u,\nonumber \\
		\phi_t-\phi_{xx}=&\alpha^\frac{1}{\alpha} q^\frac{1}{\alpha} -\phi.
\end{align}

 We find a blow-up profile which forms a shock-type singularity. Roughly speaking, we prove the following. 
\begin{theorem}\label{intro thm}
For some smooth initial data $(q_0,u_0,\phi_0)$ with $\mathcal{O}(1)$ in $L^\infty(\mathbb{R})$ and
\begin{equation*}
\min_{x\in\mathbb{R}}(q_0+u_0)_x(x)=-\frac{1}{\epsilon},
\end{equation*}
for $\epsilon>0$ taken sufficiently small, 
there exists a smooth solution  of the HPC system (\ref{simple model}) such that the following hold true:
 \begin{enumerate}[(1)]
 	\item The blow-up profile is stable in the sense of $H^m(\mathbb{R})$ topology for any $m\geq 5$ prior to the occurrence of the singularity.
 	\item The solution $(q,u)$ has a unique blow-up point $(x^*,T^*)$, where $|x^*|=\mathcal{O}(\epsilon)$ and $T^*= \mathcal{O}(\epsilon)$:
 	\begin{equation*}
 		 q_x(x^*,T^*)=u_x(x^*,T^*)=-\infty.
 	\end{equation*}
 	\item At the first blow-up point, the solution $(q,u)$ has a cusp singularity with H\"older $C^\frac{1}{3}$ regularity.
  \item  The solution $(q,u)$ remains $C^1$ regular away from the blow-up point $(x^*,T^*)$ for $0\leq t\leq T^*$.  In contrast, the solution $\phi$ remains $C^2$ regular up to time $T^*$.
 \end{enumerate}
\end{theorem} 

\begin{remark}
 	This blow-up profile suggests that a sufficient accumulation of the density gradient can lead to blow-up. Moreover, the behavior of $q$ and $u$ near the blow-up point is characterized by $-(x-x^*)^{\frac{1}{3}}$ in a similar manner. See Remark \ref{structure}.
\end{remark}

\begin{remark}
	Theorem \ref{intro thm} also holds for $\beta=0$.  For $\beta>0$,  the absolute value of the maximally negative slope increases as $\beta$ increases. For further details, see Remark \ref{rmk;relation}. 
\end{remark}

\begin{remark}
 The precise assumption of initial data and the detailed statement of the Theorem \ref{intro thm} are given in Section \ref{sec;main result}.
\end{remark}

 There is literature concerning the global existence and asymptotic stability of classical solutions for  the  system \eqref{original}. First, Russo and Sepe \cite{russo2013existence} proved the global-in-time existence and asymptotic stability for solutions near a constant ground state $(\bar\rho, 0 , \frac{a}{b}\bar\rho)$ in the  Sobolev space $H^m$($m\geq d/2+1$) where the constant ground state is also small. Liu et al. \cite{liu2022asymptotic} demonstrated the convergence of linear diffusion waves in 3D. Crin-Barat et al \cite{crin2023hyperbolic} showed the global-in-time existence and temporal decay estimates for solutions of (HPC) near a constant ground state $(\bar\rho, 0, \frac{a}{b}\bar\rho)$ in the  hybrid Besov space $\dot{B}^\frac{d}{2}_{2,1}\cap\dot{B}^{\frac{d}{2}+1}_{2,1}$. Under the assumption $bP'(\rho)-a\mu\rho>0$, Hong et al. \cite{hong2021nonlinear} established the global existence and asymptotic stability of steady-states and  Liu et al. \cite{liu2022convergence} proved convergence of nonlinear diffusion waves in 1D.
 
To the best of our knowledge, this is the first analytical blow-up analysis  for the HPC system \eqref{original}. We  not only demonstrate  $H^m$ norm inflation, but also provide a detailed description of the scenario for the first singularity. In order to provide this blow-up profile, we adapt the approach proposed by Buckmaster, Vicol, and Shkoller in \cite{buckmaster2023formation, buckmaster2022formation}  for the construction of shock solutions to the multidimensional isentropic compressible Euler equations.

 Before the emergence of this approach, the theory of shock formation for the compressible Euler equations had already been extensively developed. Riemann \cite{riemann1860fortpflanzung} first proved shock formation for the 1D isentropic Euler equations by introducing Riemann invariants and employing the method of characteristics. In multiple dimensions, Sideris \cite{sideris1985formation} proved $C^1$ norm inflation in Sobolev spaces for the 3D case via a moment method combined with a contradiction argument. For 3D irrotational flows, Christodoulou \cite{christodoulou2007formation} first established shock formation for the relativistic Euler equations by introducing the inverse foliation density, whose vanishing signals shock formation.  Christodoulou and Miao \cite{chrismiao} subsequently extended a corresponding framework for the non-relativistic case. In the isentropic case with nonzero vorticity, Luk and Speck \cite{luk2016shock} gave the first proof of the 2D shock formation by incorporating the ``vorticity'' transport equation into the framework in \cite{christodoulou2007formation, chrismiao}.

 In contrast to these earlier multidimensional results, Buckmaster, Vicol, and Shkoller obtain a sharp characterization of the singularity for shock solutions with nontrivial vorticity, including a uniquely determined blow-up time and location with explicit estimates, as well as the formation of a cusp singularity. They are inspired from the self-similar modulation analysis, which was developed to study the singularity formation for Schr\"odinger equations and nonlinear heat equations \cite{merle1996asymptotics, merle2005blow, merle2020strongly, merle1997stability}.  In their modulation analysis, they analyze the stability of a self-similar variable $W$ near the steady-state Burgers profile $\overline{W}$ in Section \ref{sec;burger}, where this solution is useful to study the self-similar blow-up for Burgers equations (see \cite{cassel1996onset, collot2018singularity, eggers2008role}).  Through this process, they prove the shock formation to the 2D compressible Euler equations with azimuthal symmetry \cite{buckmaster2022formation} and point shock formation to the 3D compressible Euler equations \cite{buckmaster2023formation}. Asymptotically self-similar Burgers-type profiles have been employed to investigate the singularity formation of both hyperbolic and non-hyperbolic PDEs, including the non-isentropic Euler equations, the Euler-Poisson system in ion dynamics, and the Burgers-Hilbert equation etc in \cite{bae2024structure, buckmaster2023shock, oh2024gradient, yang2021shock}.

 \subsection{Outline of the Strategy.} The initial step is to reformulate the HPC system \eqref{simple model} by introducing Riemann-type variables $(w,z)$ in \eqref{riemann;type}. Subsequently, self-similar transformation $\left(x,\tfrac{2}{1+\alpha} t\right)\mapsto(y,s)$ and modulation variables ($\tau,\kappa,\xi$) in \eqref{self trans} and \eqref{self similar} are employed and the variables $(w,z,\phi,\rho,u)$ correspond to the self-similar variables $(e^{-\frac{s}{2}}W+\kappa, Z,\Phi, \sigma, U)$ in \eqref{self similar} and \eqref{dense,vel} to transform the finite-time blow-up problem into a global stability analysis. We remark that Russo and Sepe \cite{russo2013existence} have already demonstrated the local well-posedness in Sobolev spaces. However, they transformed the HPC system into a more tractable structure by proposing entropy variables. In order to guarantee the $H^m(\mathbb{R})$-stability before the singularity formation, it is necessary for us to prove the local well-posedness and derive a continuation criterion \eqref{continucri} for the system \eqref{simple model} in Theorem \ref{simple result}, where the criterion is related to the accumulation of the gradients of $q$ and $u$.
 
   We define the ODE system for modulation variables by constraining some spatial derivatives of a self-similar variable $W$ at 0 in Section \ref{sec;formmodul}. We assume the initial data in Section \ref{assump ini} with the corresponding to bootstrap assumptions in Section \ref{boot;assum}. In the bootstrap argument with bootstrap assumptions, the modulation variables $(\tau,\xi)$ estimates are closed and the self-similar variable $\Phi$ estimates is established using classical heat kernel estimates. To close the self-similar variables ($W,Z$) estimates, we derive the several bounds of particle trajectories for the velocities of self-similar variables $(W,Z)$ in Section \ref{subsec;traj}. Combining these bounds and bootstrap assumptions, we partially close the $(W, Z)$ estimates in Section \ref{Sec z} and \ref{sec w}. 
  
To get information such as the uniqueness of the blow-up point, the presence of a cusp singularity and their regularity, we carefully analyze  the $L^\infty$-stability of the spatial gradient of the self-similar variable $W$ near the spatial gradient of the steady-state Burgers profile $\overline{W}$ which satisfies the following equation:
 \begin{equation*}
	-\frac{1}{2}\overline{W}+\left(\frac{3}{2}y+\overline{W}\right)\overline{W}_y=0, \quad \hbox{for all} \quad y\in\mathbb{R}. 
\end{equation*}
 First, we get the $L^\infty$-stability of $W_y$ near $\overline{W}_y$ near $0$ by standard bootstrap method in Section \ref{sec;near 0}. Next, we investigate the asymptotic behavior of $W_y$ away from $0$. Since the parabolic component $\phi$ in the system \eqref{simple model} does not exhibit finite speed of propagation, we establish the local well-posedness of the system \eqref{simple result} for weighted variables in Sobolev spaces and derive a corresponding continuity criterion in Theorem \ref{weight local}, analogous to that in Theorem \ref{simple result}. This allows us to analyze not only uniform spatial decay rates of the perturbations in $q$, $u$, and the spatial derivatives of $\phi$ in Section \ref{sec;uniform}, but also weighted estimates for the self-similar variables $(W_y, Z_y)$. In addition, owing to the damping effect of $\beta u$, which prevents singularity formation in the system, we perform a bootstrap argument to obtain a weighted estimate for $Z_y$. More specifically, for the self-similar variable $W$, we establish the $L^\infty$-stability of $W_y$ near the steady-state solution $\overline{W}_y$ on the middle interval and the asymptotic behavior of $W_y$ as $y\to\infty$. By closing the bootstrap argument, we obtain pointwise estimates for $w$, $z$, and $\phi$ at the blow-up time $T^*$, providing insights into the singularity scenario.
  
%  Our analysis takes into account a structural aspect of the HPC system namely, that the interaction between its hyperbolic and parabolic components does not satisfy the finite speed of propagation.
 \subsubsection*{Notations} All generic positive constants are denoted by $C$ which is independent of $t$. $C(p_1, p_2, \dots )$ implies that a constant depends on $p_1, p_2, \dots$.  The notations $f \ls g$ and $f\thicksim g$ mean that there exists a positive constant $C$ such that $f \leq C g$ and $C^{-1}f\leq g\leq C f$, respectively. The notation $f\lesssim_{p_1,p_2,\cdots}g $ is used to indicate that $f\leq  C g$ for some $C=C(p_1,p_2,\cdots)>0$. For simplicity, we denote $\int_{\mathbb{R}}\,dx$, $L^p(\R)$, and $H^m(\R)$ as  $\int \,dx$, $L^p$, and $H^m$, respectively.  We write $\|g(u_1,\cdots, u_n)\|^2_X$ to mean $\|gu_1\|^2_X+\cdots+\|gu_n\|^2_X$, where $\|\cdot\|_X$ is the norm of the function space $X$.  $\partial^k_x$ represents a derivative of order $k$ with respect to a variable $x$. We sometimes omit the time dependence of differential operators, for example: $\dot{f}=\frac{df}{dt}$.

\subsection*{Organization of this paper}  The rest of this paper is organized as follows. In Section \ref{sec;prelimi}, we present  the self-similar Burgers profile and reformulate the HPC system.  In Section \ref{sec;main result}, we state the assumptions of initial data and provide the detailed statement of the main theorem. In Section \ref{sec;localconti}, we discuss the proof of the local well-posedness and its associated continuation criterion. In Section \ref{sec;formmodul}, we define the modulation variables. In Section \ref{boot;assum},  we set the bootstrap assumptions. Section \ref{sec;closboot} and \ref{sec;stabilW} are devoted to  closing the bootstrap argument to analyze the global stability of modulation variables and self-similar variables. Finally, in Section \ref{sec;finalmain}, we prove the main theorem. 

\section{Preliminaries}\label{sec;prelimi}
\subsection{Stable self-similar Burgers profile.}\label{sec;burger}
The steady-state and self-similar Burgers profile to study singularity formation of the Burgers equation \cite{cassel1996onset,collot2018singularity,eggers2008role} is as follows:
 \begin{equation}\label{burgers}
	\overline{W}(y)=\left(-\frac{y}{2}+\left(\frac{1}{27}+\frac{y^2}{4}\right)^\frac{1}{2}\right)^\frac{1}{3}-\left(\frac{y}{2}+\left(\frac{1}{27}+\frac{y^2}{4}\right)^\frac{1}{2}\right)^\frac{1}{3},
\end{equation}
which solves
\begin{equation}\label{eq;burger}
	-\frac{1}{2}\overline{W}+\left(\frac{3}{2}y+\overline{W}\right)\partial_y\overline{W}=0.
\end{equation}
We present quantitative properties of the following $\overline{W}$, which are  useful for the bootstrap argument in Section \ref{sec;closboot} and \ref{sec;stabilW}. The proofs can be found in \cite{buckmaster2022formation, yang2021shock}, where they are derived from the explicit solution \eqref{burgers}. 
\begin{lemma}\label{lem;sbur} The profile $\overline{W}$  satisfies the following properties$:$
\begin{enumerate}[(1)]
\item The values of the derivatives of $\overline{W}$ at $y = 0$ are given by 
\begin{equation*} \overline{W}(0) = 0, \qquad \partial_y \overline{W}(0) = -1, \qquad \partial^2_y \overline{W}(0) = 0, \qquad \text{and} \qquad \partial^3_y \overline{W}(0) = 6. 	
\end{equation*}

\item For all  $y \in \mathbb{R}$,
\begin{equation*} 
   \ \ \ \ |\overline{W}(y)| \leq (1 + y^2)^{\frac{1}{6}}, \quad |\partial_y \overline{W}(y)| \leq (1 + y^2)^{-\frac{1}{3}}, \quad \text{and} \quad |\partial^2_y \overline{W}(y)| \leq (1 + y^2)^{-\frac{5}{6}}.
\end{equation*}

\item For all $|y| \leq \frac{1}{5}$,
\begin{equation*}
 |\partial^3_y \overline{W}(y)| \leq 6, \quad |\partial^4_y \overline{W}(y)| \leq 30, \quad \text{and} \quad |\partial^5_y \overline{W}(y)| \leq 360.
\end{equation*}
\item For all  $|y| \geq 100$, 
\begin{equation*}
 -\frac{7}{20} y^{-\frac{2}{3}}  \leq \partial_y\overline{W}(y) \leq -\frac{1}{4}  y^{-\frac{2}{3}}.
 \end{equation*}
\item For all  $y \in \mathbb{R}$,
\begin{equation*}
 \frac{5}{2} + 3 \partial_y \overline{W} + \frac{1}{y(1 + y^2)} \left(\frac{3y}{2} + \overline{W}\right) \geq \frac{y^2}{1 + y^2}. 
 \end{equation*}
\end{enumerate}
\end{lemma}

\subsection{Reformulation.}\label{sec;reform} To study the blow-up profile of the HPC system \eqref{simple model}, we reformulate it using Riemann-type and modulation variables. Since the HPC system \eqref{simple model} is not Galilean invariant---namely, it is not preserved under the transformation \(x \mapsto x - vt\), \(u \mapsto u + v\), which leaves the compressible Euler equations invariant---we introduce the following Riemann-type variables.
%Since the HPC system \eqref{simple model} lacks Galilean invariance, unlike the compressible Euler equations, we introduce the following Riemann-type variables:
% In contrast to the compressible Euler equations, the HPC system \eqref{simple model} lacks Galilean invariance, which motivates the introduction of the following Riemann-type variables:
% Since in contrast to compressible Euler equations, the HPC system \eqref{simple model} does not possess Galilean invariance, we introduce the following Riemann-type variables: 
\begin{equation}\label{riemann;type}
	w=u+ q+\frac{\kappa_0}{2},\ \ \text{and} \ \ z=u-q+\frac{\kappa_0}{2},
\end{equation}
for some constant $\kappa_0>0$, which represents the initial data of the modulation variable $\kappa$. Then we can rewrite the system (\ref{simple model}) as
\begin{equation*}
		\begin{split}
		\partial_t w+\frac{\beta}{2}w+\left(\frac{1+\alpha}{2}w+\frac{1-\alpha}{2}z-\frac{\kappa_0}{2} \right) w_x =&\phi_x+\frac{\beta\kappa_0}{2}-\frac{\beta}{2}z,\\
		\partial_t z+\frac{\beta}{2}z+\left(\frac{1-\alpha}{2}w+\frac{1+\alpha}{2}z-\frac{\kappa_0}{2}\right) z_x =& \phi_x+\frac{\beta\kappa_0}{2}-\frac{\beta}{2}w,\\
		\partial_t \phi- \phi_{xx}+\phi=&\alpha^{\frac{1}{\alpha}}q^{\frac{1}{\alpha}}.
		\end{split}
\end{equation*}
By defining  $t=\frac{2}{1+\alpha}\tilde{t}$, and then setting $t=\tilde{t}$, the preceding system becomes 
\begin{subequations} \label{riemann}
\begin{align}
\label{subrie1}		\partial_{t} w+\frac{\beta}{1+\alpha}w+\left(w+\frac{1-\alpha}{1+\alpha}z-\frac{\kappa_0}{1+\alpha} \right) w_x =&\frac{2}{1+\alpha}\left(\phi_x+\frac{\beta\kappa_0}{2}-\frac{\beta}{2}z\right),\\
\label{subrie2}		\partial_{t} z +\frac{\beta}{1+\alpha}z+\left(\frac{1-\alpha}{1+\alpha}w+z-\frac{\kappa_0}{1+\alpha}\right) z_x=&\frac{2}{1+\alpha}\left( \phi_x+\frac{\beta\kappa_0}{2}-\frac{\beta}{2}w\right), \\		
\label{subrie3}\partial_{t} \phi-\frac{2}{1+\alpha} \phi_{xx}+\frac{2}{1+\alpha}\phi=&\frac{2\alpha^\frac{1}{\alpha} }{1+\alpha}q^\frac{1}{\alpha}. 
\end{align}
\end{subequations}
We employ the following self-similar transformation  
\begin{align}\label{self trans}
	y(x,t)&=\frac{x-\xi(t)}{(\tau(t)-t)^\frac{3}{2}}, \qquad  s(t)=-\ln(\tau(t)-t),
\end{align}
and define self-similar variables $(W,Z,\Phi)$ by
\begin{align}\label{self similar}
	w(x,t)&=e^{-\frac{s}{2}}W(y,s)+\kappa(t), \qquad z(x,t)=Z(y,s), \qquad \phi(x,t)=\Phi(y,s),
\end{align}
where the modulation variables $\xi(t)$, $\tau(t)$, and $\kappa(t)$ are determined by an ODE system introduced in Section~\ref{sec;formmodul}. Here, $\xi$, $\tau$, and $\kappa$ represent the blow-up location, the blow-up time, and the wave amplitude of the cell, respectively. We will establish that $\tau(t) > t$ for all $t \in [0, T^*)$, and the blow-up time is given by $T^*=\tau(T^*)$.  The following identities are used to obtain the self-similar form of the system:
%where the modulation variables  $\xi(t) ,  \tau(t)$, and $\kappa(t)$ are determined by an ODE system with the initial data $(\xi(0), \tau(0),\kappa(0))=(0,\epsilon,\kappa_0)$ introduced in Section \ref{sec;formmodul}. Here, $\kappa$, $\xi$, and  $\tau$ represent wave amplitude of cells, the blow-up location, and the blow-up time, respectively. In particular, $\tau(t)>t$ for all $t\in[0,T^*)$ and the blow-up time is given by $T^*=\tau(T^*)$ in our analysis.
%The following identities are used to obtain the self-similar form of the system: $t\in[0,T^*)$ maps to $s\in[-\ln\epsilon, \infty)$.
\begin{align}\label{identity}
		&e^{s}=\frac{1}{\tau(t)-t},\qquad \frac{ds}{dt}=e^s(1-\dot{\tau}),\nonumber\\
	y=e^{\frac{3}{2}s}(x-\xi(t)), &\qquad  \partial_x y=e^{\frac{3}{2}s},\qquad \hbox{and}  \qquad \partial_t y=-e^{\frac{3}{2}s}\dot{\xi}+\frac{3}{2}e^{s}(1-\dot{\tau})y.
\end{align}
Applying these identities \eqref{identity} to the equations \eqref{subrie1} and \eqref{subrie2} yields the following system with respect to ($W,Z,\Phi$):
\begin{align} \label{no ordero}
	    &\left(\partial_s-\frac{1}{2}+\frac{\beta e^{-s}}{(1-\dot{\tau})(1+\alpha)}\right)W \nonumber\\
	    &+\frac{1}{1-\dot{\tau}} \left\{e^{\frac{s}{2}}\left(\kappa+\frac{1-\alpha}{1+\alpha}Z-\dot{\xi}-\frac{\kappa_0}{1+\alpha}\right)
		+\frac{3}{2}(1-\dot{\tau})y+W\right\}W_y\nonumber\\
		&=\frac{e^{-\frac{s}{2}}}{1-\dot{\tau}}\left(\frac{2}{1+\alpha} \left(e^{\frac{3}{2}s} \partial_y\Phi+\frac{\beta\kappa_0}{2}-\frac{\beta}{2}Z\right)-\dot{\kappa}-\frac{\beta\kappa}{1+\alpha}\right),\nonumber\\
		&\left(\partial_s+\frac{\beta e^{-s}}{(1-\dot{\tau})(1+\alpha)}\right)Z\nonumber\\
		&+\frac{1}{1-\dot{\tau}}\left\{e^{\frac{s}{2}}\left(\frac{1-\alpha}{1+\alpha}\kappa -\dot{\xi}-\frac{\kappa_0}{1+\alpha} \right)+\frac{1-\alpha}{1+\alpha}W
		+\frac{3}{2}(1-\dot{\tau})y+e^\frac{s}{2} Z\right\}Z_y\nonumber\\
		&=\frac{2e^{-s}}{(1-\dot{\tau})(1+\alpha)}\left(e^{\frac{3}{2}s} \partial_y\Phi+\frac{\beta \kappa_0}{2}-\beta\frac{e^{-\frac{s}{2}}W+\kappa}{2}\right).
\end{align}
We rewrite the system \eqref{no ordero} as
\begin{subequations} \label{no order}
	\begin{align}
	    &\left(\partial_s-\frac{1}{2}+\frac{\beta e^{-s}}{(1-\dot{\tau})(1+\alpha)}\right)W+ \left(G_W
		+\frac{3}{2}y+\frac{W}{1-\dot{\tau}}\right)W_y=F_W,\label{no order W}\\
		&\left(\partial_s+\frac{\beta e^{-s}}{(1-\dot{\tau})(1+\alpha)}\right)Z+\left(G_Z+\frac{3}{2}y+\frac{ e^\frac{s}{2}}{1-\dot{\tau}} Z\right)Z_y=F_Z,\label{no order Z}
	\end{align}
\end{subequations}
where transport components and forcing terms are defined as follows:
\begin{align}\label{g and force} 
	&G_W=\frac{e^{\frac{s}{2}}}{1-\dot{\tau}}\left(\kappa+\frac{1-\alpha}{1+\alpha}Z-\dot{\xi}-\frac{\kappa_0}{1+\alpha} \right),\nonumber\\
	&G_Z=\frac{e^{\frac{s}{2}}}{1-\dot{\tau}}\left(\frac{1-\alpha}{1+\alpha}\kappa-\dot{\xi}-\frac{\kappa_0}{1+\alpha} \right)+\frac{1-\alpha}{(1-\dot\tau)(1+\alpha)}W,\nonumber\\
 &F_W=\frac{e^{-\frac{s}{2}}}{1-\dot{\tau}}\left(\frac{2}{1+\alpha} \left(e^{\frac{3}{2}s} \partial_y\Phi+\frac{\beta\kappa_0}{2}-\frac{\beta}{2}Z\right)-\dot{\kappa}-\frac{\beta\kappa}{1+\alpha}\right),\nonumber\\
   &F_Z=\frac{2e^{-s}}{(1-\dot{\tau})(1+\alpha)}\left(e^{\frac{3}{2}s} \partial_y\Phi+\frac{\beta \kappa_0}{2}-\beta\frac{e^{-\frac{s}{2}}W+\kappa}{2}\right).
\end{align}
For all $n\in\mathbb{N}$, applying the differential operator $\partial^n_y$ and the Leibniz rule to the system \eqref{no order} yields
\begin{align}\label{higher order}
			&\left(\partial_s+\frac{\beta e^{-s}}{(1-\dot{\tau})(1+\alpha)}+\frac{3n-1}{2}+\frac{n+\chi_{n\geq 2}}{1-\dot{\tau}}\partial_yW+n\partial_yG_W\right)\partial^n_yW\nonumber\\
			&+\left(G_W+\frac{3}{2}y+\frac{1}{1-\dot{\tau}}W\right)\partial^{n+1}_yW\nonumber\\
			&=\partial^n_yF_W-\chi_{n\geq 2}\partial^n_y G_W\partial_yW-\chi_{n\geq 3}\sum^{n-1}_{k=2}\begin{pmatrix}
			n\\
			k
		\end{pmatrix}\left(\frac{1}{1-\dot{\tau}}\partial^k_yW +\partial^k_y G_W\right)\partial^{n-k+1}_yW,\nonumber\\
			&\left(\partial_s+\frac{\beta e^{-s}}{(1-\dot{\tau})(1+\alpha)}+\frac{3n}{2}+\frac{n+\chi_{n\geq 2}}{1-\dot{\tau}}e^\frac{s}{2}\partial_yZ +n\partial_y G_Z \right)\partial^n_yZ\nonumber\\
			&+\left(G_Z+\frac{3}{2}y+\frac{e^\frac{s}{2}}{1-\dot{\tau}}Z\right)\partial^{n+1}_yZ\nonumber\\
		&=\partial^n_yF_Z-\chi_{n\geq 2}\partial^n_y G_Z\partial_yZ-\chi_{n\geq 3}\sum^{n-1}_{k=2}\begin{pmatrix}
			n\\
			k
		\end{pmatrix}\left(\frac{e^\frac{s}{2}}{1-\dot{\tau}}\partial^k_y Z+\partial^k_y G_Z \right)\partial^{n-k+1}_yZ,	
\end{align}
where $\chi$ is the characteristic function. For the sake of simplicity, we rewrite the systems \eqref{no order} and \eqref{higher order} as
\begin{subequations} 
	\begin{align}
	    \partial_s\partial^n_yW+D^{(n)}_W\partial^n_yW+\mathcal{V}_W\partial^{n+1}_y W=& F^{(n)}_W,\label{higher order W}\\
		\partial_s\partial^n_yZ+D^{(n)}_Z\partial^n_yZ+\mathcal{V}_Z\partial^{n+1}_y Z=& F^{(n)}_Z, \qquad \hbox{for} \ n=0,1,2,\cdots \label{higher order Z}
	\end{align}
\end{subequations}
where damping terms, transport velocity terms, and forcing terms are given by: 
\begin{align}\label{damp,vel}
	D^{(n)}_W=&\frac{\beta e^{-s}}{(1-\dot{\tau})(1+\alpha)}+\frac{3n-1}{2}+\frac{n+\chi_{n\geq 2}}{1-\dot{\tau}}\partial_yW+n\partial_yG_W,\nonumber\\
		D^{(n)}_Z=&\frac{\beta e^{-s}}{(1-\dot{\tau})(1+\alpha)}+\frac{3n}{2}+\frac{n+\chi_{n\geq 2}}{1-\dot{\tau}}e^\frac{s}{2}\partial_yZ +n\partial_y G_Z,\nonumber\\
	\mathcal{V}_W=&G_W+\frac{3}{2}y+\frac{1}{1-\dot\tau}W, \nonumber\\
	\mathcal{V}_Z=&G_Z+\frac{3}{2}y+\frac{e^\frac{s}{2}}{1-\dot\tau}Z, \nonumber\\
			F^{(n)}_W=&\partial^n_yF_W-\chi_{n\geq 2}\partial^n_y G_W\partial_yW-\chi_{n\geq 3}\sum^{n-1}_{k=2}\begin{pmatrix}
			n\\
			k
		\end{pmatrix}\left(\frac{1}{1-\dot{\tau}}\partial^k_yW +\partial^k_y G_W\right)\partial^{n-k+1}_yW,\nonumber\\
	    F^{(n)}_Z=&\partial^n_yF_Z-\chi_{n\geq 2}\partial^n_y G_Z\partial_yZ-\chi_{n\geq 3}\sum^{n-1}_{k=2}\begin{pmatrix}
			n\\
			k
		\end{pmatrix}\left(\frac{e^\frac{s}{2}}{1-\dot{\tau}}\partial^k_y Z+\partial^k_y G_Z \right)\partial^{n-k+1}_yZ. 
\end{align}
 In Section \ref{sec;stabilW}, we analyze the $L^\infty$-stability of $\partial_yW$ near $\partial_y\overline{W}$. To facilitate this, we define a perturbation as
 \begin{equation*}
 	\widetilde{W}=W-\overline{W}.
 \end{equation*}
By subtracting the $\overline{W}$ equation  \eqref{eq;burger} from the  equation \eqref{no order W}, $\widetilde{W}$ satisfies the following equation: 
\begin{equation*}
	\begin{split}
		 &\left(\partial_s-\frac{1}{2}+\frac{\beta e^{-s}}{(1-\dot{\tau})(1+\alpha)}+\frac{\partial_y\overline{W}}{1-\dot\tau}\right)\widetilde{ W}+ \left(G_W
		+\frac{3}{2}y+\frac{1}{1-\dot{\tau}} W\right)\partial_y\widetilde{ W}\\
		&=F_W-\left(G_W+\frac{\dot\tau}{1-\dot\tau} \overline{W}\right)\partial_y\overline{W}-\frac{\beta e^{-s}}{(1-\dot\tau)(1+\alpha)}\overline{W}.    
			\end{split}
\end{equation*}
In a similar manner to obtain the system \eqref{higher order}, taking the derivative $\partial^n_y$ $(n=1,2,3,\cdots)$ of the above equation yields that the function $\partial_y^n\widetilde{W}$ $(n=0,1,2,\cdots)$ obeys

\begin{align}
	    &\partial_s\partial^n_y\widetilde{W}+D^{(n)}_{\widetilde{W}}\partial^n_y\widetilde{W}+\mathcal{V}_W\partial^{n+1}_y\widetilde{ W}=F^{(n)}_{\widetilde{W}},\label{higher order pertW}
	\end{align}
where damping terms and forcing terms are given by 
\begin{align}\label{damp vel pert}
D^{(n)}_{\widetilde{W}}=& \frac{3n-1}{2}+\frac{\beta e^{-s}}{(1-\dot{\tau})(1+\alpha)}+n\partial_yG_W+\frac{n\partial_y\widetilde{W}+\chi_{n\geq2}\partial_y\widetilde{W}+(n+1)\partial_y\overline W}{1-\dot\tau},\nonumber\\
F^{(n)}_{\widetilde{W}}=& \partial^n_y F_{W}-\frac{\beta e^{-s}}{(1-\dot\tau)(1+\alpha)}\partial^n_y\overline{W} -G_W\partial^{n+1}_y\overline{W}-\frac{\dot\tau}{1-\dot{\tau}}\partial_y\overline{W}\partial^n_y\overline{W}\nonumber\\
& -\chi_{n\geq1}\sum^n_{k=1}\begin{pmatrix}
		 	n\\
		 	k
		 \end{pmatrix} \partial^{k}_y G_W\partial^{n-k+1}_y\overline{W} 
		 -\chi_{n\geq 1}\sum^n_{k=1}\begin{pmatrix}
		 	n\\
		 	k
		 \end{pmatrix}\left(\frac{\dot{\tau}\partial^{n-k}_y\overline{ W}+ \partial^{n-k}_y\widetilde{ W} }{1-\dot\tau}\right) \partial^{k+1}_y\overline{W}\nonumber\\
		 &-\chi_{n\geq 2} \left(\partial^{n}_y G_W+\frac{\partial^n_y\overline{W}}{1-\dot{\tau}} \right)\partial_y\widetilde{W}\nonumber\\
		 & -\chi_{n\geq 3}\sum^{n-1}_{k=2}\begin{pmatrix}
		 	n\\
		 	k
		 \end{pmatrix} \left(\partial^{k}_y G_W+\frac{\partial^k_y\widetilde{W}+\partial^k_y\overline{W}}{1-\dot{\tau}} \right)\partial^{n-k+1}_y\widetilde{W}. 	 
\end{align}	
Finally, we derive the system for the self-similar variables for the density and velocity $(\rho,u)$ to establish the uniform spatial decay rates near constant states. 
To begin, we recall the equations for $\rho$ and $u$ from the systems \eqref{hpc2} and \eqref{simple model},  respectively, written in the rescaled time variable introduced before \eqref{riemann}:
\begin{align}\label{eq;rhou}
		 \rho_t+\frac{2}{1+\alpha} u\rho_x+\frac{2}{1+\alpha}\rho u_x=&0,\nonumber\\
		u_t+\frac{2}{1+\alpha}u u_x +\frac{2\alpha}{1+\alpha} q q_x +\frac{2\beta}{1+\alpha} u=&\frac{2}{1+\alpha}\phi_x.
\end{align}
We define the self-similar variables $(\sigma, U)$ as follows: 
\begin{equation}\label{dense,vel}
	\sigma(y,s)=\rho(x,t), \qquad \qquad U(y,s)=u(x,t),
\end{equation}
and obtain the following identities:
\begin{equation}\label{self similar rho u}
	\frac{1}{\alpha} \sigma^\alpha=\frac{e^{-\frac{s}{2}}W+\kappa-Z}{2}, \qquad U=\frac{e^{-\frac{s}{2}}W+\kappa+Z-\kappa_0}{2}.
\end{equation}
By applying the identities in \eqref{identity} and \eqref{self similar rho u} to the system \eqref{eq;rhou}, we get the following system for $(\sigma,U)$:
\begin{align}\label{eq;sigmaU}
	 &\partial_s\sigma+ \left(-e^{\frac{s}{2}}\frac{\dot{\xi}}{(1-\dot\tau)}+\frac{3}{2}y+\frac{2e^{\frac{s}{2}}}{(1+\alpha)(1-\dot\tau)}U \right)\partial_y\sigma +\sigma \frac{\partial_yW+e^{\frac{s}{2}} \partial_yZ}{(1+\alpha)(1-\dot\tau)}=0,\nonumber\\
	&\partial_s U+\frac{2\beta e^{-s}}{(1+\alpha)(1-\dot\tau)} U+ \left(-e^{\frac{s}{2}}\frac{\dot{\xi}}{(1-\dot\tau)}+\frac{3}{2}y+\frac{2 e^{\frac{s}{2}}}{(1+\alpha)(1-\dot\tau)}U +\frac{2e^{\frac{s}{2}}}{(1+\alpha)(1-\dot{\tau})}\sigma^\alpha\right)\partial_yU\nonumber\\
	&=\frac{ 2 e^{\frac{s}{2}} }{(1+\alpha)(1-\dot\tau)}\sigma^\alpha\partial_yZ+\frac{2 e^{\frac{s}{2}}}{(1+\alpha)(1-\dot\tau)}\partial_y\Phi.	
\end{align}
Let
\begin{equation*}
	\widetilde{\sigma}=(1+e^{-3s}y^2)^\frac{1}{3}\left( \sigma-\bar{\rho}\right), \qquad \hbox{and} \qquad \widetilde{U}=(1+e^{-3s}y^2)^\frac{1}{3} U.
\end{equation*} 
where $(1+e^{-3s}y^2)^\frac{1}{3}=(1+x^2)^\frac{1}{3}$ and a positive constant $\bar\rho>0$. From the below equalities for  self-similar  variable $R(y,s)$
\begin{equation*}
   \begin{split}
   \partial_s\left[(1+e^{-3s}y^2)^\frac{1}{3} R\right]=&-\frac{e^{-3s}y^2}{(1+e^{-3s}y^2)^\frac{2}{3}}R+(1+e^{-3s}y^2)^\frac{1}{3}\partial_sR, 
   \end{split}
   \end{equation*}
and
\begin{equation*}
   \begin{split}
   \partial_y\left[(1+e^{-3s}y^2)^\frac{1}{3} R\right] =&\frac{2e^{-3s}y}{3(1+e^{-3s}y^2)^\frac{2}{3}}R+(1+e^{-3s}y^2)^\frac{1}{3} \partial_yR,
   \end{split}	
\end{equation*}
the system \eqref{eq;sigmaU} takes the form of a transport equation with a forcing term for $(\widetilde{\sigma}, \widetilde{U})$, expressed as:
\begin{align}\label{system;sigmaU}
	 &\partial_s\widetilde{\sigma}+ \left(-e^{\frac{s}{2}}\frac{\dot{\xi}}{(1-\dot\tau)}+\frac{3}{2}y+\frac{2e^{\frac{s}{2}}}{(1+\alpha)(1-\dot\tau)}U \right)\partial_y\widetilde{\sigma} \nonumber\\
	 =& \left(-e^{\frac{s}{2}}\frac{\dot{\xi}}{(1-\dot\tau)}+\frac{2e^{\frac{s}{2}}}{(1+\alpha)(1-\dot\tau)}U \right)\frac{2e^{-3s}y}{3(1+e^{-3s}y^2)} \widetilde{\sigma} -\sigma(1+e^{-3s}y^2)^\frac{1}{3} \frac{\partial_yW+e^{\frac{s}{2}} \partial_yZ}{(1+\alpha)(1-\dot\tau)},\nonumber\\
	&\partial_s\widetilde{ U}+\frac{2\beta e^{-s}}{(1+\alpha)(1-\dot\tau)}\widetilde{ U}+ \left(-e^{\frac{s}{2}}\frac{\dot{\xi}}{(1-\dot\tau)}+\frac{3}{2}y+\frac{2 e^{\frac{s}{2}}}{(1+\alpha)(1-\dot\tau)}U +\frac{2e^{\frac{s}{2}}}{(1+\alpha)(1-\dot{\tau})}\sigma^\alpha\right)\partial_y\widetilde{U}\nonumber\\
	=&  \left(-e^{\frac{s}{2}}\frac{\dot{\xi}}{(1-\dot\tau)}+\frac{2 e^{\frac{s}{2}}}{(1+\alpha)(1-\dot\tau)}U +\frac{2e^{\frac{s}{2}}}{(1+\alpha)(1-\dot{\tau})}\sigma^\alpha\right)\frac{2e^{-3s}y}{3(1+e^{-3s}y^2)} \widetilde{U} \nonumber\\
	&+\frac{ 2 e^{\frac{s}{2}} }{(1+\alpha)(1-\dot\tau)}\sigma^\alpha(1+e^{-3s}y^2)^\frac{1}{3} \partial_yZ+\frac{2 e^{\frac{s}{2}}}{(1+\alpha)(1-\dot\tau)}(1+e^{-3s}y^2)^\frac{1}{3} \partial_y\Phi.	
\end{align}
%\begin{align}
%   	 &\partial_s\widetilde{\sigma} + \left(-e^{\frac{s}{2}}\frac{\dot{\xi}}{(1-\dot\tau)}+\frac{3}{2}y+\frac{e^{\frac{s}{2}}}{(1-\dot\tau)} U\right)\partial_y\widetilde{\sigma}\nonumber\\
%   	 =&\left(-e^{\frac{s}{2}}\frac{\dot{\xi}}{(1-\dot\tau)}+\frac{e^{\frac{s}{2}}}{(1-\dot\tau)}U \right)\frac{2 e^{-3s}y}{3(1+e^{-3s}y^2 )} \widetilde{\sigma} -\sigma(1+e^{-3s} y^2 )^\frac{1}{3}  \frac{\partial_yW+e^{\frac{s}{2}} \partial_yZ}{2(1-\dot\tau)},\nonumber\\
%   	 &\partial_s\widetilde{ U}+\frac{\beta e^{-s}}{(1-\dot\tau)}\widetilde{U} +\left(-e^{\frac{s}{2}}\frac{\dot\xi}{(1-\dot\tau)}+\frac{3}{2}y+\frac{e^{\frac{s}{2}}}{(1-\dot\tau)} U+\frac{e^{\frac{s}{2}}}{(1-\dot\tau)}\sigma^\alpha \right)\partial_y\widetilde{U}\nonumber\\
%	=&\left(-e^{\frac{s}{2}}\frac{\dot\xi}{(1-\dot\tau)}+\frac{e^{\frac{s}{2} }}{(1-\dot\tau)} U+\frac{e^{\frac{s}{2}}}{(1-\dot\tau)} \sigma^\alpha \right)\frac{2e^{-3s}y}{3(1+e^{-3s}y^2 )}\widetilde{U}\nonumber\\
%	&+\frac{e^{\frac{s}{2} }}{(1-\dot\tau)}\sigma^\alpha (1+e^{-3s}y^2)^\frac{1}{3}\partial_yZ+\frac{e^{\frac{s}{2} }}{(1-\dot\tau)}(1+e^{-3s}y^2)^\frac{1}{3}\partial_y \Phi.
%\end{align}
For convenience, we rewrite the system \eqref{system;sigmaU} as
\begin{equation}\label{finalrhou}
	\begin{split}
		\partial_s\widetilde{\sigma}+\mathcal{V}_{\widetilde{\sigma}} \partial_y\widetilde{\sigma}=F_{\widetilde{\sigma}}, \qquad \hbox{and} \qquad \partial_s\widetilde{U}+\frac{2\beta e^{-s}}{(1+\alpha)(1-\dot{\tau})}\widetilde{U}+\mathcal{V}_{\widetilde{U}} \partial_y\widetilde{U}=F_{\widetilde{U}},
	\end{split}
\end{equation}
where transport velocity terms and forcing terms are given by:
\begin{align}\label{rho u vel}
	\mathcal{V}_{\widetilde{\sigma}}=& -e^{\frac{s}{2}}\frac{\dot{\xi}}{(1-\dot\tau)}+\frac{3}{2}y+\frac{2e^{\frac{s}{2}}}{(1+\alpha)(1-\dot\tau)}U ,\nonumber\\
	 \mathcal{V}_{\widetilde{U}}=&-e^{\frac{s}{2}}\frac{\dot{\xi}}{(1-\dot\tau)}+\frac{3}{2}y+\frac{2 e^{\frac{s}{2}}}{(1+\alpha)(1-\dot\tau)}U +\frac{2e^{\frac{s}{2}}}{(1+\alpha)(1-\dot{\tau})}\sigma^\alpha ,\nonumber\\
		 F_{\widetilde{\sigma}}=&\left(-e^{\frac{s}{2}}\frac{\dot{\xi}}{(1-\dot\tau)}+\frac{2e^{\frac{s}{2}}}{(1+\alpha)(1-\dot\tau)}U \right)\frac{2e^{-3s}y}{3(1+e^{-3s}y^2)} \widetilde{\sigma} -\sigma(1+e^{-3s}y^2)^\frac{1}{3} \frac{\partial_yW+e^{\frac{s}{2}} \partial_yZ}{(1+\alpha)(1-\dot\tau)},\nonumber\\
		 F_{\widetilde{U}}=&\left(-e^{\frac{s}{2}}\frac{\dot{\xi}}{(1-\dot\tau)}+\frac{2 e^{\frac{s}{2}}}{(1+\alpha)(1-\dot\tau)}U +\frac{2e^{\frac{s}{2}}}{(1+\alpha)(1-\dot{\tau})}\sigma^\alpha\right)\frac{2e^{-3s}y}{3(1+e^{-3s}y^2)} \widetilde{U} \nonumber\\
	&+\frac{ 2 e^{\frac{s}{2}} }{(1+\alpha)(1-\dot\tau)}\sigma^\alpha(1+e^{-3s}y^2)^\frac{1}{3} \partial_yZ+\frac{2 e^{\frac{s}{2}}}{(1+\alpha)(1-\dot\tau)}(1+e^{-3s}y^2)^\frac{1}{3} \partial_y\Phi.	
\end{align}

\section{Main result}\label{sec;main result}

In this section, we impose the assumptions of initial data to construct the shock-type singularity in the system \eqref{simple model} and state the precise main result. We take a  sufficiently large parameter $M>0$ and a sufficiently small parameter $\epsilon>0$. We define $l=\frac{1}{M}$.
\subsection{Assumptions on the initial data.}\label{assump ini}
 To define the modulation variables, we prescribe their initial data as follows:
\begin{equation}\label{ini mod}
	\kappa_0=\kappa(0), \qquad \tau_0=\tau(0)=\epsilon,\qquad \hbox{and} \qquad \xi_0=\xi(0)=0,
\end{equation}
where
\begin{equation}\label{kappa constant}
	\kappa_0\geq \frac{5(1+\alpha)}{\alpha}.
\end{equation}
Let 
\begin{equation}
\bar{\phi}=\left(\frac{\alpha\kappa_0}{2}\right)^{\frac{1}{\alpha}}, 
\qquad 
\tilde{\alpha}=\frac{1}{\alpha}-\left[\frac{1}{\alpha}\right]-2,
\end{equation}
where \( \left[\frac{1}{\alpha} \right] \) denotes the Gauss bracket, i.e., the greatest integer less than or equal to \( \frac{1}{\alpha} \).  
Assume that the smooth initial data \( (q_0, u_0, \phi_0) \) satisfy
\begin{equation}\label{initial sobolev}
q_0-\frac{\kappa_0}{2}, \ \
q_0^{\tilde{\alpha}}-\left(\frac{\kappa_0}{2}\right)^{\tilde{\alpha}}, \ \
u_0 \in H^m, 
\ \  \text{and} \ \ 
\phi_0-\bar{\phi} \in H^{m+1}, \quad \text{for some $m \in \mathbb{N}$ with $m \geq 5$.}
\end{equation}
%Let $\bar{\phi}=\left(\frac{\alpha\kappa_0}{2}\right)^\frac{1}{\alpha}$, and $\tilde{\alpha}=\frac{1}{\alpha}-\left[\frac{1}{\alpha}\right]-2$. Here, $\left[\frac{1}{\alpha}\right]$ denotes the Gauss bracket,  representing the greatest integer less than or equal to $\frac{1}{\alpha}$.  The smooth initial data for ($q,u,\phi$) belongs to the following Sobolev space. 
%\begin{equation}\label{initial sobolev}
%	q_0-\frac{\kappa_0}{2},\ q^{\tilde{\alpha}}_0-\left(\frac{\kappa_0}{2}\right)^{\tilde{\alpha}}, \ u_0 \in H^m, \qquad \text{and}\qquad  \ \phi_0-\bar{\phi} \in H^{m+1}, \quad \text{$m\in\mathbb{N}$ with $m\geq5$.} 
%\end{equation}
% Let $m\in\mathbb{N}$ with $m\geq5$. Then the smooth initial data for ($q,u,\phi$) belongs to the following Sobolev space. 
%\begin{equation}\label{initial sobolev}
%	q_0-\frac{\kappa_0}{2},\ q^{\tilde{\alpha}}_0-\left(\frac{\kappa_0}{2}\right)^{\tilde{\alpha}}, \ u_0 \in H^m, \qquad \text{and}\qquad  \ \phi_0-\bar{\phi} \in H^{m+1},
%\end{equation}
%where $\bar{\phi}=\left(\frac{\alpha\kappa_0}{2}\right)^\frac{1}{\alpha}$, and $\tilde{\alpha}=\frac{1}{\alpha}-\left[\frac{1}{\alpha}\right]-2$. Here, $\left[\frac{1}{\alpha}\right]$ denotes the Gauss bracket,  representing the greatest integer less than or equal to $\frac{1}{\alpha}$.\\
In addition, the weighted initial data also satisfies the following: 
\begin{equation}\label{initial weight sobolev}
	\left(q_0-\frac{\kappa_0}{2}\right)x,\ u_0x, \ (\rho_0-\bar{\rho})x\in H^{3}, \qquad \hbox{and} \qquad \left(\phi_0-\bar{\phi}\right)x\in H^{4},   
\end{equation}
where $\bar{\rho}=\bar{\phi}$. We impose constraints of $w_0$ and its several derivatives at $x=0$ and their bounds:
\begin{align}\label{ini 1}
w_0(0) = \kappa_0 ,\qquad \partial_x w_0(0) = - \epsilon^{-1}, \qquad \partial_x^2 w_0(0) =0, \qquad \hbox{and} \qquad\partial_x^3 w_0(0) =  6 \epsilon^{-4}, 
\end{align}
and
\begin{align}\label{ini 2}
 \|\partial_x w_0\|_{L^\infty}\leq \epsilon^{-1}, \quad \|\partial_x^2 w_0\|_{L^\infty}\leq \epsilon^{-\frac{5}{2}}, \ \ \|\partial_x^3 w_0\|_{L^\infty}\leq 7 \epsilon^{-4},\ \ \ \hbox{and} \ \ \|\partial_x^4 w_0\|_{L^\infty}\leq \epsilon^{-\frac{11}{2}},
\end{align}
which imply that the  global minimum of $\partial_x w_0$ is at $x=0$.
Furthermore, to guarantee the positivity of the initial density, we assume
\begin{align}\label{ini W}
\|w_0 (\cdot) - \kappa_0\|_{L^\infty} \leq   \frac{\kappa_0}{8}.
\end{align}
Additionally, the initial data of $z$ and $\phi$ satisfies
\begin{align}\label{ini phi Z}
  \|z_0\|_{C^4}+\|\phi_0-\bar{\phi}\|_{C^5}   \leq 1.
\end{align}
To get the structure of a unique singularity and its regularity, we assume a more detailed behavior of $w_0$. \\
For $|x|\leq \epsilon^\frac{3}{2}l$ with  $l=\frac{1}{M}$ and for $n=1,2,3,4:$
\begin{align}\label{ini;weight0}
        &\left|\epsilon^{-\frac{1}{2}}w_0(x)-\overline{W}\left(\frac{x}{\epsilon^\frac{3}{2}}\right)-\epsilon^{-\frac{1}{2}}\kappa_0\right|\leq \epsilon^\frac{1}{5}l^4,\ \ \hbox{and} \ \  \left|\epsilon^{\frac{3}{2}n-\frac{1}{2}}\partial^n_xw_0(x)-\partial^n_y\overline{W}\left(\frac{x}{\epsilon^\frac{3}{2} }\right)\right|\leq  \epsilon^\frac{1}{5} l^{4-n}.
\end{align}
For  $\epsilon^\frac{3}{2} l\leq|x|\leq  \frac{5}{4}M\epsilon+1$:
\begin{align}\label{ini;weight1}
		&(1+\epsilon^{-3} x^2)^{-\frac{1}{6}}\left|\epsilon^{-\frac{1}{2}} w_0(x)-\overline{W}\left(\frac{x}{\epsilon^\frac{3}{2}} \right)-\epsilon^{-\frac{1}{2}}\kappa_0\right|\leq \frac{1}{3} \epsilon^\frac{1}{5},\nonumber\\
\text{and}	\quad	& (1+\epsilon^{-3} x^2)^{\frac{1}{3}}\left|\epsilon\partial_x w_0(x)-\partial_y\overline{W}\left(\frac{x}{\epsilon^\frac{3}{2}}\right)\right|\leq \frac{1}{3}\epsilon^\frac{1}{5}.
\end{align} 
For   $1\leq |x|$: 
\begin{equation}\label{ini;weight2}
	(1+\epsilon^{-3} x^2)^{\frac{1}{3}}\left|\epsilon\partial_x w_0(x)\right|\leq \frac{1}{2}.
\end{equation}
Here,  $\overline W$ is the self-similar Burgers profile in \eqref{eq;burger}.
Finally, we impose the bounds for the weighted initial data of $\rho, u,\phi$ and $\partial_xz$: For all  $x \in\mathbb{R}$,
\begin{equation}\label{ini final}
	(1+x^2)^\frac{1}{3}\max\left\{|\rho-\bar{\rho}|,|u|,|\phi-\bar{\phi}|,|\phi_x|,|\phi_{xx}|\right\}(x,0)\leq 1, 
\end{equation}
and
\begin{align}\label{ini 3'}
 \left(1 + \left(\frac{x}{\epsilon^\frac{3}{2}}\right)^2\right)^\frac{1}{3}\epsilon^\frac{3}{2} \left|\partial_x z_0(x) \right| 
&\leq  \frac{1}{2}.
\end{align}
\begin{remark}
The assumptions in \eqref{ini phi Z}--\eqref{ini 3'} imply that the initial data are small perturbations of the steady states. In the modulation analysis, treating them as perturbations around these steady states allows us to track the evolution of w and capture its behavior up to the formation of blow-up.
\end{remark}

\subsection{Statement of the main result.} 
Our main objective is to establish the following theorem on the shock-type singularity of the HPC system  \eqref{simple model}. To align the time scaling in the theorem with the system \eqref{riemann}, we interpret the time variable $t$  in the system \eqref{simple model} as $\frac{1+\alpha}{2}t$.  
\begin{theorem}[\bf Shock-type singularity] \label{delicate main}
 Suppose that the smooth initial data $(q_0, u_0, \phi_0)$ satisfy assumptions \eqref{kappa constant}--\eqref{ini 3'}. For a sufficiently large $M=M(\alpha,\beta,\kappa_0)>0$ and a sufficiently small $\epsilon(\alpha,\beta, \kappa_0, M)>0$, there exists a unique smooth solution to the HPC system (\ref{simple model}) which blows up at finite time $T^*<\infty$ such that the following holds: 
\begin{enumerate}[(1)]
\item $\left(q-\frac{\kappa_0}{2} ,u\right)\in C([0,T^*);H^m),$ and $ \left(\phi-\bar{\phi}\right)\in C([0,T^*);H^{m+1}) \cap L^2((0,T^*); H^{m+2}).$
\item The blow-up time $T^*$ satisfies $0<T^*\leq\frac{3}{2} \epsilon$, and the position  $x^*  = \lim_{t\to T^*} \xi(t)$ satisfies $|x^*|\leq 6M \epsilon $.

\item $\lim_{t \to T^*}    w_x(\xi(t),t)  = -\infty $ and the blow-up rate satisfies
\begin{align*}
	 \| w_x(\cdot,t)\|_{L^\infty} \thicksim \frac{1}{T^*-t} \quad  as  \quad t \to T^*.
\end{align*} 
Moreover, we have
\begin{equation*}
 		\lim_{t\to T^*}  q_x(\xi(t),t)=\lim_{t\to T^*} u_x(\xi(t),t)=-\infty.
 	\end{equation*}
\item $w( \cdot , T^*)$ has a unique blow-up point $x^*$ with cusp singularity. Specifically, we have 
\begin{equation*}
	 w_x(x, T^*)\thicksim -|x-x^*|^{-\frac{2}{3}}\qquad   \hbox{for all}\qquad  0<|x-x^*|< 1.
\end{equation*} 
\item $\sup_{t\in[0, T^*]} \left(   \| w(t)\|_{ C^\frac{1}{3}}+\| z(t)\|_{ C^1}+\| \phi(t)\|_{ C^{2}} \right) < +\infty$. $w$ is $C^1$ differentiable at $x\neq x^*$. 
\end{enumerate}   
\end{theorem} 

\begin{remark}\label{structure}	From Theorem \ref{delicate main} $(4)$, $ w(x,T^*)$ behaves as $-(x-x^*)^\frac{1}{3}$ for $|x-x^*|< 1$. Since it can be readily shown that the perturbation for $z$ remains small, with $\|z_0\|_{L^\infty}$ sufficiently small, until the singularity forms, we can capture the structure of $q$ and $u$ near the blow-up point as $-(x-x^*)^\frac{1}{3}$ for $|x-x^*|< 1$.
\end{remark}

\begin{remark}\label{rmk;relation}
	To ensure the formation of  a shock-type singularity, $\epsilon>0$ need to be sufficiently small, as $\beta$ increases. This can be seen in Section \ref{sec;closboot} and \ref{sec;stabilW}, where the terms $F_Z$ and $F_W$ in \eqref{g and force} contain $\beta e^{-\frac{s}{2}} \frac{W}{2}$ and $\beta\frac{Z}{2}$ respectively. These terms contribute to a reduction in $\epsilon$ as $\beta$ increases.
\end{remark}

\section{Local well-posedness and continuation criterion}\label{sec;localconti}

 In this section, we study the local well-posedness for classical solutions to the system \eqref{simple model} and its continuation criterion. Moreover, we prove the local well-posedness for weighted variables \eqref{weight variable} within the system \eqref{simple model}.
\begin{theorem}\label{simple result}
	Let $m\geq 3$ be an integer and $c>0$. If $(q_0-c, q_0^{\tilde{\alpha}}-c^{\tilde\alpha}, u_0)\in H^m$ and $\phi_0-(\alpha c)^\frac{1}{\alpha} \in H^{m+1}$ where $\tilde{\alpha}=\frac{1}{\alpha}-\left[\frac{1}{\alpha}\right]-2$, then there exists $T>0$ such that the system  \eqref{simple model} has a  unique solution $(q,u,\phi)$ with
	\begin{equation} \label{local well1}
	(q-c, q^{\tilde{\alpha}}-c^{\tilde\alpha}, u)\in C([0,T);H^m),\ \ and\ \ \phi-(\alpha c)^\frac{1}{\alpha}\in C([0,T);H^{m+1})\cap L^2((0,T);H^{m+2}).	
	\end{equation}
 Moreover, the solution continues to exist on $[0,T]$ if and only if
	\begin{equation}\label{continucri}
		\int^{T}_0 \|(q_x,u_x)(t')\|_{L^\infty} dt'<\infty.
	\end{equation}
\end{theorem}
\begin{remark}
The aforementioned theorem extends to both the case $\beta\geq0$, and to  general spatial dimensions d, where $m>\frac{d}{2}+2$.
\end{remark}

\begin{remark}
Due to  the assumption that $q_0^{\tilde{\alpha}}-c^{\tilde{\alpha}}\in H^m$ in the aforementioned theorem, the initial density is bounded away from zero. In particular,
\begin{equation}\label{density positive condition}
    \left\|\frac{1}{q_0}\right\|_{L^\infty}<\infty.
\end{equation}
\end{remark}

 We need the following lemma:
\begin{lemma}\label{big lemma} Let $u,v,w\in \dot{H}^k\cap L^\infty$ with $k,n\in\mathbb{N}$. Then the following holds:
 \begin{enumerate}[(1)]
\item  Homogeneous Sobolev product inequality:
 	\begin{equation*}
	\|\partial^k_x(uv)\|_{L^2}\lesssim_{k} \|u\|_{L^\infty}\|\partial^k_x v\|_{L^2}+\|v\|_{L^\infty}\|\partial^k_xu\|_{L^2}. 
\end{equation*}
    \item   Commutator estimates: 
\begin{equation*}
	\|[\partial^k_x,u\partial]v\|_{L^2}\lesssim_{k}  \|\partial_x u\|_{L^\infty}\|\partial^k_xv\|_{L^2}+\|\partial_x v\|_{L^\infty}\|\partial^k_x u\|_{L^2},
\end{equation*}
where $[\partial^k_x,u\partial_x]v=\partial^k_x(u\partial_x v)-u \partial^{k+1}_xv$.
  \item Auxiliary estimates:
  \begin{equation*}
	\|[\partial^k_x[(v+1)^n w]\|_{L^2}\lesssim_{k,n}  \left(\|v\|^n_{L^\infty}+1\right)\|\partial^k_x w\|_{L^2}+\| w\|_{L^\infty}(\|v\|^{n-1}_{L^\infty}+1)\|\partial^k_x v\|_{L^2},
\end{equation*}
and
\begin{equation*}
\begin{split}
	\|[\partial^k_x[u(v+1)^n w]\|_{L^2}\lesssim_{k,n} & \left(\|v\|^n_{L^\infty}+1\right)\left(\|u\|_{L^\infty}\|\partial^k_x w\|_{L^2}+\| w\|_{L^\infty}\|\partial^k_x u\|_{L^2}\right)\\
	&+\|u\|_{L^\infty}\| w\|_{L^\infty}\|v\|^{n-1}_{L^\infty}\|\partial^k_x v\|_{L^2}.
\end{split}	
\end{equation*}
\end{enumerate}
\end{lemma}
\begin{proof}
	 We recall these classical  inequalities in $(1)$ and $(2)$ from \cite{bahourifourier}. The proof of $(3)$ is based on $(1)$. Straightforward computations yield
	 	  \begin{align*}
	\|[\partial^k_x[(v+1)^n w]\|_{L^2}&\leq \|\partial^k_x[((v+1)^n-1) w]\|_{L^2}+\|\partial^k_x w\|_{L^2}\\
	&\lesssim_k \|(v+1)^n-1\|_{L^\infty}\|\partial^k_x w \|_{L^2}+\|w\|_{L^\infty}\|\partial^k_x[(v+1)^n-1 ]\|_{L^2}+\|\partial^k_x w\|_{L^2}\\
	&\lesssim_{k,n} \left(\|v\|^n_{L^\infty}+1\right)\|\partial^k_x w\|_{L^2}+\| w\|_{L^\infty}(\|v\|^{n-1}_{L^\infty}+1)\|\partial^k_x v\|_{L^2},
\end{align*}
	 and
	\begin{align*}
			\|[\partial^k_x[u(v+1)^n w]\|_{L^2}&\lesssim_k  \|u\|_{L^\infty}\|\partial^k_x[(v+1)^n w]\|_{L^2}+\|(v+1)^n w\|_{L^\infty}\|\partial^k_x u\|_{L^2}\\
			&\lesssim_{k,n}   (\|u\|_{L^\infty}(\left(\|v\|^n_{L^\infty}+1\right)\|\partial^k_x w\|_{L^2}+\| w\|_{L^\infty}(\|v\|^{n-1}_{L^\infty}+1)\|\partial^k_x v\|_{L^2})\\
			&\quad+(\|v\|^n_{L^\infty}+1)\| w\|_{L^\infty}\|\partial^k_x u\|_{L^2}\\
			&\lesssim_{k,n}   \left(\|v\|^n_{L^\infty}+1\right)\left(\|u\|_{L^\infty}\|\partial^k_x w\|_{L^2}+\| w\|_{L^\infty}\|\partial^k_x u\|_{L^2}\right)\\
			&\quad+\|u\|_{L^\infty}\| w\|_{L^\infty}\|v\|^{n-1}_{L^\infty}\|\partial^k_x v\|_{L^2}.
	\end{align*}
\end{proof}

\begin{proof}\textbf{(Proof of Theorem \ref{simple result}.)}
The constants $C$ may be dependent on the parameters $m$ and $\alpha$, but not on the solution. The proof is divided into several steps. For simplicity, let $c=1$. To get a suitable priori estimates, we consider the following system from  the system \eqref{simple model} by denoting $\eta=q-1$, $h=q^{\tilde{\alpha}}-1$: 

\begin{align}\label{local system}
		\eta_t+u \eta_x +\alpha \eta u_x+\alpha  u_x=&0,\nonumber\\
		 u_t+u u_x +\alpha \eta \eta_x+\alpha \eta_x+\beta u=& \phi_x,\nonumber\\
		 h_t+u h_x +c_\alpha h u_x+c_\alpha  u_x=&0,\nonumber\\
	    \phi_{xt}-\phi_{xxx} +\phi_x=&\tilde{c}_\alpha (h+1)(\eta+1)^n\eta_x,
\end{align}
 where $\tilde{\alpha}=\frac{1}{\alpha}-\left[\frac{1}{\alpha}\right]-2$, $c_\alpha=\alpha\tilde{\alpha}$, $\tilde{c}_\alpha=\alpha^{\frac{1}{\alpha}-1}$, and $n=\left[\frac{1}{\alpha}\right]+1$.\\
\noindent \textbullet \textbf{ (A priori estimates)} We begin with $L^2$ estimates:
\begin{equation*}
	\begin{split}
		\frac{1}{2}\frac{d}{dt}\| (\eta, u)\|^2_{L^2}+\beta \|u\|^2_{L^2}=&  -\int (u \eta_x) \eta\, dx -\alpha\int (\eta  u_x) \eta\, dx-\alpha\int  u_x\eta\, dx\\
		&-\int (u  u_x) u\, dx-\alpha\int (\eta \eta_x)u\, dx-\alpha\int \eta_x u\, dx + \int \phi_x  u\, dx\\
		\leq & C(\|u_x\|_{L^\infty}+\|\eta_x\|_{L^\infty}+1)(\|u\|^2_{L^2}+\|\eta\|^2_{L^2})+\frac{1}{4}\|\phi_x\|^2_{L^2}.
	\end{split}
\end{equation*}
 For $H^k$ estimates with $k=1,\cdots, m$, we take the derivative $\partial^k_x$ and $\partial^{k-1}_x$ for $(\eta,u)$ and $(h,\phi_x)$ in the system \eqref{local system} respectively. Taking the $L^2$ inner product and using Lemma \ref{big lemma} (2), (3)  gives
\begin{align*}
		&\frac{1}{2}\frac{d}{dt}\left(\|\partial^k_x (\eta, u)\|^2_{L^2}+\|\partial^{k-1}_x (h, \phi_x)\|^2_{L^2}\right)+\|\partial^k_x \phi\|^2_{H^1}+\beta\|\partial^k_x u\|^2_{L^2}\\
		=& -\int\partial^k_x\eta\cdot[\partial^k_x,u\partial_x]\eta\, dx 
		   -\alpha\int\partial^k_x\eta\cdot[\partial^k_x,\eta\partial_x]u\, dx 
		    -\int\partial^k_x u\cdot[\partial^k_x,u\partial_x] u\, dx \\
		    &-\alpha\int\partial^k_x u\cdot[\partial^k_x,\eta\partial_x] \eta\, dx  -\int\partial^{k-1}_x h\cdot[\partial^{k-1}_x,u\partial_x] h\, dx -c_\alpha\int\partial^{k-1}_x h\cdot[\partial^{k-1}_x,h\partial_x]u\, dx\\
		    & -c_\alpha\int\partial^{k-1}_x h\cdot\partial^{k-1}_xu_x\, dx+\frac{1}{2}\int \partial_x u\left(|\partial^k_x\eta|^2+|\partial^k_x u|^2+|\partial^{k-1}_x h|^2\right)\, dx\\
		    &+\alpha\int \partial_x \eta\partial^k_x u\partial^k_x \eta\, dx-c_\alpha\int h\partial^{k-1}_x h\partial^k_x u\, dx +\int \partial^k_x\phi_x\partial^k_x u\, dx\\
		    & +\tilde{c}_\alpha\int  \partial^{k-1}_x \phi_x \partial^{k-1}_x\left((h+1)(\eta+1)^n\eta_x\right)\, dx\\
		\lesssim & \left(\|\partial_x u\|_{L^\infty}\|\partial^k_x \eta\|_{L^2}+\|\partial_x \eta\|_{L^\infty}\|\partial^{k}_x u\|_{L^2}\right)\|\partial^k_x \eta\|_{L^2} +\|\partial_x u\|_{L^\infty}\|\partial^k_x u\|^2_{L^2}\\
		&+\left\|\partial^k_x u\right\|_{L^2}\left\|\partial^k_x\phi_x\right\|_{L^2}+(\|\partial_x h\|_{L^\infty}\|\partial^{k-1}_x u\|_{L^2}+\|\partial_x u\|_{L^\infty}\|\partial^{k-1}_xh\|_{L^2})\|\partial^{k-1}_xh\|_{L^2}\\
		&+(\|h\|_{L^\infty}+1)(\|\partial^k_x u\|^2_{L^2}+\|\partial^{k-1}_x h\|^2_{L^2})+\left\|\partial^{k-1}_x  \left((h+1)(\eta+1)^n\eta_x\right)
		 \right\|_{L^2}\|\partial^{k-1}_x\phi_x\|_{L^2}\\
		\leq &  C\left(\| u_x\|_{L^\infty}+\| \eta_x\|_{L^\infty}+\|h\|_{W^{1,\infty}}+1\right)\left(\|\partial^{k-1}_x u\|^2_{H^1}+\|\partial^k_x \eta\|^2_{L^2}+\|\partial^{k-1}_x h\|^2_{L^2}\right)\\
		&+\frac{1}{4}\|\partial^{k+1}_x\phi\|^2_{L^2}+C\bigg[(\|\eta\|^n_{L^\infty}+1)\left(\|h\|_{L^\infty}\|\partial^k_x\eta\|_{L^2}+\|\eta_x\|_{L^\infty}\|\partial^{k-1}_x h\|_{L^2}+\|\partial^k_x\eta\|_{L^2}\right)\\
		&\quad\qquad\qquad\qquad \quad\qquad+(\|h\|_{L^\infty}+1)\|\eta_x\|_{L^\infty}\|\eta\|^{n-1}_{L^\infty}\|\partial^{k-1}_x\eta\|_{L^2}\bigg]\|\partial^{k-1}_x\phi_x\|_{L^2}\\
		\leq &  C(\|h\|_{L^\infty}+1)(\|\eta\|^n_{L^\infty}+1)\left(\|( u_x, \eta_x,h_x)\|_{L^\infty} +1\right)\left(\|\partial^{k-1}_x (u,\eta)\|^2_{H^1}+\|\partial^{k-1}_x (h,\phi_x)\|^2_{L^2}\right)\\
		&+\frac{1}{4}\|\partial^{k+1}_x\phi\|^2_{L^2}.
\end{align*}
Summing these estimates yields
\begin{align}\label{cont crit}
		&\frac{1}{2}\frac{d}{dt}\left(\| (\eta, u)\|^2_{H^m}+\|(h,\phi_x)|^2_{H^{m-1}}\right)+\frac{1}{2}\| \phi_x\|^2_{H^m}+\beta\|u\|^2_{H^m}\nonumber\\
		\lesssim & (\|h\|_{L^\infty}+1)(\|\eta\|^n_{L^\infty}+1)\left(\|( u_x, \eta_x,h_x)\|_{L^\infty} +1\right)\left(\| (\eta, u)\|^2_{H^m}+\| (h,\phi_x)\|^2_{H^{m-1}}\right).
\end{align}
For the quantity
\begin{equation*}
	X_m=\| (\eta, u)\|^2_{H^m}+\|(h,\phi_x)|^2_{H^{m-1}}+1,
\end{equation*}
we obtain the differential inequality:
\begin{equation*}
	\frac{d}{dt}X_m\leq C X^N_m,
\end{equation*}
for some $N\in\mathbb{N}$.\\

\noindent \textbullet  \textbf{ (Existence)} From the aforementioned a priori estimates, it is possible to fix some small time interval $[0,T]$, on which we have $X(t)\leq 2X(0)$. We obtain a sequence of smooth solutions $(\eta^\epsilon, u^\epsilon,  h^\epsilon, \phi_x^\epsilon)$ on $[0,T]$ by adapting the mollification method in \cite{majda2002vorticity}, that is, by applying the mollifier $\mathcal{J}^\epsilon$
  to the system. Specifically, we mollify each variable in \eqref{local system}, except for the time derivative, and apply the mollification to each nonlinear term as a whole. The sequence $(\eta^\epsilon, u^\epsilon)$ and the sequence  $(h^\epsilon, \phi_x^\epsilon)$  are uniformly bounded in $C([0,T];H^m)$  and $C([0,T];H^{m-1})$, respectively. A consequence of these uniform bounds and the weak * compactness theorem is that a sequence $(\eta^\epsilon, u^\epsilon,  h^\epsilon, \phi_x^\epsilon)$ has a subsequence such that converges weakly $*$ to 
\begin{equation*}
	 (\eta,u) \in  L^\infty([0,T];H^m), \quad \hbox{and} \quad 
	(h,\phi_x) \in  L^\infty([0,T];H^{m-1}).
\end{equation*}
Moreover, by using  Lemma \ref{big lemma} (1), we can show that for all $g\in C^\infty_c(\mathbb{R})$, the sequence $(g\eta^\epsilon, gu^\epsilon)$, and the sequence  $(gh^\epsilon, g\phi_x^\epsilon)$  are  bounded in $L^\infty([0,T];H^{m-1})$  and $L^\infty([0,T];H^{m-2})$, respectively.
Therefore, applying the Aubin-Lions lemma, we can extract a subsequence of $(g\eta^\epsilon, gu^\epsilon,  gh^\epsilon, g\phi_x^\epsilon)$ for each $m'<m$ so that converges to
\begin{equation*}
	 (g\eta, gu) \in  L^\infty([0,T];H^{m'}), \quad  \hbox{and} \quad (gh, g\phi_x)  \in L^\infty([0,T];H^{m'-1}).
\end{equation*}
So, we conclude the  $(\eta,u,h,\phi_x)$ satisfies the system \eqref{local system} for $[0,T]\times\mathbb{R}$ a.e. Since the equation for $\phi_x$ in the system \eqref{local system} can be considered as  a heat equation with zeroth order term and force, we have $\phi_x\in C([0,T];H^m)\cap L^2([0,T];H^{m+1})$. Sequentially, we obtain	
\begin{equation*}
	\eta,u\in C([0,T];H^m),\ \ \text{and} \ \  h\in C([0,T];H^{m-1}).
\end{equation*}
Applying Fa\`a di Bruno's formula to $h(=q^{\tilde{\alpha}}-1)$ yields
\begin{equation*}
	h\in C([0,T];H^m).
\end{equation*}
We note that
\begin{align*}
  \rho=&\alpha^{\frac{1}{\alpha}}\sum^{\left[\frac{1}{\alpha}\right]+2}_{k=1}\begin{pmatrix}
		 	\left[\frac{1}{\alpha}\right]+2\\
		 	k
		 \end{pmatrix}(q-1)^kq^{\tilde{\alpha}}+\alpha^{\frac{1}{\alpha}}\left(q^{\tilde{\alpha}}-1\right)+\alpha^{\frac{1}{\alpha}}.
\end{align*}
Using Leibniz rule  and Sobolev embedding, we have 
\begin{align}\label{rho space}
  \rho-\alpha^{\frac{1}{\alpha}}\in C([0,T];H^m).
\end{align}
So, we conclude
\begin{equation*}
	\eta,u,h\in C([0,T];H^m), \quad  \text{and} \quad \phi-\alpha^{\frac{1}{\alpha}} \in C([0,T];H^{m+1})\cap L^2([0,T];H^{m+2}).
\end{equation*}
 
\noindent \textbullet \textbf{ (Uniqueness)} Let $(q_1,u_1,\phi_1)$ and $(q_2,u_2,\phi_2)$ be two solutions of the system (\ref{simple model}) with the same initial data. We denote $(q,u,\phi)=(q_1-q_2,u_1-u_2,\phi_1-\phi_2)$. Then $(q,u,\phi)$ satisfies the following equations: 
\begin{equation*}
		\begin{split}
		\partial_t q+u\partial_x q_1+u_2\partial_x q +\alpha q \partial_x u_1+\alpha q_2 \partial_x u=&0,\\
		\partial_t u+u\partial_x u_1+u_2\partial_x u +\alpha q \partial_x q_1+\alpha q_2\partial_x q+\beta u=& \phi_x,\\
		\partial_t \phi-\partial^2_x\phi+\phi=&\alpha^\frac{1}{\alpha}\left( q^\frac{1}{\alpha}_1-q^\frac{1}{\alpha}_2\right). 
		\end{split}
\end{equation*}
Taking the $L^2$ inner product to the above system with $(q,u,\phi)$, and using that
\begin{equation*}
	q^\frac{1}{\alpha}_1-q^\frac{1}{\alpha}_2=\frac{1}{\alpha}q\int^1_0(q_2+s q)^{\frac{1}{\alpha}-1}ds,
\end{equation*}
we have
\begin{equation*}
	\begin{split}
		&\frac{1}{2}\frac{d}{dt}\| (q, u, \phi)\|^2_{L^2}+\|\phi\|^2_{H^1}+\beta\|u\|^2_{L^2}\\
		=&  -\int (u \partial_x q_1) q\, dx-\int (u_2 \partial_x q ) q\, dx  -\alpha\int (q\partial_x u_1) q\, dx-\alpha\int (q_2 \partial_x u) q\, dx-\int (u \partial_x u_1) u\, dx\\
		&-\int (u_2 \partial_x u) u\, dx-\alpha\int (q\partial_x q_1)u\, dx-\alpha\int (q_2\partial_x q)u\, dx+\int \phi_x u\, dx +\alpha^\frac{1}{\alpha} \int (q_1^\frac{1}{\alpha}-q_2^\frac{1}{\alpha})\phi\, dx\\
		\leq & C\bigg(\|\left(\partial_x q_1,\partial_x q_2,\partial_x u_1,\partial_x u_2\right)\|_{L^\infty}+\|(q_1,q_2)\|_{L^\infty}\\
		&\quad\quad\qquad\qquad\qquad\qquad\qquad \quad+\left\|\int^1_0(q_2+s q)^{\frac{1}{\alpha}-1}ds\right\|_{L^\infty}+1\bigg)\|(q,u,\phi)\|^2_{L^2}+ \frac{1}{4}\|\phi_x\|^2_{L^2}.
	\end{split}
\end{equation*}
Applying Gr\"onwall inequality yields uniqueness.\\

\noindent \textbullet \textbf{ (Continuation criterion)} From the equation of $q$,
\begin{equation*}
	q_t+uq_x+\alpha qu_x=0,
\end{equation*}
we get
\begin{equation*}
	\left(\frac{1}{q}\right)_t+u\left(\frac{1}{q}\right)_x-\alpha u_x\left(\frac{1}{q}\right)=0.
\end{equation*}
Evaluating the preceding equation  along the flow produced by $u$, together with \eqref{density positive condition}, yields
\begin{equation}\label{cont1}
	\left\|\frac{1}{q}(t)\right\|_{L^\infty}\lesssim \exp\left(C\int^t_0\|u_x(s)\|_{L^\infty}ds\right)\left\|\frac{1}{q_0}\right\|_{L^\infty},
\end{equation}
By a similar argument, the following estimates also hold for $h+1$, and $q$
\begin{equation}\label{cont2}
	\left\|(h+1)(t)\right\|_{L^\infty}\lesssim \exp\left(C\int^t_0\|u_x(s)\|_{L^\infty}ds\right)\left\|h_0+1\right\|_{L^\infty}, 
\end{equation}
and 
\begin{equation}\label{cont3}
 \left\|q(t)\right\|_{L^\infty}\lesssim \exp\left(C\int^t_0\|u_x(s)\|_{L^\infty}ds\right)\left\|q_0\right\|_{L^\infty}.
\end{equation}
By integrating \eqref{cont crit} in time and using \eqref{cont1}, \eqref{cont2}, and \eqref{cont3}, we get
\begin{equation*}
\begin{split}
X_m(t)+\int_0^t \|\phi_x(s)\|_{H^m}^2 \, ds\lesssim\exp\!\left(\exp\!\left(C\bigl(X_m(0)\bigr)\int_0^t \|(q_x,u_x)(s)\|_{L^\infty} \, ds\right)\right)X_m(0).
\end{split}
\end{equation*}
The above estimate yields the desired continuation criterion. 
\end{proof}

From \eqref{rho space} in \textbf{(Existence)}  of Theorem \ref{simple result}, we immediately obtain the following corollary.
\begin{corollary}\label{cor local}
For the classical solutions constructed in Theorem \ref{simple result}, we have
\begin{equation*}
  \rho-(\alpha c)^{\frac{1}{\alpha}} \in C([0,T);H^m).
\end{equation*}
\end{corollary}
 In order to apply a bootstrap argument of weighted self-similar variables, we establish the following theorem. For simplicity in the theorem's statement, we introduce the notations:
\begin{equation}\label{weight variable}
	\tilde{\eta}=(q-c)x,\qquad \tilde{\rho}=(\rho-(\alpha c)^\frac{1}{\alpha} )x,\qquad \tilde{u}=ux, \qquad \hbox{and}\qquad \tilde{\phi}=(\phi-(\alpha c)^\frac{1}{\alpha})x. 
\end{equation}
 
\begin{theorem}\label{weight local}
Let $l$ and $m$ be integers with $3\leq l\leq m-1$. If $\tilde{\eta}_0, \tilde{u}_0,\tilde{\phi}_0\in H^{l}$ and $\tilde{\rho}_0\in H^{l-1}$  with  the same assumption in Theorem \ref{simple result}, then the unique solution $(q,u,\phi)$ of the system \eqref{simple model} has 
	\begin{equation}\label{weight conti}
		\tilde{\eta},\tilde{u}\in C([0,T); H^{l}), \quad \tilde{\rho}\in C([0,T); H^{l-1}), \quad \hbox{and} \quad  \tilde{\phi}\in  C([0,T); H^l)\cap L^2((0,T);H^{l+1}).
	\end{equation}
The continuation criterion in \eqref{weight conti} coincides with that in \eqref{continucri} from Theorem \ref{simple result}.
\end{theorem}
\begin{proof}
 Following the proof of Theorem \ref{simple result}, we set $c=1$ and denote the positive constant $C(\alpha, m)$ as $C$.   We also consider the following system, which naturally arises from the system \eqref{local system} for $(\tilde{\eta},\tilde{u},\tilde{\phi})$  by denoting  $\rho_p=\rho-\alpha^\frac{1}{\alpha}$ and $\phi_p=\phi-\alpha^\frac{1}{\alpha}$:
\begin{equation}\label{weight system}
		\begin{split}
		\tilde{\eta}_t+u\tilde{ \eta}_x +\alpha \eta \tilde{ u}_x+\alpha \tilde{ u}_x=&(\alpha+1)u\eta+\alpha u,\\
		\tilde{u}_t+u\tilde{ u}_x+\alpha \eta \tilde{\eta}_x+\alpha \tilde{\eta}_x+\beta\tilde{u}=&u^2+\alpha\eta^2+\alpha\eta+ \tilde{\phi}_x-\phi_p,\\
	\tilde{ \phi}_{xt}-\tilde{ \phi}_{xxx}+\tilde{\phi}_x =&\tilde{c}_\alpha (h+1)(\eta+1)^n\tilde{\eta}_x-\tilde{c}_\alpha (h+1)(\eta+1)^n\eta+\rho_p-2\phi_{xx}. 
		\end{split}
\end{equation}
 where $u,\eta,\rho_p$, and $\phi$ are already constructed in Theorem \ref{simple result} and Corollary \ref{cor local}.
 
  Let us establish a priori estimates for a solution $(\tilde{\eta},\tilde{u},\tilde{\phi})$. We begin with $L^2$ estimates:
\begin{align}\label{local weight estimate1}
		\frac{1}{2}\frac{d}{dt}&\| (\tilde{\eta},\tilde{u})\|^2_{L^2}+\beta\|\tilde{u}\|^2_{L^2}\nonumber\\
		= & -\int (u \tilde{\eta}_x) \tilde\eta\, dx -\alpha\int (\eta \tilde{ u}_x) \tilde\eta\,  dx-\alpha \int \tilde{u}_x\tilde{\eta}\,  dx+\int ((\alpha+1)u\eta+\alpha u)\tilde{\eta}\,  dx\nonumber\\
		&-\int (u \tilde{ u}_x)\tilde{ u}\,  dx-\alpha \int (\eta \tilde{ \eta}_x)\tilde{ u}\,  dx-\alpha \int  \tilde{ \eta}_x\tilde{ u}\,  dx+\int (u^2+\alpha\eta^2+\alpha \eta)\tilde{ u}\, dx \nonumber\\
		&+ \int \tilde{\phi}_x  \tilde{ u}\,  dx- \int \phi_p  \tilde{ u}\,  dx\nonumber\\
		\leq &C (\|(\eta_x,u_x)\|_{L^\infty}+1)\|(\tilde{\eta},\tilde{u})\|^2_{L^2}+C(\|(\eta,u)\|^2_{L^\infty}+1)\|(\eta,u)\|^2_{L^2}+\frac{1}{4}\|\tilde{\phi}_x\|^2_{L^2}+C\|\phi_p\|^2_{L^2}.
\end{align}
 For $H^k$ estimates with $k=1,\dots, l$, we get the following estimates by using Lemma \ref{big lemma} (2): 
\begin{align}\label{local weight estimate2}
		\frac{1}{2}\frac{d}{dt}&\left(\|\partial^k_x(\tilde{\eta},\tilde{u} )\|^2_{L^2}+\|\partial^{k-1}_x\tilde{  \phi}_x\|^2_{L^2}\right)+\|\partial^{k-1}_x\tilde{ \phi}_x\|^2_{H^1}+\beta\|\partial^k_x\tilde{ u}\|^2_{L^2}\nonumber\\
		=& -\int\partial^k_x\tilde{\eta}\cdot[\partial^k_x,u\partial_x]\tilde{\eta}\, dx 
		   -\alpha\int\partial^k_x\tilde{\eta}\cdot[\partial^k_x,\eta\partial_x]\tilde{u}\, dx 
		    -\int\partial^k_x\tilde{u}\cdot[\partial^k_x,u\partial_x] \tilde{u}\, dx\nonumber\\ 
		    &-\alpha\int\partial^k_x\tilde{ u}\cdot[\partial^k_x,\eta\partial_x]\tilde{\eta}\, dx +\frac{1}{2} \int \partial_x u\left(|\partial^k_x\tilde{\eta}|^2+|\partial^k_x\tilde{u}|^2\right)\, dx+\alpha\int \partial_x \eta\partial^k_x\tilde{u}\partial^k_x\tilde{\eta}\, dx\nonumber\\
		    &-\int \partial^k_x \tilde{u}\partial^{k}_x\tilde{\phi}_x\, dx-2\int \partial^{k-1}_x \tilde{\phi}_x\partial^{k-1}_x\phi_{xx}\, dx +\int\partial^k_x((\alpha+1)u\eta+\alpha u)\partial^k_x\tilde{\eta}\,dx\nonumber\\
		    &+\int\partial^k_x(u^2+\alpha\eta^2+\alpha\eta)\partial^k_x\tilde{u}\,dx-\int\partial^k_x \tilde{u}\partial^k_x\phi_p\, dx+\int \partial^{k-1}_x\tilde\phi_x\partial^{k-1}_x \rho_p\, dx\nonumber \\
		    &+\tilde{c}_\alpha\int \partial^{k-1}_x \tilde{\phi}_x\, \partial^{k-1}_x\left((h+1)(\eta+1)^n\tilde{\eta}_x\right)\, dx-\tilde{c}_\alpha\int \partial^{k-1}_x \tilde{\phi}_x\,  \partial^{k-1}_x\left((h+1)(\eta+1)^n\eta\right)\, dx\nonumber\\
	\lesssim & \left(\|\partial_x u\|_{L^\infty}\|\partial^k_x \tilde{\eta}\|_{L^2}+\|\partial_x\tilde{ \eta}\|_{L^\infty}\|\partial^k_x u\|_{L^2}\right)\|\partial^k_x\tilde{ \eta}\|_{L^2}\nonumber\\
	&+\left(\|\partial_x \eta\|_{L^\infty}\|\partial^k_x \tilde{u}\|_{L^2}+\|\partial_x\tilde{ u}\|_{L^\infty}\|\partial^k_x \eta\|_{L^2}\right)\|\partial^k_x\tilde u\|_{L^2}\nonumber\\
		& +\left(\|\partial_x u\|_{L^\infty}\|\partial^k_x \tilde{u}\|_{L^2}+\|\partial_x\tilde{ u}\|_{L^\infty}\|\partial^k_x u\|_{L^2}\right)\|\partial^k_x\tilde{ u}\|_{L^2}\nonumber\\
		&+\left(\|\partial_x \eta\|_{L^\infty}\|\partial^k_x \tilde{\eta}\|_{L^2}+\|\partial_x\tilde{ \eta}\|_{L^\infty}\|\partial^k_x \eta\|_{L^2}\right)\|\partial^k_x\tilde{ u}\|_{L^2}\nonumber\\
		&+(\|\eta\|_{L^\infty}+\|u\|_{L^\infty}+1)\left(\|\partial^k_x(\eta, u)\|^2_{L^2}+\|\partial^k_x(\tilde{\eta},\tilde{u})\|^2_{L^2}\right)\nonumber\\
		& +\left(\|\partial^{k}_x\tilde{\phi}_x\|_{L^2}+\|\partial^k_x\phi_p\|_{L^2}\right)\|\partial^k_x\tilde{u}\|_{L^2}+\left(\|\partial^{k-1}_x\phi_{xx}\|_{L^2}+\|\partial^{k-1}_x\rho_p\|_{L^2}\right)\|\partial^{k-1}_x\tilde{\phi}_x\|_{L^2}\nonumber\\
		&+\tilde{c}_\alpha\int \partial^{k-1}_x \tilde{\phi}_x\, \partial^{k-1}_x\left((h+1)(\eta+1)^n\tilde{\eta}_x\right)\, dx-\tilde{c}_\alpha\int \partial^{k-1}_x \tilde{\phi}_x\,  \partial^{k-1}_x\left((h+1)(\eta+1)^n\eta\right)\, dx\nonumber\\
		\leq &  C\left(\|(\eta, u) \|_{H^l}+1\right)\left(\| \tilde{\eta},\tilde{ u})\|^2_{H^l}+\|(\eta ,u,\rho_p,\phi_p)\|^2_{H^{l+1}}\right)+\frac{1}{4}\|\partial^{k-1}_x\tilde{\phi}_x\|^2_{H^l}\nonumber\\
		&+\tilde{c}_\alpha\int \partial^{k-1}_x \tilde{\phi}_x\, \partial^{k-1}_x\left((h+1)(\eta+1)^n\tilde{\eta}_x\right)\, dx-\tilde{c}_\alpha\int \partial^{k-1}_x \tilde{\phi}_x\,  \partial^{k-1}_x\left((h+1)(\eta+1)^n\eta\right)\, dx.
\end{align}
For the last two term on the right-hand side of the preceding inequality, we estimate these terms by using  Lemma \ref{big lemma} (3):
\begin{align}\label{local weight estimate3}
	&\tilde{c}_\alpha\int \partial^{k-1}_x \tilde{\phi}_x\, \partial^{k-1}_x\left((h+1)(\eta+1)^n\tilde{\eta}_x\right)\, dx-\tilde{c}_\alpha\int \partial^{k-1}_x \tilde{\phi}_x\,  \partial^{k-1}_x\left((h+1)(\eta+1)^n\eta\right)\, dx\nonumber\\
		\lesssim & \big[(\|\eta\|^n_{L^\infty}+1)\left(\|h\|_{L^\infty}\|\partial^{k-1}_x\tilde{\eta}_x\|_{L^2}+\|\tilde{\eta}_x\|_{L^\infty}\|\partial^{k-1}_xh\|_{L^2}\right)\nonumber\\
		&\qquad\qquad\qquad\qquad\qquad\qquad\qquad+\|h\|_{L^\infty}\|\eta\|^{n-1}_{L^\infty}\|\partial_x\tilde{\eta}\|_{L^\infty}\|\partial^{k-1}_x\eta\|_{L^2}\big]\|\partial^{k-1}_x\tilde{\phi}_x\|_{L^2}\nonumber\\
		& +\left[(\|\eta\|^n_{L^\infty}+1)\|\partial^{k-1}_x\tilde{\eta}_x\|_{L^2}+\|\tilde{\eta}_x\|_{L^\infty}(\|\eta\|^{n-1}_{L^\infty}+1)\|\partial^{k-1}_x\eta\|_{L^2}\right]\|\partial^{k-1}_x\tilde{\phi}_x\|_{L^2}\nonumber\\
		&+\big[(\|\eta\|^n_{L^\infty}+1)\left(\|h\|_{L^\infty}\|\partial^{k-1}_x\eta\|_{L^2}+\|\eta\|_{L^\infty}\|\partial^{k-1}_xh\|_{L^2}\right)\nonumber\\
		&\qquad\qquad\qquad\qquad\qquad\qquad\qquad+\|h\|_{L^\infty}\|\eta\|^{n-1}_{L^\infty}\|\eta\|_{L^\infty}\|\partial^{k-1}_x\eta\|_{L^2}\big]\|\partial^{k-1}_x\tilde{\phi}_x\|_{L^2}\nonumber\\
		& +\left[(\|\eta\|^n_{L^\infty}+1)\|\partial^{k-1}_x \eta\|_{L^2}+\|\eta\|_{L^\infty}(\|\eta\|^{n-1}_{L^\infty}+1)\|\partial^{k-1}_x\eta\|_{L^2}\right]\|\partial^{k-1}_x\tilde{\phi}_x\|_{L^2}\nonumber\\
		\lesssim &  (\|\eta\|^n_{H^l}+1)(\|h\|_{H^{l}}+1)\left(\|(\eta,\tilde{\eta})\|^2_{H^{l}}+\|\tilde\phi_x \|^2_{H^{l-1}}\right).
\end{align}
Summing these estimates \eqref{local weight estimate1}, \eqref{local weight estimate2}, and \eqref{local weight estimate3} yields
\begin{equation*}
	\begin{split}
		&\frac{1}{2}\frac{d}{dt}\left(\| (\tilde{\eta},\tilde{ u})\|^2_{H^{l}}+\|\tilde{\phi}_x|^2_{H^{l-1}}\right)+\frac{1}{2}\|\tilde{ \phi}_x\|^2_{H^l}+\beta\|\tilde{u}\|^2_{H^{l}}\\
	\lesssim &  \left(\|(\eta, u,h)\|^{n+1}_{H^l}+1\right)\left(\|(\tilde{\eta},\tilde{ u})\|^2_{H^{l}}+\|\tilde{\phi}_x\|^2_{H^{l-1}}+\|(\eta, u,\rho_p,\phi_p)\|_{H^{l+1}}\right).
	\end{split}
\end{equation*}
For the quantity
\begin{equation*}
	Y_l=\| (\tilde{\eta}, \tilde{u})\|^2_{H^{l}}+\|\tilde{\phi}_x\|^2_{H^{l-1}}+1,
\end{equation*}
we obtain the differential inequality:
\begin{equation*}
	\frac{d}{dt}Y_l\leq C \left(\|(\eta, u,h)\|^{n+1}_{H^l}+1\right)(Y_l+\|(\eta, u,\rho_p,\phi_p)\|_{H^{l+1}}).
\end{equation*}

 We will briefly address the proof of existence, and uniqueness.  The system \eqref{weight system} can be treated as a linear system with a simpler structure than \eqref{local system}. Therefore, existence and uniqueness can be demonstrated in a straightforward manner, following the approach utilized in the proof of Theorem \ref{simple result}. So we can get 
 \begin{equation}\label{weight exis}
 	\begin{split}
 		\tilde{\eta}, \tilde{u}\in C([0,T); H^l) \ \ \ \hbox{and} \ \ \ \tilde{\phi}_x\in C([0,T); H^{l-1})\cap L^2((0, T);H^l).
 	\end{split}
 \end{equation}
Next, the equation about $\tilde{\rho}$ is stated as follows:
\begin{align*}
		\tilde{\rho}_t+u \tilde{\rho}_x=&-\rho_p\tilde{ u}_x-\alpha^\frac{1}{\alpha} \tilde{u}_x+2u\rho_p+\alpha^\frac{1}{\alpha} u.
\end{align*}
From \eqref{local well1} and \eqref{weight exis}, we consider the above equation as a transport equation with force. So we conclude
\begin{equation*} 	\begin{split}
 		\tilde{\eta}, \tilde{u}\in C([0,T); H^l), \qquad \tilde{\rho}\in C([0,T); H^{l-1}),  \ \ \hbox{and} \ \ \tilde{\phi}\in C([0,T); H^l)\cap L^2((0,T);H^{l+1}).
 	\end{split}
 \end{equation*}
Since \eqref{weight conti} also holds as long as  \eqref{local well1} are satisfied, we deduce that the continuation criterion aligns with
the result in Theorem \ref{simple result}.
\end{proof}

\section{Constraints of $W$ and the system of modulation  variables}\label{sec;formmodul}

To capture the emergence of the shock-type singularity  through the stability of $W$ around $\overline{W}$, we impose the following constraints on $W$ at $y=0$, which are consistent with the structure of the blow-up profile.
\begin{equation} \label{W constraint}
	W(0,s)=0,\qquad \partial_yW(0,s)=-1, \qquad  \partial^2_yW(0,s)=0.
\end{equation}
From these constraints and the identity $w=e^{-\frac{s}{2}}W+\kappa$ in \eqref{self similar}, we obtain
\begin{equation}\label{form kappa}
\kappa(t)=w(\xi(t),t).
\end{equation}
Putting $y=0$ and the constraints \eqref{W constraint} into the equation \eqref{no order W}, we have
\begin{equation}\label{modul form 1}
	\begin{split}
		&G^0_W=-F^0_W=-\frac{e^{-\frac{s}{2}}}{1-\dot{\tau}}\left(\frac{2}{1+\alpha} \left(e^{\frac{3}{2}s}\partial_y \Phi^0+\frac{\beta\kappa_0}{2}-\frac{\beta}{2}Z^0\right)-\dot{\kappa}-\frac{\beta\kappa}{1+\alpha}\right),
	\end{split}
\end{equation}
where we denote $f^0(s)=f(0,s)$.

Next, applying  $n=1$, $y=0$ and   the constraints \eqref{W constraint} into the equation \eqref{higher order W}, we obtain
\begin{equation}\label{form tau}
	\begin{split}
		&\dot{\tau}=\frac{\beta e^{-s}}{(1+\alpha)}+\frac{1-\alpha}{1+\alpha }e^\frac{s}{2}\partial_y Z^0+e^{-\frac{s}{2}}\frac{2}{1+\alpha} \left(e^{\frac{3}{2}s}\partial^2_y \Phi^0 -\frac{\beta\partial_y Z^0}{2}\right).
	\end{split}
\end{equation}
Similarly, setting $n=2$, $y=0$ and the constraints \eqref{W constraint} into the equation \eqref{higher order W},  we get 
\begin{equation}\label{modul form 3}
	G^0_W=\frac{F^{(2),0}_W}{\partial^{3}_yW^0}.
\end{equation}
From \eqref{modul form 3} and the definition of $G_W$ in \eqref{g and force}, we arrive at the following equations:
\begin{equation} \label{form xi}
	\dot{\xi}=-(1-\dot{\tau})e^{-\frac{s}{2}}\frac{F^{(2),0}_W}{\partial^{3}_yW^0}+\kappa+\frac{1-\alpha}{1+\alpha}Z^0 -\frac{\kappa_0}{1+\alpha}.
\end{equation}
From \eqref{form tau} and \eqref{form xi}, we derive the system governing the modulation variables \( \tau, \xi : [0, T^*) \to \mathbb{R} \), where both variables are expressed as polynomials and rational functions whose coefficients depend on the derivatives of the self-similar variables \( (W^0, Z^0, \Phi^0) \) and on \eqref{form kappa}, with  \( (\tau(0), \xi(0)) = (\epsilon, 0) \) in \eqref{ini mod}. Since the local well-posedness for $(q,u,\phi)$ in Sobolev space was proved in Theorem \eqref{simple result}, we now demonstrate the local well-posedness of the system consisting of \eqref{form tau} and \eqref{form xi}. Then $\kappa$ in \eqref{form kappa} is determined by $\xi$.

 In addition, we get the following identity from \eqref{modul form 1} and \eqref{modul form 3} to investigate the $L^\infty$-stability of $(1+y^2)^{-\frac{1}{6}}\widetilde{W}$ in Section \ref{sec;interW}:

%From \eqref{form tau} and  \eqref{form xi},  we derive the system governing modulation variables $\tau, \xi : [0,T^*)\to \mathbb{R}$, which are expressed as polynomials and rational functions with coefficients determined by the derivatives of the self-similar variables $(W^0,Z^0,\Phi^0)$ and \eqref{form kappa} with $(\tau(0), \xi(0))=(\epsilon,0)$ in \eqref{ini mod}. Since the local well-posedness for $(q,u,\phi)$ in Sobolev space was proved in Theorem \eqref{simple result}, we now demonstrate the local well-posedness of the system consisting of \eqref{form tau} and \eqref{form xi}. Then $\kappa$ in \eqref{form kappa} is determined by $\xi$. In addition, we get the below identity from \eqref{modul form 1} and \eqref{modul form 3} to investigate the $L^\infty$ stability of $(1+y^2)^\frac{1}{3}\partial_y\widetilde{W}$ in Section \ref{sec;interW}:
\begin{equation}\label{FW}
F_W=-\frac{F^{(2),0}_W}{\partial^{3}_yW^0}+\frac{2 e^{-\frac{s}{2}}}{(1-\dot{\tau})(1+\alpha)} \left(e^{\frac{3}{2}s}\partial_y \Phi-\frac{\beta}{2}Z-e^{\frac{3}{2}s}\partial_y \Phi^0+\frac{\beta}{2}Z^0\right).
\end{equation}

\section{Bootstrap assumptions}\label{boot;assum}

In this section, we list the bootstrap assumptions used for closing the bootstrap arguments.  Combining the initial data \eqref{ini mod}--\eqref{ini 3'} and the continuity criterion \eqref{continucri}, which is connected to Theorems \ref{simple result}, \ref{weight local}, and Corollary \ref{cor local}, allows us to impose these bootstrap assumptions.  Throughout Sections \ref{boot;assum}, \ref{sec;closboot}, and \ref{sec;stabilW}, all assumptions and propositions concerning variables with respect to $s$ hold for all $s\geq -\ln\epsilon$ unless stated otherwise.  To close the bootstrap argument, we fix $M>0$ sufficiently large, and then choose $\epsilon>0$ sufficiently small depending on M. For the modulation variables,  we assume that
\begin{subequations}
\begin{align}
&|\dot \tau(t)| \leq 5Me^{-s}, \qquad \qquad \qquad \   \qquad \qquad \qquad \qquad \quad \ \  |\dot \xi(t)| \leq 8M, \label{mod accel}\\
&|\tau(t)-\epsilon| \leq 10M\epsilon^{2}, \qquad \qquad  T^*\leq 2\epsilon,  \qquad  \hbox{and} \qquad \quad  |\xi(t)| \leq  8M \epsilon, 
\label{mod spped} 
\end{align}
\end{subequations}
for all $t < T^*$. 
From $s=-\ln(\tau-t)$ in \eqref{self trans}, $\tau_0=\epsilon$ in  \eqref{ini mod}, and \eqref{mod accel} we have 
\begin{equation}\label{useep}
e^{-s}\leq \epsilon.
\end{equation} 
Using $e^{-s}\leq\epsilon$ and absorbing M into $e^{-s}$, we obtain:
\begin{equation}\label{rmk;tau}
 \frac{1}{1-\dot{\tau}}\leq \frac{9}{8}.
\end{equation}
Our bootstrap assumptions for the derivatives $ \partial_y^n W$ with $ n = 0, 2, 3, 4$ are stated below. 
\begin{subequations}\label{boot;wholeW}
\begin{align}
&\frac{3}{4}\kappa_0\leq e^{-\frac{s}{2}}W(y,s)+\kappa\leq \frac{5}{4}\kappa_0, \qquad \hbox{for all} \ y\in\mathbb{R},\label{boot;W}\\
& |\partial^3_yW(0,s) - 6|  \leq 1,\label{third W} \\
&\left|\partial^2_y W(y,s)\right|\leq\frac{40|y|}{(1+y^2)^{\frac{1}{2}}},   \qquad \qquad \  \ \hbox{for all} \ y\in\mathbb{R}, \label{boot;2W}\\
&\|\partial_y^3 W(\cdot,s)\|_{L^{\infty}} \leq M^\frac{3}{4}, \qquad \hbox{and} \qquad \|\partial_y^4 W(\cdot,s)\|_{L^\infty} \leq M. \label{boot;dW}
\end{align} 
\end{subequations}
Specifically, \eqref{boot;W} describes the amplitude of $w$, and \eqref{third W} measures the perturbation of $\partial_y^3 W(0,s)$ around $\partial_y^3\overline{W}(0)=6$.\\
From \eqref{kappa constant} and \eqref{boot;W} with $W(0,s)=0$ in \eqref{W constraint}, we observe that
\begin{equation}\label{kappa estimate}
\kappa_0 \geq\frac{ 5(1+\alpha)}{\alpha} > 5 \qquad \hbox{and}\qquad	 \frac{3}{4}\kappa_0\leq  \kappa(t)\leq \frac{5}{4}\kappa_0, \qquad \hbox{for all} \quad t<T^*.
\end{equation}
For the other self-similar variables $Z$ and $\Phi$,  we assume the following estimates:
\begin{align}
	& \|Z(\cdot,s)\|_{L^\infty} \leq 1+8 M\epsilon
 ,\quad \|\partial_y^n Z(\cdot,s)\|_{L^\infty}    \leq 2Me^{- \frac{(7-n)}{4}s}, \quad \hbox{for all}\  n=1,2,3, 4, \label{boot;Z} \\
&\|\partial_y\Phi(\cdot, s)\|_{L^\infty} \leq 2Me^{-\frac{3}{2}s},  \ \ \|\partial^2_y\Phi(\cdot,s)\|_{L^\infty}\leq Me^{-3s},  \ \   \|\partial_y^n \Phi(\cdot,s)\|_{L^{\infty}}\leq 2e^{-\frac{3}{2}s}, \label{boot;phi} \\
&\qquad\qquad\qquad\qquad\qquad\qquad\qquad\qquad\qquad\qquad\qquad\qquad\qquad  \hbox{for all} \quad n=3,4, 5.\nonumber
\end{align}
 To prove the $L^\infty$-stability of $\partial_yW$ in Section \ref{sec;stabilW}, it is necessary to assume some uniform spatial decay rates of perturbations for $\rho$, $u$, and $\phi$:
\begin{equation}\label{eq;weightboot;phi}
	\begin{split}
		&(1+x^2)^\frac{1}{3}\max\left\{ |\rho-\bar\rho|,|u|,|\phi-\bar\phi|,|\phi_x|,|\phi_{xx}|\right\}(x,t)\leq M^\frac{1}{2}, \qquad \hbox{for  all} \ x\in\mathbb{R} \ \hbox{and} \ t<T^*.
	\end{split}
\end{equation} 
It can be inferred from \eqref{eq;weightboot;phi} that  for sufficiently small $\epsilon>0$, the following inequalities hold
\begin{equation}\label{rmk;phi}
 (1+y^2)^\frac{1}{3}\left|\partial_y\Phi\right|(y,s)\leq M^{\frac{1}{2}}e^{-\frac{s}{2}}\quad \hbox{and} \quad (1+y^2)^\frac{1}{3}\left|\partial^2_y\Phi\right|(y,s)\leq M^{\frac{1}{2}}e^{-2s},\ \ \  \hbox{for all} \ \ \  y\in\mathbb{R}.
\end{equation}
Moreover, we also assume
\begin{equation}\label{eq;weightboot;Z}
	\begin{split}
	&(1+y^2)^{\frac{1}{3}}\left|\partial_y Z(y,s)\right|\leq 1, \quad \hbox{for all} \ y\in\mathbb{R}.
	\end{split}
\end{equation} 
Finally, the assumptions for $\partial_yW$ and several derivatives of the perturbation $\widetilde{W}(:=W-\overline{W})$ are as follows.
For $|y|\leq l$:
\begin{equation}\label{eq;tildeW}
	\begin{split}
	&\left|\partial^4_y\widetilde{W}(y,s)\right|\leq 2\epsilon^\frac{1}{5},\qquad \hbox{and}\qquad \left|\partial^n_y\widetilde{W}(y,s)\right|\leq  3\epsilon^\frac{1}{5} l^{4-n},  \qquad n=0,1,2,3.
	\end{split}
\end{equation}
For  $l\leq|y|\leq e^{\frac{3}{2}s}$:
\begin{equation}\label{eq;tildeW;inter}
	\begin{split}
		&(1+y^2)^{-\frac{1}{6}}\left|\widetilde{W}(y,s)\right|\leq \epsilon^\frac{1}{6} \qquad \hbox{and}\qquad (1+y^2)^{\frac{1}{3}}\left|\partial_y\widetilde{W}(y,s)\right|\leq \epsilon^\frac{1}{7}. 
	\end{split}
\end{equation} 
For $e^{\frac{3}{2}s}\leq|y|$:
\begin{equation}\label{eq;W;final}
	(1+y^2)^{\frac{1}{3}}|\partial_y W(y,s)|\leq 1.
\end{equation} 
We begin by establishing several estimates that follow from the bootstrap assumptions.  Combining Lemma \ref{lem;sbur} (4) and the bound \eqref{eq;tildeW;inter} yields
\begin{equation}\label{eq;indfinalrevision}
 (1+y^2)^{\frac{1}{3}}|\partial_yW(y,s)|\leq \frac{1}{2} +\epsilon^\frac{1}{7}	\quad \hbox{for all}\quad |y|\geq 100. 
\end{equation}
For sufficiently small $\epsilon>0$,  applying Lemma \ref{lem;sbur} (2), \eqref{eq;tildeW}, \eqref{eq;tildeW;inter} for $|y|\leq e^{\frac{3}{2}s}$ and \eqref{eq;W;final} for $e^{\frac{3}{2}s}\leq|y|$, we obtain
\begin{equation}\label{eq;Wy}
   (1+y^2)^{\frac{1}{3}}|\partial_yW(y,s)|\leq (1+\epsilon^\frac{1}{7})	\quad \hbox{for all}\quad y\in\mathbb{R}. 
\end{equation}
Integrating the preceding estimate between $0$ and $y$, and using the constraint $W(0,s)=0$ in \eqref{W constraint}, we have
\begin{align*}
|W(y,s)|
&=\left|\int_0^y \partial_y W(\eta,s),d\eta\right|\leq
\left(1+\epsilon^{\frac17}\right)
\int_0^{|y|}\eta^{-\frac23}\,d\eta=3\left(1+\epsilon^{\frac17}\right)|y|^{\frac13}.
\end{align*}
Therefore,
\begin{equation}\label{eq;finalrevision}
|W(y,s)|
\leq
3\left(1+\epsilon^{\frac17}\right)|y|^{\frac13} \quad \hbox{for all}\quad y\in\mathbb{R}.
\end{equation}
Furthermore, since  $M>\kappa_0>0$, subtracting the $\kappa$ estimate in \eqref{kappa estimate} from \eqref{boot;W} and then multiplying by $e^{\frac{s}{2}}$, we get
\begin{equation}\label{wholeW}
\|W(\cdot,s)\|_{L^\infty}\leq M e^{\frac{s}{2}}.
\end{equation}
\section{Closing bootstrap argument I}\label{sec;closboot}

In Section \ref{sec;closboot} and  \ref{sec;stabilW}, we close the bootstrap argument to get blow-up results. Throughout these sections, $f\lesssim g$ represents $f\leq C g$ for some positive constants $C$  which depend on $\alpha$ but are independent of $s,\beta,M,\epsilon$. The parameter $M>0$ is chosen to be sufficiently large, such that it exceeds the combination of the sum and product of $\alpha,\beta$, and $\kappa_0$. The bound $e^{-s}\leq\epsilon$ in \eqref{useep} will be used repeatedly.

In this section, we establish the bootstrap estimates before studying the stability of $\partial_y W$ near $\partial_y \overline{W}$ and the corresponding weighted estimates in the next section. We begin by deriving bounds for the modulation variables.

\subsection{Modulation variables $\tau,\xi$ estimates.}
\begin{proposition}\label{prop;mod}We have 
\begin{align*}
&|\dot \tau(t)| \leq 4M e^{-s}, \qquad \qquad \qquad \   \qquad \qquad \qquad \qquad \quad \ \  |\dot \xi(t)| \leq 4M,\\
&|\tau(t)-\epsilon| \leq 8M\epsilon^2, \qquad \qquad  T^*\leq\frac{3}{2}\epsilon,  \qquad  \hbox{and} \qquad \quad  \ |\xi(t)| \leq  6M \epsilon,
\end{align*}
for all $t<T^*$.
\end{proposition}

\begin{proof} 
Applying the bounds  \eqref{boot;Z} and \eqref{boot;phi} to the form of $\dot{\tau}$ in \eqref{form tau} yields
\begin{equation*}
	\begin{split}
		|\dot{\tau}|\leq & \beta e^{-s}+e^{\frac{s}{2}}|\partial_y Z^0 |+e^{-\frac{s}{2}} \left|2 e^{\frac{3}{2}s}\partial^2_y\Phi^0 -\beta\partial_yZ^0\right|\\
		\leq &Me^{-s}+2 M e^{- s}+2Me^{-2s}+ 2M^2e^{-2s}\leq 4M e^{-s},
	\end{split}
\end{equation*}
where we absorb  $3M$ into $e^{-\frac{s}{2}}$. Thus, from the above estimate,  the bound $T^*\leq2\epsilon$ in \eqref{mod spped}, and the initial data $\tau(0)=\epsilon$ in \eqref{ini mod},  we obtain
\begin{equation*}
	-8M\epsilon^2  \leq \int^{t}_0\dot{\tau}(t')dt'=\tau(t)-\epsilon\leq 8M\epsilon^2, \quad \hbox{for all} \quad t<T^*.
\end{equation*}
In a similar manner, using the identity $\tau(T^*)=T^*$, we get
\begin{equation*}
	(1-4M\epsilon)T^* \leq \int^{T^*}_0 1-\dot\tau(t)dt= \epsilon,
\end{equation*}
which implies
\begin{equation*}
        T^*\leq \frac{\epsilon}{1-4M\epsilon}\leq \frac{3}{2}\epsilon,
\end{equation*}
where $\epsilon>0$ is sufficiently small and depends on $M$. From the $F^{(2)}_W$ estimate in \eqref{damp,vel}, the bounds  \eqref{rmk;tau}, \eqref{third W}, \eqref{boot;Z}, and  \eqref{boot;phi}, we obtain
\begin{align}\label{sec force W}
	\left|\frac{ F^{(2),0}_W}{\partial^3_yW^0}\right|\leq &\frac{1}{5}\left|\frac{2e^{-\frac{s}{2}}}{(1-\dot{\tau})(1+\alpha)}\left( e^{\frac{3}{2}s}\partial^3_y \Phi^0-\frac{\beta}{2}\partial^2_y Z^0\right)+\frac{e^\frac{s}{2}(1-\alpha)}{(1-\dot{\tau})(1+\alpha)}\partial^2_yZ^0\right|\nonumber\\
	\leq &\frac{1}{5}\left(\frac{9}{2} e^{-\frac{s}{2} }+2M^2e^{-\frac{7}{4} s}+\frac{9}{4} Me^{-\frac{3}{4} s}\right)\leq 2 e^{-\frac{s}{2} },
\end{align}
where $3M$ is absorbed into $e^{-\frac{s}{4}}$. Applying the bounds \eqref{mod accel}, \eqref{kappa estimate}, \eqref{boot;Z}, and \eqref{sec force W} with $\xi_0=0$ in \eqref{ini mod} to the form of  $\dot{\xi}$ \eqref{form xi} implies\begin{equation*}
	\begin{split}
		|\dot{\xi}|=&\left|-(1-\dot{\tau})e^{-\frac{s}{2}}\frac{F^{(2),0}_W}{\partial^{3}_yW^0}+\kappa+\frac{1-\alpha}{1+\alpha}Z^0 -\frac{ \kappa_0}{1+\alpha}\right|\leq  2e^{-s}  +M+M+M\leq 4M.
	\end{split}
\end{equation*}
Since $T^*\leq\frac{3}{2}\epsilon$, we obtain
\begin{equation*}
	|\xi(t)|\leq 6M\epsilon.
\end{equation*}
\end{proof}

\subsection{$\Phi$ estimates.}\label{sec;phi} 
\begin{proposition}\label{prop;phi} We have
	\begin{equation*}
	     \|\partial_y\Phi(\cdot, s)\|_{L^\infty}\leq Me^{-\frac{3}{2} s}      \qquad      \hbox{and}       \qquad	\|\partial^n_y\Phi(\cdot,s)\|_{L^\infty}\leq e^{-\frac{3}{2}s},
	\end{equation*}
	for all $n=2,3,4,5$.
\end{proposition}
\begin{remark}
We close the $\partial^2_y\Phi$ estimate  in Proposition \ref{lem;bound phixx}. 
\end{remark}
We state some classical $L^\infty$ properties of the following 1-dimensional heat kernel $H$:
\begin{equation*}
	H_t(x)=(4\pi t)^{-\frac{1}{2}}e^{-\frac{x^2}{4t}}.
\end{equation*}
\begin{lemma}\label{heat kernel}  Let $f:\mathbb{R}\to\mathbb{R}$ be a smooth function. Then, we have
\begin{equation*}
		 \|H_t*f\|_{L^\infty}\leq C\|f\|_{L^\infty}\qquad \hbox{and} \qquad	 \|H_t*\partial_xf\|_{L^\infty}\leq \frac{C}{t^{\frac{1}{2}}} \|f\|_{L^\infty},
\end{equation*}
for all $t>0$ and some constant $C>0$.
\end{lemma}

\begin{lemma}\label{prop;rho}  We have
\begin{equation*}
	\left(\frac{1}{16}\alpha\kappa_0\right)^\frac{1}{\alpha} \leq \rho \leq \left( \alpha\kappa_0\right)^\frac{1}{\alpha}. 
\end{equation*}	
\end{lemma}
\begin{proof} Using the bounds \eqref{boot;W}, \eqref{kappa estimate},  and \eqref{boot;Z} gives
	\begin{equation*}
		\frac{\kappa_0}{8} \leq \frac{3}{4}\kappa_0-\frac{\kappa_0}{5}-8M\epsilon\leq e^{-\frac{s}{2}}W+\kappa-Z\leq \frac{5}{4}\kappa_0+\frac{ \kappa_0}{5}+8M\epsilon\leq 2\kappa_0,
	\end{equation*}
for sufficiently small $\epsilon>0$, dependent on $M>0$. Then we obtain the desired conclusion from the identity of $\sigma$ for the density in \eqref{self similar rho u}.
\end{proof}

\begin{proof}\textbf{(Proof of Proposition \ref{prop;phi}.)} We express derivatives of $\phi$  in the system  \eqref{subrie3} as 
\begin{equation*}
	  \partial^n_x\phi=e^{-\frac{2}{1+\alpha}t}H_{\frac{2}{1+\alpha} t}*\partial^n_x\phi_0+\frac{2}{1+\alpha} \int^t_0 e^{-\frac{2}{1+\alpha}(t-t')}H_{\frac{2}{1+\alpha}(t-t')}*\partial^n_x\rho( t')dt',
\end{equation*}
where we use the identity $\rho=(\alpha q)^\frac{1}{\alpha}$. From Lemma \ref{heat kernel}, we obtain 
\begin{align*}
		 \|\partial^n_x\phi(\cdot,t)\|_{L^\infty}\lesssim &\left\|e^{-\frac{2}{1+\alpha} t}H_{\frac{2}{1+\alpha}t}*\partial^n_x\phi_0\right\|_{L^\infty}+\int^t_0\left\|e^{-\frac{2}{1+\alpha}(t-t')}H_{\frac{2}{1+\alpha}(t-t')}*\partial^n_x\rho(\cdot,t')\right\|_{L^\infty}dt'\\
		\lesssim & \|\partial^{n}_x\phi_0\|_{L^\infty}+ \int^t_0\frac{1}{ (t-t')^\frac{1}{2}}\|\partial^{n-1}_x\rho(\cdot,t')\|_{L^\infty}dt',
\end{align*}
for $n=1,2,3,4,5$.  Using Leibniz rule and Fa\`a di Bruno's formula for the density $\sigma$ in \eqref{self similar rho u}, and then applying Lemma \ref{prop;rho}, the bounds \eqref{boot;wholeW}, \eqref{boot;Z} and \eqref{eq;Wy}, we obtain
\begin{equation*}\label{eq;derden}
\begin{split}
	|\partial^n_y\sigma|&\lesssim \left|\partial^n_y (e^{-\frac{s}{2}}W+\kappa-Z)^\frac{1}{\alpha}\right| \\
	  &\lesssim  \sum_{ \sum k_i i=n, \sum k_i=k}  |e^{-\frac{s}{2}}W+\kappa-Z|^{\frac{1}{\alpha}-k}|e^{-\frac{s}{2}}\partial_yW-\partial_yZ|^{k_1}\cdots\left|e^{-\frac{s}{2}}\partial^n_yW-\partial^n_yZ\right|^{k_n}\\
	  & \leq M,
	\end{split}
\end{equation*}
for sufficiently large $M>0$ and for all $n=0,1,2,3,4$. From the two preceding  estimates  along with \eqref{ini phi Z} and the bound \eqref{mod spped}, we get
\begin{align*}
		e^{\frac{3}{2}s}\|\partial^n_y\Phi(\cdot,s)\|_{L^\infty_y}=&e^{-\frac{3}{2}(n-1)s} \|\partial^n_x\phi(\cdot,t)\|_{L^\infty}	\\
			\ls &e^{-\frac{3}{2}(n-1)s} \|\partial^{n}_x\phi_0\|_{L^\infty}+\int^t_0\frac{1}{(t-t')^\frac{1}{2}}e^{-\frac{3}{2}(n-1)s}\|\partial^{n-1}_x\rho(\cdot,t')\|_{L^\infty}dt'\\
  \ls &e^{-\frac{3}{2}(n-1)s} \|\partial^{n}_x\phi_0\|_{L^\infty}+\int^t_0\frac{1}{(t-t')^\frac{1}{2}}e^{-\frac{3}{2}(n-1)(s-s')}\|\partial^{n-1}_y\sigma(\cdot,s')\|_{L^\infty}dt'\\
     	 \ls &e^{-\frac{3}{2}(n-1)s}+M\int^t_0\frac{1}{(t-t')^\frac{1}{2}}dt'  \ls e^{-\frac{3}{2}(n-1)s}+M\epsilon^\frac{1}{2}, 
\end{align*}
for all $n=1,2,3,4,5$. 
This implies that
	\begin{equation*}
	     \|\partial_y\Phi(\cdot,s)\|_{L^\infty}\leq Me^{-\frac{3}{2} s}      \qquad      \hbox{and}       \qquad	\|\partial^n_y\Phi(\cdot,s)\|_{L^\infty}\leq e^{-\frac{3}{2} s}, \quad \hbox{for}\quad n=2,3,4,5,
	\end{equation*}
where $\epsilon>0$ is sufficiently small.
\end{proof}

\subsection{Bounds for particle trajectories.}\label{subsec;traj}
Let $\psi:\mathbb{R}\times[s_0
,\infty)\to\mathbb{R}$ be the particle trajectory with the velocity $\mathcal{V}$, satisfying the following ODE system:
\begin{equation*}
\frac{d}{ds}\psi(y_0,s)=\mathcal{V}( \psi(y_0,s),s), \qquad \psi(y_0,s_0)=y_0.
\end{equation*}
We denote $\psi_X$ as the particle trajectory corresponding to the velocity $\mathcal{V}_X$ and write $\mathcal{V}\circ\psi(y,s)$ as $\mathcal{V}(\psi(y,s),s)$.  We first present the classical equality for a particle trajectory in \cite{evans2022partial}:
\begin{lemma}\label{max prin1}
Let $f(y,s)$ be a solution of the following damped transport equation with a force
\[
\partial_s f + D f + {\mathcal V} \partial_y f =  F,
\]
 where $D$, ${\mathcal V}$,  and $F$ are smooth functions for all $(y,s )\in\mathbb{R}\times [s_0,\infty)$. Then we have
\begin{equation*}
	f\circ\psi(y,s)=f(y,s_0)e^{-\int^s_{s_0}D\circ\psi(y,s')ds'}+\int^s_{s_0}e^{-\int^s_{s'}D\circ\psi(y,s'')ds''}F\circ\psi(y,s')ds',
\end{equation*}
for all $y\in\mathbb{R}$ and $s\geq s_0$.
\end{lemma}

 We will list several bounds with respect to trajectories $\psi_W$,  and $\psi_Z$.
 \begin{lemma}\label{lowertraj}
For any $y_0\in\mathbb{R}$ with $|y_0|\geq l$ and $s_0\geq -\ln\epsilon$, the trajectory $\psi_W$ with $\psi_W(y_0,s_0)=y_0$  propagates to infinity at an exponential rate, satisfying the following lower bound: 
\begin{align*}
 |\psi_W(y_0,s)|\geq |y_0|e^{\frac{1}{8}(s-s_0)},
\end{align*}
for all $s\geq s_0$.
\end{lemma}
\begin{proof}
Using the mean value theorem, the constraint of $W$ \eqref{W constraint}, and the bound \eqref{eq;Wy} yields 
\begin{align}\label{eq;boundW}
|W(y,s)|\leq |W(0,s)|+ \|\partial_yW(\cdot,s)\|_{L^\infty}|y|\leq \left(1+\epsilon^\frac{1}{7}\right)|y|.
\end{align}
Applying the mean value theorem to $G_W$ in \eqref{g and force}, and subsequently using the equation \eqref{modul form 3}, the bound \eqref{boot;Z}, and the estimate $F^{0,(2)}_W/\partial^3_yW^0$ \eqref{sec force W}, we get
\begin{equation}\label{estimate g}
\begin{split}
	|G_W|= \left|\frac{e^{\frac{s}{2}}}{1-\dot{\tau}}\left(\kappa+\frac{1-\alpha}{1+\alpha}Z-\dot{\xi}-\frac{\kappa_0}{1+\alpha} \right)\right|\leq & \left|\frac{e^{\frac{s}{2}}}{1-\dot{\tau}}\left(\kappa+\frac{1-\alpha}{1+\alpha}Z^0-\dot{\xi}-\frac{\kappa_0}{1+\alpha} \right)\right|+|\partial_yZ||y|\\
	\leq &\left|\frac{F^{0,(2)}_W}{\partial^3_yW^0}\right|+2Me^{-\frac{3}{2} s}|y|\leq 2e^{-\frac{s}{2}}+2Me^{-\frac{3}{2} s}|y|.
\end{split}
\end{equation}
Then, for $|y|\geq l$, using the two preceding estimates and the bound \eqref{rmk;tau}  gives
\begin{align*}
|\mathcal{V}_W(y,s)|=&\left|G_W+\frac{3}{2}y+\frac{1}{1-\dot\tau}W\right |\geq \frac{3}{2}|y|-2e^{-\frac{s}{2}}-2Me^{-\frac{3}{2} s}|y|-\frac{9}{8} \left(1+\epsilon^\frac{1}{7}\right)|y|\geq \frac{1}{8}|y|,
\end{align*}
where $\epsilon>0$ is sufficiently small. From the above bound, we obtain
\begin{align*}
\frac{1}{2}\frac{d}{ds}\psi^2_W(y_0,s)&=\mathcal{V}_W\left(\psi(y_0,s)\right)\psi_W(y_0,s)\geq \frac{1}{8}\psi^2_W(y_0,s).
\end{align*}
Solving the above ODE inequality with $\psi_W(y_0,s_0)=y_0$ yields
\begin{equation*}
|\psi_W(y_0,s)|\geq|y_0|e^{\frac{1}{8}(s-s_0)}.
\end{equation*}
\end{proof}
\begin{lemma} \label{uplw} For any $s_0\geq -\ln\epsilon$ and $y_0\in\mathbb{R}$ with $|y_0|\geq \frac{5}{4}Me^{\frac{s_0}{2}}+e^{\frac{3}{2}s_0} $, the trajectory $\psi_W$ with $\psi_W(y_0,s_0)=y_0$  satisfies the uniform lower bound: 
\begin{equation*}
  |\psi_W(y_0,s)|\geq   e^{\frac{3}{2}s} ,
\end{equation*}
for all $s\geq s_0$
\end{lemma}
\begin{proof}
Differentiating $e^{-\frac{3}{2}s}\psi_Z$ with respect to $s$ yields
\begin{equation*}
	\begin{split}
	       \frac{d}{ds}\left(e^{-\frac{3}{2}s}\psi_W(y_0,s) \right)=e^{-\frac{3}{2}s} \left(\mathcal{V}_W\circ\psi_W-\frac{3}{2}\psi_W\right)=e^{-\frac{3}{2}s} \left(G_W+\frac{W}{1-\dot\tau}\right)\circ\psi_W.
	\end{split}
\end{equation*}
Using an argument analogous to that leading to \eqref{estimate g} and the bound \eqref{boot;Z}, we obtain
\begin{equation*}
\begin{split}
	|G_W|\leq & \left|\frac{e^{\frac{s}{2}}}{1-\dot{\tau}}\left(\kappa+\frac{1-\alpha}{1+\alpha}Z^0-\dot{\xi}-\frac{\kappa_0}{1+\alpha} \right)\right|+|Z-Z^0|\leq \left|\frac{F^{0,(2)}_W}{\partial^3_yW^0}\right|+3\leq 2e^{-\frac{s}{2}}+3.
\end{split}
\end{equation*}
Integrating the above identity from $s_0$ to $s$ and using the bounds \eqref{rmk;tau} and \eqref{wholeW}, together with the preceding estimate, gives
\begin{equation*}
	\begin{split}
	     &  \left| e^{-\frac{3}{2}s}\psi_W(y_0,s)-e^{-\frac{3}{2}s_0} y_0\right|\\
	      \leq&\int^s_{s_0}\left|e^{-\frac{3}{2}s'} \left( G_W+\frac{W}{1-\dot\tau}\right)\circ \psi_W(y_0,s')\right|ds'\leq \int^s_{s_0} 2e^{-2s'}+3e^{-\frac{3}{2} s'}+\frac{9}{8}Me^{-s'}ds'\\
	       \leq & \left[-e^{-2s}+e^{-2s_0}-2\left(e^{-\frac{3}{2}s}-e^{-\frac{3}{2}s_0}\right) -\frac{9}{8}M(e^{-s}-e^{-s_0})\right]\\
	       \leq & \frac{5}{4}M(1-e^{-(s-s_0)})e^{-s_0}\leq  \frac{5}{4}Me^{-s_0}.
	\end{split}
\end{equation*}
for sufficiently small $\epsilon>0$. From this, we obtain
\begin{equation*}
	\begin{split}
	&e^{\frac{3}{2}(s -s_0)}\left[ - \frac{5}{4}Me^{\frac{s_0}{2}}+y_0\right]\leq  \psi_W(y_0,s) \leq e^{\frac{3}{2}(s -s_0)}\left[ \frac{5}{4}Me^{\frac{s_0}{2}}+y_0\right].
	\end{split}
\end{equation*}
Hence, if $y_0\geq\frac{5}{4}Me^{\frac{s_0}{2}}+e^{\frac{3}{2}s_0}$, then
\begin{equation*}
	\begin{split}
	& e^{\frac{3}{2}s } \leq \psi_W(y_0,s).
	\end{split}
\end{equation*}
Similarly, if $y_0\leq -\left(\frac{5}{4}Me^{\frac{s_0}{2}}+e^{\frac{3}{2}s_0}\right)$, then
\begin{equation*}
	\begin{split}
	      \psi_W(y_0,s) \leq& -e^{\frac{3}{2}s }.
	\end{split}
\end{equation*}
Combining the above two estimates completes the proof.
\end{proof}

\begin{lemma}\label{special lemma}For any $y_0\in\mathbb{R}$ and $s_0= -\ln\epsilon$, the trajectory $\psi_Z$ satisfies
  \begin{equation*}
   \int^\infty_{-\ln\epsilon} (1+|\psi_Z(y_0,s')|^2)^{-\frac{1}{3}}ds'\leq C,
  \end{equation*}
  for some constant $C>0$.
\end{lemma}
\begin{proof} For convenience, we denote $\psi(s)=\psi_Z(y_0,s)$ and $f\circ\psi(s)=f(\psi_Z(y_0,s),s)$. First, we prove the following claim:  If  $\psi(s)\leq e^{\frac{s}{2}}$ for all $s\in[-\ln\epsilon,\infty)$, then 
\begin{equation}\label{first ineq}
\frac{d}{ds}\psi(s)\leq -\frac{1}{2}e^{\frac{s}{2}}.
\end{equation}
From the identities of $G_W$ and $G_Z$ in \eqref{g and force}, we get
\begin{equation*}
\begin{split}
 G_Z=&\frac{e^{\frac{s}{2}}}{1-\dot{\tau}}\left(\frac{1-\alpha}{1+\alpha}\kappa-\dot{\xi}-\frac{\kappa_0}{1+\alpha} \right)+\frac{1-\alpha}{(1-\dot\tau)(1+\alpha)}W\\
  =&G_W+\frac{1-\alpha}{(1-\dot\tau)(1+\alpha)}W-\frac{2\alpha e^{\frac{s}{2}}}{(1-\dot\tau)(1+\alpha)}\kappa-\frac{(1-\alpha)e^{\frac{s}{2}}}{(1-\dot{\tau})(1+\alpha)}Z.
 \end{split}
\end{equation*}
Combining the above estimate and the definition of $\mathcal{V}_Z$ in \eqref{damp,vel} gives 
\begin{equation*}
\begin{split}
\mathcal{V}_Z=\frac{3}{2}y+G_Z+ \frac{e^{\frac{s}{2}}}{1-\dot{\tau}}Z=&\left(\frac{3}{2}y+G_W+\frac{(1-\alpha)}{(1-\dot\tau)(1+\alpha)}W\right)+\frac{(-2\alpha e^{\frac{s}{2}}\kappa+2\alpha e^{\frac{s}{2}}Z)}{(1-\dot\tau)(1+\alpha)}\\
=&I_1+I_2.
\end{split}
\end{equation*}
For sufficiently small $\epsilon>0$, using the bounds \eqref{rmk;tau},  \eqref{eq;boundW}, and \eqref{estimate g}, we have 
\begin{align*}
     I_1\circ\psi(s)=&	\left(\frac{3}{2}y+G_W+\frac{(1-\alpha)}{(1-\dot\tau)(1+\alpha)}W\right)\circ\psi(s)\\
     \leq &\frac{3}{2}\psi(s)+2e^{-\frac{s}{2}}+2Me^{-\frac{3}{2} s}|\psi|(s)+(1+\epsilon^{\frac{1}{7}})\frac{9}{8} |\psi|(s)\leq 
    \begin{cases}
2e^{-\frac{s}{2}}, & \psi(s)\leq0,\\[1mm]
3e^{\frac{s}{2}}, & \psi(s)>0,
\end{cases}
\end{align*}
while using the bounds \eqref{rmk;tau}, \eqref{kappa estimate}, and \eqref{boot;Z}, we get
\begin{align*}
     I_2\circ\psi(s)=&\frac{(-2\alpha e^{\frac{s}{2}}\kappa+2\alpha e^{\frac{s}{2}}Z)}{(1-\dot\tau)(1+\alpha)}\leq -\frac{15}{2(1+5M\epsilon)}e^\frac{s}{2}+\frac{9}{4} (1+8M\epsilon)e^{\frac{s}{2}}\leq -\frac{7}{2}e^{\frac{s}{2}}.
\end{align*}
Combining the above two estimates yields
\begin{equation*}
\begin{split}
\frac{d}{ds}\psi(s)=\mathcal{V}_Z\circ\psi(s)=I_1\circ\psi(s)+I_2\circ\psi(s) \leq  -\frac{1}{2}e^{\frac{s}{2}},
\end{split}
\end{equation*}
which proves \eqref{first ineq}. We also show the following claim: There exists  $s_*\geq -\ln\epsilon$ such that
\begin{equation}\label{claim 2}
|\psi(s)|\geq \min\left\{\left|e^{\frac{s}{2}}-e^{\frac{s_*}{2}}\right|, e^{\frac{s}{2}}\right\}, \quad \hbox{for all}\quad s\geq-\ln\epsilon.
\end{equation}
This claim is divided into the following two cases:
\begin{enumerate}[(1)]
\item $\psi(s)> e^{\frac{s}{2}}$ for all $s\in[-\ln\epsilon,\infty)$, or $y_0\leq 0$.
\item  There exists the smallest $s_1\in[-\ln\epsilon, \infty)$ such that $0<\psi(s_1)\leq e^{\frac{s}{2}}$ and $y_0>0$.
\end{enumerate} 
First, we analyze the case (1). If $\psi(s)> e^{\frac{s}{2}}$ for all $s\in[-\ln\epsilon, \infty)$, then it is obvious for \eqref{claim 2}. Otherwise, let $\psi(-\ln\epsilon)=y_0\leq 0$.  Since $y_0\leq 0$, the continuity of $\psi(s)$ together with the monotonicity of $\psi(s)$ implied by \eqref{first ineq} yields
\begin{equation*}
    \psi(s)\leq 0 \quad \text{for all} \quad s\in[-\ln\epsilon,\infty).
\end{equation*}
Therefore, integrating the differential inequality \eqref{first ineq} with respect to $s$, we obtain
\begin{equation*}
\psi(s)\leq -e^{\frac{s}{2}}+\epsilon^{-\frac{1}{2}}\leq 0 \quad  \hbox{for all} \quad s\in[-\ln\epsilon, \infty).
\end{equation*}
Thus, \eqref{claim 2} holds with $s_*=-\ln\epsilon$.

 Next, we analyze the case (2). Applying the inequality \eqref{first ineq} with the assumption of the case (2) implies that   $\frac{d}{ds}\psi(s)\leq -\frac{1}{2}e^{\frac{s}{2}}$ for all $s\geq s_1$. So we conclude that $\psi(s_*)=0$ for some $s_*>s_1$. If $s\geq s_*$, then we obtain
\begin{equation*}
\psi(s)=\int^s_{s_*}\frac{d}{ds}\psi(s')ds'\leq [- e^\frac{s'}{2}]^s_{s_*}=- e^{\frac{s}{2}}+ e^{\frac{s_*}{2}}.  
\end{equation*}
Otherwise, we get
\begin{equation*}
\psi(s)\geq e^{\frac{s_*}{2}}-e^{\frac{s}{2}}.
\end{equation*}
Combining two above estimates yields the inequality \eqref{claim 2}. This concludes the second claim. To complete  the proof, from the inequality \eqref{claim 2} and $\int^\infty_{-\ln\epsilon}e^{-\frac{1}{3}s'}ds'< \infty$,  it suffices to show that 
\begin{equation*}
\mathcal{D}:=\int^\infty_{-\ln\epsilon}(1+|e^{\frac{s'}{2}}-e^{\frac{s_*}{2}}|)^{-\frac{2}{3}}ds'<\infty.
\end{equation*}
By the change of variables $r=e^{\frac{s'}{2}}$, we get
\begin{equation*}
\begin{split}
\mathcal{D}\ls&\int^\infty_{\epsilon^{-\frac{1}{2}}}r^{-1}(1+|r-e^{\frac{s^*}{2}}|)^{-\frac{2}{3}}dr\ls\int^\infty_{\epsilon^{-\frac{1}{2}}}\frac{1}{r^\frac{5}{3}}dr+\int^\infty_{\epsilon^{-\frac{1}{2}}}\frac{1}{(1+|r-e^\frac{s_*}{2}|)^\frac{5}{3}}dr\\
\ls &-\left[r^{-\frac{2}{3}}\right]^\infty_{\epsilon^{-\frac{1}{2}}}+\int^{e^{\frac{s_*}{2}}}_{\epsilon^{-\frac{1}{2}}}\frac{1}{(1+|r-e^{\frac{s_*}{2}}|)^\frac{5}{3}}dr+\int^{\infty}_{e^{\frac{s_*}{2}}}\frac{1}{(1+|r-e^{\frac{s_*}{2}}|)^\frac{5}{3}}dr\\
\ls & \epsilon^{\frac{1}{3}}+\int^{-\epsilon^{-\frac{1}{2}}+e^\frac{s_*}{2}}_{0}\frac{1}{(1+r)^\frac{5}{3}}dr+\int^{\infty}_{0}\frac{1}{(1+r)^\frac{5}{3}}dr\ls  \epsilon^{\frac{1}{3}}+\left[\frac{-1}{(1+r)^\frac{2}{3}}\right]^\infty_0\ls 1.
\end{split}
\end{equation*}

\end{proof}

 \subsection{Z estimates.}\label{Sec z}
\begin{proposition}\label{prop;Z} We have
	\begin{equation*}
		\|Z(\cdot,s)\|_{L^\infty}\leq 1+6M\epsilon \qquad \hbox{and} \qquad  \|\partial^n_y Z(\cdot,s)\|_{L^\infty}\leq Me^{-\frac{(7-n)}{4}s},  	
	\end{equation*}
 for $n=1,2,3,4$.
\end{proposition}

\begin{proof}
$\bullet$ \textbf{($Z$ estimate)} From the equation  \eqref{no order Z}, the damping term $D^{(0)}_Z$ is bounded below from $0$:
\begin{equation*}
	D^{(0)}_Z=\frac{\beta e^{-s}}{(1-\dot{\tau})(1+\alpha)}\geq 0.
\end{equation*}
 From the bounds \eqref{rmk;tau}, \eqref{boot;W}, and  \eqref{boot;phi}, we obtain the forcing term estimates $F_Z$ in \eqref{g and force} 
 \begin{equation*}
	|F_Z|\leq \frac{ e^{-s}}{1-\dot\tau}\left(e^{\frac{3}{2}s}|\partial_y \Phi|+\beta \kappa_0+\beta |e^{-\frac{s}{2}}W+\kappa|\right)\leq \frac{3}{2} e^{-s}(2M+M+M)\leq 6M e^{-s}.
 \end{equation*}
Applying  Lemma \ref{max prin1} for $y_0\in\mathbb{R}$ with the above two estimates and \eqref{ini phi Z} to the equation \eqref{no order Z} yields
 \begin{equation*}
	\begin{split}
		|  Z\circ\psi_Z(y_0,s)|\leq  | Z(y_0, -\ln\epsilon)|+\int^s_{-\ln\epsilon} \left|F_{Z}\circ\psi_Z(y_0,s')\right|ds'\leq & 1+\int^s_{-\ln\epsilon}6Me^{-s'}ds'\\
		\leq &1+6M\epsilon.
	\end{split}
\end{equation*}
$\bullet$ \textbf{($\partial_y Z$ estimate)}
 Multiplying $e^{\frac{3}{2}s}$ into the equation \eqref{higher order Z} and putting $n=1$ yields
\begin{equation}\label{eq;Z1 estimate}
 \partial_s\left( e^{\frac{3}{2}s}\partial_y Z\right)+\left(D^{(1)}_Z-\frac{3}{2} \right)\left( e^{\frac{3}{2}s}\partial_y Z\right)+\mathcal{V}_Z\partial_y\left( e^{\frac{3}{2}s} \partial_yZ\right)=e^{\frac{3}{2}s}\partial_y F_Z.
\end{equation}
Using the bounds \eqref{rmk;tau}, \eqref{boot;Z}, and \eqref{eq;Wy} allows us to compute the  damping term $D^{(1)}_Z-\frac{3}{2}$ as
\begin{equation*}
\begin{split}
\left|D^{(1)}_Z-\frac{3}{2}\right|=&\left|\frac{\beta e^{-s}}{(1+\alpha)(1-\dot{\tau})}+\frac{e^{\frac{s}{2}}}{1-\dot{\tau}}\partial_yZ +\frac{1-\alpha}{(1+\alpha)(1-\dot{\tau})}\partial_yW\right|\\
 \ls & Me^{-s}+Me^{-s}+(1+y^2)^{-\frac{1}{3}} \ls Me^{-s}+(1+y^2)^{-\frac{1}{3}}.
\end{split}
\end{equation*}
Thus, combining the above estimate and Lemma  \ref{special lemma}, we obtain
\begin{equation*}
\exp\left(-\int^s_{-\ln\epsilon}\left(D^{(1)}_Z-\frac{3}{2}\right)\circ\psi_Z(y_0,s') ds'\right)\ls \exp\int^s_{-\ln\epsilon} Me^{-s'}+(1+\psi^2_Z(y_0,s'))^{-\frac{1}{3}}ds'\ls 1,
\end{equation*}
where $\epsilon>0$ is sufficiently small so that $M\epsilon\leq 1$.
From the bounds \eqref{rmk;tau}, \eqref{rmk;phi}, and \eqref{eq;Wy}, we estimate the weighted forcing term $(1+y^2)^\frac{1}{3}\partial_yF_Z$:
\begin{equation*}
\begin{split}
\left|(1+y^2)^{\frac{1}{3}}\partial_yF_Z\right|=&(1+y^2)^{\frac{1}{3}}\frac{2e^{-s}}{(1-\dot{\tau})(1+\alpha)}\left|e^{\frac{3}{2}s} \partial^2_y \Phi-\beta\frac{e^{-\frac{s}{2}}}{2}\partial_y W\right|\\
\ls& M^\frac{1}{2} e^{-\frac{3}{2}s}+M^\frac{1}{2} e^{-\frac{3}{2}s}\ls M^\frac{1}{2}e^{-\frac{3}{2}s}.
\end{split}
\end{equation*}
Thus  applying Lemma \ref{max prin1} for $y_0\in\mathbb{R}$ with Lemma \ref{special lemma}, the two preceding estimates, and \eqref{ini phi Z} to the equation \eqref{eq;Z1 estimate} yields 
 \begin{align}\label{Z1 estimate}
	&e^{\frac{3}{2}s}| \partial_y Z\circ\psi_Z(y_0,s)|\nonumber\\
	\leq & \epsilon^{-\frac{3}{2}} |\partial_y Z(y_0, -\ln\epsilon)|\exp\left(-\int^s_{-\ln\epsilon}\left(D^{(1)}_Z-\frac{3}{2}\right)\circ\psi_Z(y_0,s') ds'\right)\nonumber\\
		&+\int^s_{-\ln\epsilon}e^{\frac{3}{2}s'} \left| \partial_y F_{Z}\circ\psi_Z(y_0,s')\right|\exp\left(-\int^s_{s'}\left(D^{(1)}_Z-\frac{3}{2}\right)\circ\psi_Z(y_0,s'')ds''\right)ds'\nonumber \\
		\ls &  1 +\int^s_{-\ln\epsilon}M^\frac{1}{2} e^{\frac{3}{2}s'}e^{-\frac{3}{2}s'} \left(1+|\psi_Z(y_0,s')|^2\right)^{-\frac{1}{3}} ds'\leq M,
\end{align}	
for sufficiently large $M>0$.

\noindent $\bullet$ \textbf{($\partial^n_y Z$ estimates with $n=2,3,4$)} First, we deal with the $\partial^4_y Z$ estimate. Using the bounds \eqref{rmk;tau}, \eqref{boot;Z}, and \eqref{eq;Wy},   the damping term $D^{(4)}_{Z}$ in \eqref{higher order Z} is estimated as
   \begin{equation*}
   	\begin{split}
   		D^{(4)}_Z=&\frac{\beta e^{-s}}{(1-\dot{\tau})(1+\alpha)}+6+\frac{5}{1-\dot{\tau}}e^\frac{s}{2}\partial_yZ	+4\partial_y G_Z\\
		 \geq & 6-15Me^{\frac{s}{2}}e^{-\frac{3}{2}s}-4\frac{|1-\alpha|}{(1-\dot\tau)(1+\alpha)}|\partial_y W|\geq 1.
   	\end{split}
   \end{equation*}
From the bounds  \eqref{rmk;tau},  \eqref{boot;2W}, \eqref{boot;dW},  \eqref{boot;Z}, and \eqref{boot;phi}, we obtain the estimate for the forcing term $F^{(4)}_Z$ in \eqref{higher order Z}:
   \begin{equation*}
   	\begin{split}
   		|F^{(4)}_Z|\ls & |\partial^4_yF_Z|+|\partial^4_y G_Z|\cdot|\partial_y Z| +\sum^{3}_{k=2}\left(\frac{ e^\frac{s}{2}}{1-\dot\tau}|\partial^k_y Z|+|\partial^k_y G_Z|\right)|\partial^{5-k}_yZ|\\
   		\ls &  e^{-s}\left(e^{\frac{3}{2}s} |\partial^5_{y}\phi |+\beta e^{-\frac{s}{2}}| \partial^4_y W|\right)+|\partial^4_y W|Me^{-\frac{3}{2}s}\\
		& +\sum^{3}_{k=2}\left(Me^{-\frac{5-k}{4}s}+|\partial^k_yW|\right)Me^{-\frac{2+k}{4}s}\\
		\ls & e^{-s}+M^2e^{-\frac{3}{2}s}+M^2e^{-\frac{3}{2}s}+M^2e^{-\frac{7}{4}s}+M^2e^{-s} \ls M^2e^{-s}\leq e^{-\frac{3}{4}s},
   	\end{split}
   \end{equation*}
where $M^2$ is absorbed into $e^{-\frac{s}{4}}$. Applying  Lemma \ref{max prin1} for $y_0\in\mathbb{R}$, together with the above two estimates and \eqref{ini phi Z}, implies that
\begin{align}\label{Z4 estimate}
	| \partial^4_y Z\circ\psi_Z(y_0,s)|\leq &  |\partial^4_y Z(y_0, -\ln\epsilon)|\exp\left(-\int^s_{-\ln\epsilon}D^{(4)}_Z\circ\psi_Z(y_0,s') ds'\right)\nonumber\\
		&+\int^s_{-\ln\epsilon} \left|  F^{(4)}_{Z}\circ\psi_Z(y_0,s')\right|\exp\left(-\int^s_{s'}D^{(4)}_Z\circ\psi_Z(y_0,s'')ds''\right)ds'\nonumber\\
		\leq & \epsilon^6e^{-(s+\ln\epsilon)}+\int^s_{-\ln\epsilon} e^{-\frac{3}{4}s'}e^{-(s-s')} ds'\leq \epsilon^5e^{-s}+4 e^{-\frac{3}{4}s}\leq 5e^{-\frac{3}{4}s }.
\end{align}	
Applying the Gagliardo-Nirenberg-Sobolev inequality and using the bounds \eqref{Z1 estimate}, \eqref{Z4 estimate}, and the initial data $\|\partial^n_yZ(\cdot,-\ln\epsilon)\|_{L^\infty}\leq \epsilon^{\frac{3}{2}n}$ $(n=2,3)$ in \eqref{ini phi Z} gives
\begin{equation*}
\|\partial_y^n Z(\cdot,s)\|_{L^\infty}\lesssim \|\partial_y Z(\cdot,s)\|^\frac{4-n}{3}_{L^\infty}\|\partial^4_y Z(\cdot,s)\|^\frac{n-1}{3}_{L^\infty}\lesssim  M^{\frac{4-n}{3}}e^{-\frac{(7-n)}{4}s}\leq Me^{-\frac{(7-n)}{4}s},
\end{equation*}
for all $n=2,3,$ and sufficiently large $M>0$. This completes the proof.
\end{proof}
\subsection{$\partial^n_yW$ estimates for $n=0,2,3,4$.}\label{sec w}
We close the bootstrap estimates \eqref{boot;wholeW}.
\begin{proposition}\label{prop;Westi}
	We have
	\begin{equation*}
	\frac{13}{16}\kappa_0 \leq e^{-\frac{s}{2}}W+\kappa \leq \frac{19}{16}\kappa_0 .
\end{equation*}
\end{proposition}

\begin{proof}
We recall  from \eqref{subrie1} that
\begin{equation*}
\partial_t w+\frac{\beta}{1+\alpha}w+\left(w+\frac{1-\alpha}{1+\alpha}z-\frac{\kappa_0}{1+\alpha} \right) w_x =\frac{2}{1+\alpha}\left(\phi_x+\frac{\beta\kappa_0}{2}-\frac{\beta}{2}z\right).
\end{equation*}
We define the particle trajectory $\psi_w$ along the flow produced by a velocity $\left(w+\frac{1-\alpha}{1+\alpha}z-\frac{\kappa_0}{1+\alpha} \right)$:
\begin{equation*}
	\frac{d}{dt}w(\psi_w(x,t),t)=\left[-\frac{\beta}{1+\alpha}w +\frac{2}{1+\alpha}\left(\phi_x+\frac{\beta\kappa_0}{2}-\frac{\beta}{2}z\right)\right](\psi_w(x,t),t) \ \  \hbox{with} \ \  \psi_w(x,0)=x.
\end{equation*}
Integrating the above equation in time yields
\begin{equation*}
	w(\psi_w(x,t),t)=w(x,0)+\int^t_0\left[-\frac{\beta}{1+\alpha}w +\frac{2}{1+\alpha}\left(\phi_x+\frac{\beta\kappa_0}{2}-\frac{\beta}{2}z\right)\right](\psi_w(x,t'),t')dt'.
\end{equation*}
From the bounds \eqref{mod spped}, \eqref{boot;W}, \eqref{boot;Z} and \eqref{boot;phi}, we obtain
\begin{equation*}
\begin{split}
	&\int^t_0\left|\frac{\beta}{1+\alpha}w +\frac{2}{1+\alpha}\left(\phi_x+\frac{\beta\kappa_0}{2}-\frac{\beta}{2}z\right)\right|(\psi_w(x,t'),t')dt' \\
	&\qquad  \leq 2\epsilon( 2M\kappa_0+ 4M+M+M+8M^2\epsilon)\leq \frac{\kappa_0}{16},
\end{split}	
\end{equation*}
for sufficiently small $\epsilon>0$. Using the above two estimates with \eqref{ini W}, we get
\begin{equation*}
\begin{split}
	\frac{13}{16}\kappa_0 =\frac{7}{8}\kappa_0-\frac{\kappa_0}{16}  \leq |w(\psi_w(t,x),t)|&\leq \frac{9}{8}\kappa_0 +\frac{\kappa_0}{16}=\frac{19}{16}\kappa_0. 
\end{split}	
\end{equation*}
 By the time continuity of $w(\psi_w(t,x),t)$, we obtain the desired conclusion.
\end{proof}
\begin{proposition}\label{prop;W3} We have
	\begin{equation*}
	\begin{split}
		\left|\partial^3_yW^0(s)-6\right| \leq \epsilon^{\frac{1}{4}}.
	\end{split}
\end{equation*}

\end{proposition}

\begin{proof}
Putting $n=3$ with $y=0$ into the equation \eqref{higher order W} and then applying the constraint of W \eqref{W constraint}, the identity $G^0_W$ \eqref{modul form 3}, and the forcing term $F_W$ in \eqref{g and force} into $F^{(2)}_W$ in \eqref{damp,vel} yield
\begin{equation}\label{ineq;3W}
	\begin{split}
		\partial_s\partial^3_yW^0=&-\left(\frac{\beta e^{-s}}{(1-\dot{\tau})(1+\alpha)}+4\left(1-\frac{1}{1-\dot{\tau}}\right)+\frac{3(1-\alpha)}{(1-\dot\tau)(1+\alpha)} e^\frac{s}{2}\partial_y Z^0\right)\partial^3_yW^0\\
		&-\frac{F^{(2),0}_W}{\partial^3_y W^0}\partial^4_y W^0+\frac{2e^{-\frac{s}{2}}}{(1-\dot{\tau})(1+\alpha)}\left(e^{\frac{3}{2}s} \partial_y^4\Phi^0-\frac{\beta\partial^3_yZ^0}{2}\right)+\frac{(1-\alpha)}{(1-\dot\tau)(1+\alpha)} e^\frac{s}{2}\partial^3_y Z^0.	
	\end{split}
\end{equation}
From the bounds \eqref{mod accel}, \eqref{rmk;tau}, and \eqref{boot;Z}, we obtain
\begin{equation*}
	\begin{split}
		&\left| \frac{\beta e^{-s}}{(1-\dot{\tau})(1+\alpha)}-\frac{4\dot{\tau}}{1-\dot{\tau}}+\frac{3(1-\alpha)}{(1-\dot\tau)(1+\alpha)} e^\frac{s}{2}\partial_y Z^0\right|\ls  M e^{-s}+Me^{- s} +M e^{- s}\ls  M e^{- s}.
	\end{split}
\end{equation*}
Utilizing the bounds \eqref{rmk;tau}, \eqref{boot;dW}, \eqref{boot;Z}, \eqref{boot;phi}, and the $F^{(2),0}_W/\partial^3_yW^0$ estimate \eqref{sec force W}, we get
\begin{equation*}
	\begin{split}
		&\left|-\frac{F_W^{(2),0}}{\partial^3_y W^0}\partial^4_y W^0+\frac{ 2e^{-\frac{s}{2}}}{(1-\dot{\tau})(1+\alpha)}\left(e^{\frac{3}{2}s} \partial_y^4\Phi^0-\frac{\beta\partial^3_yZ^0}{2}\right)+\frac{(1-\alpha)}{(1-\dot\tau)(1+\alpha)} e^\frac{s}{2}\partial^3_y Z^0\right|\\
		\lesssim & Me^{-\frac{s}{2}}+ e^{-\frac{s}{2}}+Me^{-\frac{3}{2}s}+M e^{-\frac{s}{2} }\ls M^2e^{-\frac{ s}{2}}.	
	\end{split}
\end{equation*}
Applying the preceding two estimates to the equation \eqref{ineq;3W} gives
\begin{equation*}
	\begin{split}
		\left|\partial_s\partial^3_yW^0\right|\lesssim M^2e^{-\frac{s}{2} }.
	\end{split}
\end{equation*}
Integrating in time $s$ with the initial data $\partial^3_yW(0,-\ln\epsilon)=6$ in \eqref{ini 1} implies 
\begin{equation*}
	\begin{split}
		\left|\partial^3_yW^0(s)-6\right|\lesssim M^2\int^s_{-\ln\epsilon} e^{-\frac{ s'}{2}}ds'\ls  M^2\epsilon^{\frac{1}{2}} \leq \epsilon^{\frac{1}{4}},
	\end{split}
\end{equation*}
where we absorb $M^2$ into $\epsilon^\frac{1}{4} $.
\end{proof}
\begin{proposition}\label{prop;W2}  We have
\begin{equation*}
	 |\partial^2_yW(y,s)|\leq \frac{33|y|}{(1+y^2)^\frac{1}{2}}.
\end{equation*}	
\end{proposition}
\begin{proof}
We define
\begin{equation*}
	\begin{split}
		R(y,s)=\frac{(1+y^2)^\frac{1}{2}\partial^2_yW(y,s)}{y}.
	\end{split}
\end{equation*}
Putting $n=2$ into the equation \eqref{higher order W} and using the following identity
\begin{equation*}
   \begin{split}
   	\partial_yR=\partial_y\left(\frac{(1+y^2)^\frac{1}{2}\partial^2_yW(y,s)}{y} \right)=&\partial_y^3 W\frac{(1+y^2)^\frac{1}{2}}{y}-\frac{1}{y(1+y^2)}R, 
 \end{split}
\end{equation*}
we obtain
\begin{equation}\label{eq;tildeW2}
	\partial_s R+\left(D^{(2)}_W+\frac{1}{y(1+y^2)}\mathcal{V}_W \right)R+\mathcal{V}_W\partial_yR=\frac{(1+y^2)^\frac{1}{2}}{y}F^{(2)}_W.
\end{equation}
 By the mean value theorem, the constraint $\partial^2_yW(0,s)=0$ in \eqref{W constraint}, and the bounds \eqref{third W} and \eqref{boot;dW}, we get
\begin{equation}\label{W2smalll}
|\partial^2_yW(y,s)|\leq |y|\partial^3_yW(0,s)+\frac{y^2}{2}|\partial^4_yW|\leq 7|y|+\frac{M}{2}y^2\leq \frac{32|y|}{(1+y^2)^\frac{1}{2}},  \ \ \hbox{for all} \ \ |y|\leq l\left(=\frac{1}{M}\right).
\end{equation}
By Lemma \ref{lem;sbur} $(5)$, we obtain the lower bound of the damping term $\left(D^{(2)}_W+\frac{1}{y(1+y^2)}\mathcal{V}_W \right)$:
\begin{align}\label{damp W2}
 & D^{(2)}_W+\frac{1}{y(1+y^2)}\mathcal{V}_W\nonumber\\
  =&\frac{\beta e^{-s}}{(1+\alpha)(1-\dot{\tau})} +\frac{5}{2}+\frac{3}{1-\dot{\tau}}\partial_yW+2\partial_y G_W+\frac{1}{y(1+y^2)}\left(G_W +\frac{3}{2}y+\frac{W}{1-\dot{\tau}}\right)\nonumber\\
 \geq & \frac{5}{2}+ 3\partial_y\overline{W}+\frac{3}{1-\dot{\tau}}\partial_y\widetilde{W}+\frac{3\dot{\tau}}{1-\dot{\tau}}\partial_y\overline{W}+2\partial_y G_W\nonumber\\
 &+\frac{1}{y(1+y^2)}\left(G_W +\frac{3}{2}y+\overline{W}\right)+\frac{1}{y(1+y^2)}\frac{\widetilde{W}}{(1-\dot{\tau})}+\frac{1}{y(1+y^2)}\frac{\dot{\tau}\overline{W}}{(1-\dot{\tau})}\nonumber\\
 \geq & \frac{y^2}{1+y^2}+\frac{3\dot{\tau}}{1-\dot{\tau}}\partial_y\overline{W}+2\frac{(1-\alpha)e^{\frac{s}{2}}\partial_y Z}{(1-\dot{\tau})(1+\alpha)}+\frac{1}{(1+y^2)}\frac{G_W}{y}+\frac{1}{y(1+y^2)}\frac{\widetilde{W}}{(1-\dot\tau)}\nonumber\\
 &+\frac{3}{1-\dot{\tau}}\partial_y\widetilde{W}+\frac{1}{y(1+y^2)}\frac{\dot{\tau}\overline{W}}{(1-\dot{\tau})}.	
\end{align}
Using Lemma \ref{lem;sbur} (2), the bounds \eqref{mod accel}, \eqref{rmk;tau}, \eqref{boot;Z}, and the $G_W$ estimate \eqref{estimate g} gives
\begin{align}\label{dampsub1 W2}
&\left|\frac{3\dot{\tau}}{1-\dot{\tau}}\partial_y\overline{W}+2\frac{(1-\alpha)e^{\frac{s}{2}}\partial_y Z}{(1-\dot{\tau})(1+\alpha)}+\frac{1}{(1+y^2)}\frac{G_W}{y}+\frac{1}{y(1+y^2)}\frac{\dot{\tau}\overline{W}}{(1-\dot{\tau})}\right|\nonumber\\
\ls& Me^{-s}+Me^{- s}+\frac{e^{-\frac{s}{2}}}{|y|}+Me^{-\frac{3}{2}s}+\frac{Me^{-s}}{|y|(1+y^2)^\frac{5}{6}} \nonumber\\
\ls& Me^{-s}+Me^{-\frac{s}{2}}+Me^{-\frac{3}{2}s}+M^2e^{-s}\leq \frac{1}{8M^2}, \qquad \hbox{for all} \qquad |y|\geq l.
\end{align}
for sufficiently small $\epsilon>0$, dependent on $M>0$. From the bounds \eqref{rmk;tau} and \eqref{eq;tildeW;inter}, we obtain 
\begin{equation}\label{dampsub2 W2}
\begin{split}
\left|\frac{3}{1-\dot{\tau}}\partial_y\widetilde{W}+\frac{1}{y(1+y^2)}\frac{\widetilde{W}}{(1-\dot{\tau})}\right|\ls&\frac{ \epsilon^{\frac{1}{7}}}{(1+y^2)^\frac{1}{3} }+\frac{\epsilon^\frac{1}{6}}{|y|(1+y^2)^\frac{5}{6}} \ls \epsilon^\frac{1}{7}+M \epsilon^\frac{1}{6}\leq \frac{1}{8M^2},
\end{split}
\end{equation}
for all $e^{\frac{3}{2}s} \geq|y|\geq l$, and from Lemma \ref{lem;sbur} $(2)$, the bounds \eqref{rmk;tau},  \eqref{boot;W}, and \eqref{eq;W;final}, we get
\begin{equation}\label{dampsub3 W2}
\begin{split}
\left|\frac{3}{1-\dot{\tau}}\partial_y\widetilde{W}+\frac{1}{y(1+y^2)}\frac{\widetilde{W}}{(1-\dot{\tau})}\right|\ls& \frac{1}{(1+y^2)^\frac{1}{3}}+\frac{e^{-\frac{3}{2}s}}{(1+y^2)^\frac{5}{6}} +\frac{e^{\frac{s}{2}}\kappa_0+(1+y^2)^\frac{1}{6}}{|y|(1+y^2)}\\
 \ls& e^{-s}+e^{-4s}+M e^{-4s}+e^{-4s}\leq \frac{1}{8M^2}, \ \   \hbox{for all}\ \  |y|\geq e^{\frac{3}{2}s},
\end{split}
\end{equation}
where $M^4$ is absorbed into $\epsilon^\frac{1}{6}$. Applying the bounds \eqref{dampsub1 W2}, \eqref{dampsub2 W2}, and \eqref{dampsub3 W2} to \eqref{damp W2}, we estimate the damping term \eqref{damp W2} as
\begin{equation*}
\begin{split}
 D^{(2)}_W+\frac{1}{y(1+y^2)}\mathcal{V}_W  \geq & \frac{y^2}{1+y^2}-\frac{1}{4M^2}\geq \frac{1}{2M^2},
  \end{split}	
\end{equation*}
for  all  $|y|\geq l$. Employing the bounds  \eqref{boot;Z}, \eqref{boot;phi}, and \eqref{eq;Wy} to estimate the forcing term $\frac{(1+y^2)^\frac{1}{2}}{y}F^{(2)}_W$ in \eqref{eq;tildeW2} gives
\begin{equation*}
	\begin{split}
		\left|\frac{(1+y^2)^\frac{1}{2}}{y} F^{(2)}_W\right|\leq & \frac{(1+y^2)^\frac{1}{2}}{|y|}\frac{2 e^{-\frac{s}{2}}}{(1-\dot{\tau})(1+\alpha)}\left|e^{\frac{3}{2}s }\partial^3_y\Phi-\frac{\beta}{2}\partial^2_y Z\right|+\left|\frac{(1+y^2)^\frac{1}{2}}{y} \partial^2_yG_W\partial_y W\right|\\
		\ls& Me^{-\frac{s}{2}}+M^2e^{-\frac{7}{4}s}+M^2 e^{-\frac{3}{4} s} \leq e^{-\frac{s}{4} },
	\end{split}
\end{equation*}
for all $|y|\geq l$, where $M^2$ is absorbed into $e^{-\frac{s}{4}}$. By Lemma \ref{max prin1} and \ref{lowertraj}, applied for $|y_0|\geq l$, together with the two preceding estimates, we obtain
\begin{equation} \label{W2bigl}
\begin{split}
	&\left|\left(\frac{(1+y^2)^\frac{1}{2}\partial^2_yW }{y}\right)\circ\psi(y_0,s)\right|\\
	 \leq &\left|\frac{(1+y^2_0)^\frac{1}{2}\partial^2_yW(y_0, -\ln\epsilon) }{y_0} \right|e^{-\frac{1}{2M^2}(s+\ln\epsilon)}+\int^s_{-\ln\epsilon}e^{-\frac{s'}{4}}e^{-\frac{1}{2M^2}(s-s')}ds'\leq32+4\epsilon^\frac{1}{4} \leq 33,
	\end{split}
\end{equation}
where, using the same argument as in \eqref{W2smalll} and the initial data $\|\partial^2_yW(\cdot,-\ln\epsilon)\|_{L^\infty}=1$ from \eqref{ini 2}, we get
	\begin{equation*}
|\partial^2_yW(y,-\ln\epsilon)|\leq\min\left( 6|y|+\frac{31y^2}{2},1\right)\leq \frac{32|y|}{(1+y^2)^\frac{1}{2} }, \quad \hbox{for all} \quad y\in\mathbb{R}.
\end{equation*}
If there exists $ s_1=\min\{s'\in [-\ln\epsilon,\infty): |\psi(y_0,s')|\geq l \}$ for $|y_0|\leq l$, then from the bound \eqref{W2smalll}, we have 
\begin{align}\label{W2big2}
		&\left|\left(\frac{(1+y^2)^\frac{1}{2} \partial^2_yW}{y}\right)\circ\psi(y_0,s)\right|\nonumber\\
		\leq & \left|\left(\frac{(1+y^2)^\frac{1}{2} \partial^2_yW}{y}\right)\circ\psi(y_0,s_1)\right|e^{-\frac{1}{2M^2}(s-s_1)}+\int^s_{s_1}e^{-\frac{ s'}{4}}e^{-\frac{1}{2M^2}(s-s')}ds'\leq 33.
\end{align}
Combining the bounds \eqref{W2smalll}, \eqref{W2bigl}, and \eqref{W2big2}, we conclude
\begin{equation*}
	\|R(\cdot,s)\|_{L^\infty}\leq 33.
\end{equation*}
\end{proof}

\begin{proposition}\label{prop;W24}  We have
\begin{equation*}
	 \|\partial^3_y W(\cdot,s)\|_{L^\infty}\leq\frac{1}{2} M^\frac{3}{4} \qquad \hbox{and} \qquad   \|\partial^4_y W(\cdot,s)\|_{L^\infty}\leq \frac{M}{2}.
\end{equation*}	
\end{proposition}

\begin{proof}
First, we address the $\partial^4_yW$ estimate. Utilizing the bound \eqref{mod accel}, \eqref{rmk;tau},  \eqref{boot;Z}, and \eqref{eq;Wy}, the damping term $D^{(4)}_{W}$ in \eqref{damp,vel} is estimated as
\begin{equation*}
	\begin{split}
		D^{(4)}_W=& \frac{\beta e^{-s}}{(1+\alpha)(1-\dot{\tau})}+\frac{11}{2}+\frac{5}{1-\dot{\tau}}\partial_yW + 4\frac{1-\alpha}{(1+\alpha)(1-\dot\tau)}e^{\frac{s}{2}}\partial_y Z \\
		\geq& \frac{11}{2}-\frac{5}{1-5M\epsilon}\left(1+\epsilon^\frac{1}{7}\right)-9M e^{- s}\geq \frac{1}{4}, 
	\end{split}
\end{equation*}
where $\epsilon>0$ is sufficiently small. Employing the bounds \eqref{rmk;tau}, \eqref{boot;2W}, \eqref{boot;dW}, \eqref{boot;Z}, \eqref{boot;phi}, and \eqref{eq;Wy}, the forcing term $F^{(4)}_W$ in \eqref{damp,vel} is estimated as 
\begin{equation*}
	\begin{split}
		|F^{(4)}_W|\ls & e^{-\frac{s}{2}}\left(e^{\frac{3}{2}s} |\partial^5_y\Phi| +\beta|\partial^4_y Z|\right)+e^{\frac{s}{2}}|\partial^4_y Z||\partial_y W|+|\partial^3_y W||\partial^2_y W|+\sum^3_{k=2} |\partial^k_y Z||\partial^{5-k}_yW|\\
		\ls &e^{-\frac{s}{2}}+ M^2e^{-\frac{5}{4}s}+Me^{-\frac{s}{4}}+M^\frac{3}{4}+ Me^{-s}+M^{\frac{7}{4}}e^{-s}\leq M^\frac{5}{6},
	\end{split}
\end{equation*}
for sufficiently large $M>0$, where $M$ is absorbed into $e^{-s}$. Applying Lemma $\ref{max prin1}$ with the above two estimates and the initial data $\|\partial^4_yW(\cdot,-\ln\epsilon)\|_{L^\infty}\leq1$ from \eqref{ini 2}, we obtain
\begin{align}\label{estimate 4}
		\left|\partial^4_yW\circ\psi_W(y_0,s)\right|\leq & \left|\partial^4_yW(y_0,-\ln\epsilon)\right|e^{-\frac{1}{4}(s+\ln \epsilon)}+M^\frac{5}{6}\int^s_{-\ln\epsilon}e^{-\frac{1}{4}(s-s')}ds' \nonumber\\
		\leq & 1+4M^\frac{5}{6}\leq \frac{M}{2},
\end{align}  
where $2(1+4M^\frac{5}{6})$ is absorbed into $M$. Applying the Gagliardo-Nirenberg-Sobolev inequality with the bound \eqref{estimate 4} and the initial data $\|\partial^3_yW(\cdot,-\ln\epsilon)\|_{L^\infty}\leq7$ from \eqref{ini 2} gives
\begin{equation*}
\|\partial_y^3 W(\cdot,s)\|_{L^\infty}\lesssim\|\partial^2_y W(\cdot,s)\|^\frac{1}{2}_{L^\infty}\|\partial^4_y W(\cdot,s)\|^\frac{1}{2}_{L^\infty}\lesssim M^\frac{1}{2}\leq \frac{1}{2} M^\frac{3}{4},
\end{equation*}
for sufficiently large $M>0$. This completes the proof.
\end{proof}

\section{Closing bootstrap argument II: Stability for $\partial_y W$}\label{sec;stabilW}

In this Section,  we establish the global $L^\infty$-stability $\partial_yW$ near $\partial_y\overline{W}$ through the following estimates, where $\widetilde{W}=W-\overline{W}$:
\begin{equation*}
	\begin{split}
	&\left|\partial_y\widetilde{W}(y,s)\right|\leq 3\epsilon^\frac{1}{5}l^3, \qquad \qquad \ \ \ \ \hbox{for all} \ |y|\leq l,\\
	&\left|\partial_y\widetilde{W}(y,s)\right|\leq \epsilon^\frac{1}{7}(1+y^2)^{-\frac{1}{3}}, \ \ \quad \ \hbox{for all} \ l\leq|y|\leq e^{\frac{3}{2}s},\\
	&\left|\partial_yW(y,s)\right|\leq (1+y^2)^{-\frac{1}{3}}, \qquad \quad \hbox{for all} \ |y|\geq e^{\frac{3}{2}s}.
	\end{split}
\end{equation*}
 To achieve this, we close the bootstrap argument from \eqref{eq;weightboot;phi} to \eqref{eq;W;final}.  
\subsection{Derivatives of $\widetilde{W}$ for $|y|\leq l$.}\label{sec;near 0}
\begin{proposition}
For all $|y|\leq l$, we have
	\begin{align*}
		\left|\partial^4_y\widetilde{W}(y,s)\right|\leq & \frac{3}{2}\epsilon^\frac{1}{5} \qquad \hbox{and} \qquad	\left|\partial^n_y\widetilde{W}(y,s)\right|\leq  2\epsilon^\frac{1}{5} l^{4-n},  
	\end{align*}
	for all $n=0,1,2,3$.
\end{proposition}

\begin{proof}
 First, we address the $\partial^4_y\widetilde{W}$ estimate. Utilizing the bounds \eqref{mod accel}, \eqref{boot;Z}, and  \eqref{eq;Wy}, the damping term $D^{(4)}_{\widetilde{W}}$ in \eqref{damp vel pert} is estimated as
\begin{equation*}
	\begin{split}
		D^{(4)}_{\widetilde{W}}=&\frac{11}{2}+\frac{\beta e^{-s}}{(1-\dot{\tau})(1+\alpha)}+4\partial_yG_W+\frac{5\partial_y\widetilde{W}+5\partial_y\overline W}{1-\dot\tau}\\
		\geq &\frac{11}{2}-9Me^{-s}-5\left(\frac{1}{1-5M\epsilon}\right)\left(1+\epsilon^\frac{1}{7} \right) \geq \frac{1}{4}, 
	\end{split}
\end{equation*}
for all $|y|\leq l$ and  sufficiently small $\epsilon>0$. From Lemma \ref{lem;sbur} (2), (3), the bounds \eqref{mod accel}, \eqref{boot;Z}, \eqref{boot;phi}, \eqref{eq;tildeW}, and the $G_W$ estimate \eqref{estimate g},  the forcing term $F^{(4)}_{\widetilde{W}}$ in \eqref{damp vel pert} is estimated as
\begin{equation*}\label{eq;weight4force}
	\begin{split}
		\left|F^{(4)}_{\widetilde{W}}\right|\lesssim &e^{-\frac{s}{2}}\left|e^{\frac{3}{2}s}\partial^5_y\Phi-\frac{\beta}{2}\partial^4_y Z \right|+\beta e^{-s}\left|\partial^4_y\overline{W}\right|+|G_W||\partial^5_y\overline{W} |+|\dot{\tau}||\partial_y\overline{W}||\partial^4_y\overline{W}|\\
		&+\sum^4_{ k=1}|\partial^{k}_y G_W||\partial^{5-k}_y\overline{W}| +\sum^4_{k=1}\left(|\dot{\tau}|\left|\partial^{4-k}_y\overline{ W}\right|+\left|\partial^{4-k}_y\widetilde{ W}\right|\right)\left| \partial^{k+1}_y\overline{W}\right|\\
		&+\left(|\partial^4_yG_W |+|\partial^4_y\overline{W}|\right)\left|\partial_y\widetilde{W}\right|+\sum^3_{k=2} \left(\left|\partial^{k}_y G_W\right|+\left|\partial^k_y\widetilde{W}\right|+\left|\partial^k_y\overline{W}\right| \right)\left|\partial^{5-k}_y\widetilde{W}\right| \\
		 \lesssim &\left( e^{-\frac{s}{2}}+M^2 e^{-\frac{5}{4} s}\right)+Me^{-s}+\left(e^{-\frac{s}{2}}+Me^{-\frac{3}{2}s}|y|\right)+Me^{-s}\\
		 &+Me^{-\frac{s}{4}}+(M e^{-s}+\epsilon^\frac{1}{5}l) +\left(M e^{-\frac{s}{4} }+1\right)\epsilon^\frac{1}{5}l^3 +(Me^{-\frac{s}{2} }+\epsilon^\frac{1}{5}l + 1)\epsilon^\frac{1}{5}l\leq \frac{1}{8}\epsilon^{\frac{1}{5}},
	\end{split}
\end{equation*}
for all  $|y|\leq l$ and sufficiently small $\epsilon>0$,  depending on $M>0$.
Applying Lemma \ref{max prin1} with the above two estimates yields
\begin{equation*}
	\left|\partial^4_y\widetilde{W}\circ\psi_W(y_0,s)\right|\leq \left|\partial^4_y\widetilde{W}(y_0,-\ln\epsilon)\right|e^{-\frac{1}{4}(s+\ln\epsilon)}+\frac{1}{8}\epsilon^{\frac{1}{5}} \int^s_{-\ln\epsilon}e^{-\frac{1}{4}(s-s')}ds' \leq \frac{3}{2}\epsilon^\frac{ 1}{5},
\end{equation*}
for $|y_0|\leq l$ and as long as $|\psi_W(y_0,s)|\leq l$, where the initial data from \eqref{ini;weight0} are
\begin{equation*}
\begin{split}
	&\left|\partial^4_y\widetilde{W}(y,-\ln\epsilon)\right|\leq \epsilon^\frac{1}{5}, \qquad \hbox{for all}  \ \ |y|\leq l.
\end{split}
\end{equation*}
Due to Lemma \ref{lowertraj}, the initial position $y_0 > l$ satisfies $ \psi_W(y_0,s) > l$, thereby establishing the case $n = 4$.
 By the fundamental theorem of calculus and Proposition \ref{prop;W3}, we obtain
\begin{equation*}
	\begin{split}
		\left|\partial^3_y\widetilde{W}(y,s)\right|	\leq & \left|\partial^3_y\widetilde{W}(0,s)\right|+\int^y_0\left|\partial^4_y\widetilde{W}(y',s)\right|dy'\leq \epsilon^{\frac{1}{4}}+\frac{3}{2} \epsilon^\frac{1}{5} l\leq 2 \epsilon^\frac{1}{5}l ,
	\end{split}
\end{equation*}
for $|y|\leq l$ and sufficiently small $\epsilon>0$.
Using an iterative process from $n=2$ to $n=0$, along with Lemma \ref{lem;sbur} (1), the fundamental theorem of calculus, and the constraints on $W$ in \eqref{W constraint}, we get
\begin{equation*}
	\begin{split}
		\left|\partial^n_y\widetilde{W}(y,s)\right|	\leq &\left|\partial^n_y\widetilde{W}(0,s)\right|+\int^y_0\left|\partial^{n+1}_y\widetilde{W}(y',s)\right|dy'\leq 2  \epsilon^\frac{1}{5} l^{4-n},
	\end{split}
\end{equation*}
for $|y|\leq l$. This completes the proof.
\end{proof}

\subsection{Uniform spatial decay estimates.}\label{sec;uniform}
 In this subsection, we close the uniform spatial decay estimates for  $\rho-\bar{\rho}$, $u$,  $\phi-\bar{\phi}$ and spatial derivatives of $\phi$ in \eqref{eq;weightboot;phi}.

\begin{proposition}\label{lem;rhou} We have
\begin{equation*}
\sup_{t\in[0,T^*), x\in\mathbb{R}}\left(1+x^2\right)^\frac{1}{3}\max\left\{\left|\rho-\bar\rho\right|, |u|  \right\}(x)\leq M^\frac{1}{4}.
\end{equation*}
\end{proposition} 
\begin{proof}
Employing  the bounds \eqref{boot;W}, \eqref{boot;Z}, and Lemma \ref{prop;rho} to estimate \eqref{self similar rho u} gives
\begin{equation}\label{bound rho u}
\sup_{t\in[0,T^*)}\|\rho-\bar\rho(\cdot,t)\|_{L^\infty}\leq \left(\alpha\kappa_0\right)^\frac{1}{\alpha} \quad \hbox{and} \quad \sup_{t\in[0,T^*)}\|u(\cdot,t)\|_{L^\infty}\leq \frac{\kappa_0}{4}+1+8M\epsilon,
\end{equation}
where $\bar\rho=\left(\frac{\alpha\kappa_0}{2}\right)^\frac{1}{\alpha}$, as defined below \eqref{initial weight sobolev}. Thus, it is enough to show that
\begin{equation*}
\left(1+e^{-3s}y^2\right)^\frac{1}{3}\max\left\{\left|\sigma-\bar\rho\right|, |U|  \right\}(y,s)\leq M^\frac{1}{4}, \ \ \text{for all}  \ \ |y|\geq e^{\frac{3}{2}s}.
\end{equation*} 
$\bullet$ \textbf{($\widetilde{\sigma}$ estimate)}
We split the forcing term $F_{\widetilde{\sigma}}$ in \eqref{rho u vel} as
\begin{align*}
F_{\widetilde{\sigma}}=&\left(-e^{\frac{s}{2}}\frac{\dot{\xi}}{(1-\dot\tau)}+\frac{2e^{\frac{s}{2}}}{(1+\alpha)(1-\dot\tau)}U \right)\frac{2e^{-3s}y}{3(1+e^{-3s}y^2)} \widetilde{\sigma} -\sigma(1+e^{-3s}y^2)^\frac{1}{3} \frac{\partial_yW+e^{\frac{s}{2}} \partial_yZ}{(1+\alpha)(1-\dot\tau)}\\
=& I_1+I_2.
\end{align*}
For $|y|\geq e^\frac{3s}{2}$, applying the bounds \eqref{mod accel} and \eqref{eq;weightboot;phi} to $I_1$ gives
\begin{equation*}
\begin{split}
|I_1|=&\left|\left(-e^{\frac{s}{2}}\frac{\dot{\xi}}{(1-\dot\tau)}+\frac{2e^{\frac{s}{2}}}{(1+\alpha)(1-\dot\tau)}U \right)\frac{2e^{-3s}y}{3(1+e^{-3s}y^2)} \widetilde{\sigma}\right|\\
\ls &(Me^{\frac{s}{2}}+M^{\frac{1}{2}}e^{\frac{s}{2}})\frac{M^\frac{1}{2}}{|y|}\ls M^\frac{3}{2}e^{-s}\leq \frac{1}{2}e^{-\frac{s}{4}},
\end{split}
\end{equation*}
and from Lemma \ref{prop;rho}, the bounds \eqref{eq;weightboot;Z}, and \eqref{eq;Wy}, we obtain
\begin{equation*}
	\begin{split}
		|I_2|=&\sigma(1+e^{-3s}y^2)^\frac{1}{3} \frac{|\partial_yW+e^{\frac{s}{2}} \partial_yZ|}{(1+\alpha)(1-\dot\tau)}\ls  Me^{-s}\left(1+y^2\right )^{-\frac{1}{3}}+Me^{-\frac{s}{2} }(1+y^2)^\frac{1}{3}|\partial_yZ|\leq \frac{1}{2}e^{-\frac{s}{4}},
	\end{split}
\end{equation*}
where we absorb $M$ into $e^{-\frac{s}{4}}$. Summing these estimates yields
\begin{equation*}
	\begin{split}
		 |F_{\widetilde\sigma}|\leq &|I_1|+|I_2|\leq e^{-\frac{s}{4}}, \quad \hbox{for} \quad |y|\geq e^\frac{3s}{2}.
	\end{split}
\end{equation*}
First, suppose that $|y_0|\geq \epsilon^{-\frac32}$ and $\psi_{\widetilde{\sigma}}(y_0,-\ln\epsilon)=y_0$. Then we set
\begin{equation*}
s_1:=-\ln\epsilon.
\end{equation*}
Applying  Lemma \ref{max prin1} with  the initial data \eqref{ini final} to the $\widetilde{\sigma}$ equation in \eqref{finalrhou}, we obtain
\begin{equation*}
	\begin{split}
		|\widetilde{\sigma}\circ\psi_{\widetilde{\sigma}}(y_0,s)|\leq &\left|\widetilde{\sigma}(y_0, -\ln\epsilon)\right|+\left|\int^s_{-\ln\epsilon}F_{\widetilde{\sigma}} \circ \psi_{\widetilde{\sigma}}(y_0,s')ds'\right|\leq  1+ \int^s_{-\ln\epsilon} e^{-\frac{s'}{4}} ds'\leq 2,
	\end{split}
\end{equation*}
provided that
\[
|\psi_{\widetilde{\sigma}}(y_0,s')|
\ge
e^{\frac{3}{2} s'}, \qquad s'\in [s_1,s],
\]
for sufficiently small $\epsilon>0$.\\
Next, suppose that $|y_0|<\epsilon^{-\frac32}$ and $\psi_{\widetilde{\sigma}}(y_0,s_0)=y_0$. If the trajectory $\psi_{\widetilde{\sigma}}(y_0,s')$ enters the outer region $\{|y|\ge e^{\frac32 s'}\}$, then we define
\begin{equation*}
s_1:=\min\left\{
s'\in[s_0,\infty):
\left|\psi_{\widetilde{\sigma}}(y_0,s')\right|
\ge e^{\frac32 s'}
\right\}.
\end{equation*}
From Lemma \ref{max prin1}  with the bound \eqref{bound rho u}, we have  
\begin{equation*}
	\begin{split}
		|\widetilde{\sigma}\circ\psi_{\widetilde{\sigma}}(y_0,s)|\leq &\left|\widetilde{\sigma}\circ\psi_{\widetilde{\sigma}}(y_0,s_1)\right|+\left|\int^s_{s_1}F_{\widetilde{\sigma}} \circ \psi_{\widetilde{\sigma}}(y_0,s')ds'\right|
		\leq 2^\frac{1}{3} \left(\alpha\kappa_0\right)^\frac{1}{\alpha}+ \int^s_{s_1} e^{-\frac{s'}{4}} ds'\leq M^\frac{1}{4},
	\end{split}
\end{equation*}
provided that
\[
|\psi_{\widetilde{\sigma}}(y_0,s')|
\ge
e^{\frac{3}{2} s'}, \qquad s'\in [s_1,s],
\]
for sufficiently small $\epsilon>0$. 
Since the trajectories $\psi_{\widetilde{\sigma}}$ arising from the above two cases cover the region $|y|\geq e^{\frac{3}{2}s}$, we conclude the desired $\widetilde{\sigma}$ estimate.

\noindent $\bullet$ \textbf{($\widetilde{U}$ estimate)}
We split the function $F_{\widetilde{U}}$ in \eqref{rho u vel} as
\begin{equation*}
   \begin{split}
	 F_{\widetilde{U}}=&\left(-e^{\frac{s}{2}}\frac{\dot{\xi}}{(1-\dot\tau)}+\frac{2 e^{\frac{s}{2}}}{(1+\alpha)(1-\dot\tau)}U +\frac{2e^{\frac{s}{2}}}{(1+\alpha)(1-\dot{\tau})}\sigma^\alpha\right)\frac{2e^{-3s}y}{3(1+e^{-3s}y^2)} \widetilde{U} \nonumber\\
	&+\frac{ 2 e^{\frac{s}{2}} }{(1+\alpha)(1-\dot\tau)}\sigma^\alpha(1+e^{-3s}y^2)^\frac{1}{3} \partial_yZ+\frac{2 e^{\frac{s}{2}}}{(1+\alpha)(1-\dot\tau)}(1+e^{-3s}y^2)^\frac{1}{3} \partial_y\Phi\\
	=&I_1+I_2+I_3	
   \end{split}	
\end{equation*}
For $|y|\geq e^\frac{3s}{2}$, from Lemma \ref{prop;rho}, the bounds \eqref{mod accel}, \eqref{rmk;tau}, and \eqref{eq;weightboot;phi},   we obtain
\begin{equation*}
	\begin{split}
		|I_1|=&\left|\left(-e^{\frac{s}{2}}\frac{\dot{\xi}}{(1-\dot\tau)}+\frac{2 e^{\frac{s}{2}}}{(1+\alpha)(1-\dot\tau)}U +\frac{2e^{\frac{s}{2}}}{(1+\alpha)(1-\dot{\tau})}\sigma^\alpha\right)\frac{2e^{-3s}y}{3(1+e^{-3s}y^2)} \widetilde{U}\right| \\
		\ls& (M e^{\frac{s}{2}}+M^\frac{1}{2} e^{\frac{s}{2}}+M e^{\frac{s}{2}} )\frac{M^\frac{1}{2}}{|y|} \ls M^{\frac{3}{2}} e^{-s},
	\end{split}
\end{equation*}
For $|y|\geq e^\frac{3s}{2}$, applying Lemma \ref{prop;rho}, the bounds \eqref{rmk;tau}, \eqref{rmk;phi},  and  \eqref{eq;weightboot;Z} leads to
\begin{equation*}
	\begin{split}
		|I_2|=&\frac{ 2 e^{\frac{s}{2}} }{(1+\alpha)(1-\dot\tau)}\sigma^\alpha(1+e^{-3s}y^2)^\frac{1}{3} |\partial_yZ|\ls e^{\frac{s}{2}}M e^{-s}(1+y^2)^\frac{1}{3}|\partial_yZ| \ls M e^{-\frac{s}{2}},
	\end{split}
\end{equation*}
and
\begin{equation*}
	\begin{split}
		|I_3|=\left|\frac{2 e^{\frac{s}{2}}}{(1+\alpha)(1-\dot\tau)}(1+e^{-3s}y^2)^\frac{1}{3} \partial_y\Phi\right|\ls  e^{-\frac{s}{2}}(1+y^2)^{\frac{1}{3}}|\partial_y\Phi|\ls M^\frac{1}{2} e^{-s}.
	\end{split}
\end{equation*}
Summing these estimates yields
\begin{equation}\label{eq;forceU}
	|F_{\widetilde{U}}|\leq |I_1|+|I_2|+|I_3|\ls M^{\frac{3}{2}}e^{- s}+Me^{-\frac{s}{2} }+M^\frac{1}{2} e^{-s}\leq e^{-\frac{s}{4}},
\end{equation}
for $|y|\geq e^{\frac{3}{2}s }$ and  sufficiently small $\epsilon>0$.\\
Applying a similar argument as in the proof of the $\widetilde{\sigma}$ estimate in Proposition
\ref{lem;rhou} to the $\widetilde{U}$ equation in \eqref{finalrhou}, with $s_1$ defined for the trajectory $\psi_{\widetilde{U}}$, together with the forcing estimate \eqref{eq;forceU}, and the bounds \eqref{ini final} and \eqref{bound rho u}, we obtain  
\begin{equation*}
	\begin{split}
		\left|\widetilde{U}\circ\psi_{\widetilde{U}}(y_0,s)\right|\leq &\left|\widetilde U\circ\psi_{\widetilde{U}}(y_0,s_1)\right|\exp\left(-\int^s_{s_1}\frac{\beta e^{-s'}}{1-\dot{\tau}}ds'\right)\\
		&+\left|\int^s_{s_1}F_{\widetilde{U}}\circ \psi_{\widetilde{U}}(y_0,s')\exp\left(-\int^s_{s'}\frac{\beta e^{-s''}}{1-\dot{\tau}}ds''\right)ds'\right|\\
		\leq & \max \left\{1, \frac{\kappa_0}{4}+2\right\}+ \left|\int^s_{s_1} e^{-\frac{s'}{4}} ds'\right|\leq M^\frac{1}{4}, 
    \end{split}
\end{equation*}
provided that
\[
|\psi_{\widetilde{U}}|(y_0,s')|\geq e^{\frac{3}{2}s'}, \qquad s'\in [s_1,s],
\]
for sufficiently small $\epsilon>0$. This completes the proof.

\end{proof}

In the lemma below, we establish the spatial decay rates of the $1D$ heat kernel $H_t$.
\begin{lemma}\label{decay heat}
	For fixed $a,T>0$, let $f:\mathbb{R}\to\mathbb{R}$ be a smooth function such that
	\begin{equation*}
	\begin{split}
	  \sup_{x\in\mathbb{R}}|(1+|x|^2)^a f(x)|=C_0<+\infty.
    \end{split}
	\end{equation*}
	Then, for all $0< t\leq T$  and $x\in\mathbb{R}$, we have
	\begin{equation}\label{weight heat}
	\begin{split}
		&|H_t*f(x)|\leq\frac{C_1}{(1+|x|^2)^a} \qquad \hbox{and} \qquad |H_t*f_x(x)|\leq \frac{ C_1}{\sqrt{t}(1+|x|^2)^a},   \ 
	\end{split}
	\end{equation}
	for some constant $C_1>0$, depending only on $a$, $T$, and $C_0$.
\end{lemma}
\begin{remark} 
The constant $C_1$ increases as the time $T$ grows. So, we choose $C_1>0$ with $T^*\leq1$ in the bootstrap argument.
\end{remark}

\begin{proof}
By Lemma \ref{heat kernel}, it suffices to consider the inequalities \eqref{weight heat} hold for $x\geq 16\sqrt{T} $. First, from the assumption of Lemma \ref{decay heat}, we split the function $|H_t*f|$ as
	\begin{equation*}
	\begin{split}
		|H_t*f(x)|\leq&\left|\int_{\mathbb{R}} \frac{1}{\sqrt{t}}  e^{-\frac{|x-z|^2}{4t}}f(z)dz\right|\\
		\leq &C_0\int_{\mathbb{R}}\frac{2 e^{-z^2}}{(|x-2\sqrt{t}z|^2+1)^a}dz
		=\int^{\frac{x}{4\sqrt{ T}}}_{-\frac{x}{4\sqrt{ T}} }+\int^\infty_{\frac{ x}{4\sqrt{ T}}}+\int^{-\frac{ x}{4\sqrt{T}}}_{-\infty}=I+II+III.
		\end{split}
	\end{equation*}
	 The first term $I$ is bounded by
	\begin{equation*}
		I\lesssim_{a,C_0} \frac{1}{(1+|x|^2)^a}\int^{\frac{x}{4\sqrt{T}} }_{-\frac{x}{4\sqrt{T}}} e^{-z^2}dz\lesssim_{a,C_0}\frac{1}{(1+|x|^2)^a}.
	\end{equation*}
	We estimate the second and third terms, $II$ and $III$, by
	\begin{equation*}
	\begin{split}
		II\lesssim_{T,C_0}& \int^\infty_{\frac{x}{4\sqrt{T}}} e^{-z^2}dz\lesssim_{T,C_0} \int^\infty_{\frac{ x}{4\sqrt{T}}} |z|e^{-z^2}dz\lesssim_{T,C_0}
		 e^{-\frac{x^2}{16 T}}, \\
		  III\ls_{T,C_0} &\int^{-\frac{x}{4\sqrt{T}}}_{-\infty} e^{-z^2}dz\ls_{T,C_0} e^{-\frac{ x^2}{16 T}}.
		  \end{split}
	\end{equation*}
	Summing these estimates yields
	\begin{equation*}
		|H_t*f(x)|\lesssim_{a,T,C_0} \frac{1}{(1+|x|^2)^a}.
	\end{equation*}
Next, from the assumption of Lemma \ref{decay heat} and integration by parts, we split the function $|H_t*f_x|$ as
	\begin{equation*}
	\begin{split}
		|H_t*f_x(x)|\leq&\left|\int_{\mathbb{R}} \frac{1}{\sqrt{t}}  e^{-\frac{|x-z|^2}{4t}}\partial_z f(z)dz\right|=\left|\int_{\mathbb{R}} \frac{(x-z)}{2t\sqrt{t}}  e^{-\frac{|x-z|^2}{4t}}f(z)dz\right|\\
		\leq & \frac{C_0}{\sqrt{t}} \int_{\mathbb{R}} |z| e^{-z^2}\frac{2}{(|x-2\sqrt{t}z|^2+1)^a}dz\\
		=&\int^{\frac{x}{4\sqrt{T}}}_{-\frac{x}{4\sqrt{T}} }+\int^\infty_{\frac{ x}{4\sqrt{T}}}+\int^{-\frac{x}{4\sqrt{T}}}_{-\infty}=I+II+III.
		\end{split}
	\end{equation*}
	 The first term $I$ is bounded by
	\begin{equation*}
		I\lesssim_{a,C_0}\frac{1}{ \sqrt{t}(1+|x|^2)^a} \int^{\frac{ x}{4\sqrt{T}}}_{-\frac{x}{4\sqrt{T}}}|z| e^{-z^2}dz\lesssim_{a,C_0} \frac{1}{\sqrt{t}(1+|x|^2)^a}.
	\end{equation*}
	We estimate the second and third terms, $II$ and $III$, by
	\begin{equation*}
	\begin{split}
		II\lesssim_{T,C_0}& \frac{1}{\sqrt{t}} \int^\infty_{\frac{ x}{4\sqrt{T}}}|z| e^{-z^2}dz\lesssim_{T,C_0} \frac{1}{\sqrt{t}}e^{-\frac{x^2}{16T}}, \\
		III\ls_{T,C_0} &\frac{1}{\sqrt{t}}\int^{-\frac{x}{4\sqrt{T}} }_{-\infty}|z| e^{-z^2}dz\ls_{T,C_0}  \frac{1}{\sqrt{t}}e^{-\frac{x^2}{16T}}.
		\end{split}
	\end{equation*}
Summing these estimates yields 
	\begin{align*}
		&|H_t*f_x(x)|\lesssim_{a,T,C_0} \frac{1}{\sqrt{t} (1+|x|^2)^a}.\\
	\end{align*}
This completes the proof.
	
\end{proof}

\begin{proposition}\label{lem;bound phixx}
We have
	\begin{equation*}
		\sup_{t\in[0,T^*), x\in\mathbb{R}}\left(1+x^2\right)^{\frac{1}{3}}\max\left\{|\phi-\bar{\phi}|, |\partial_x\phi|, |\partial^2_x\phi|, |\partial_t\phi|\right\}(x,t)\leq \frac{M^\frac{1}{2}}{2}. 
	\end{equation*}
\end{proposition}
\begin{remark}By Proposition \ref{lem;bound phixx}, we close the $\partial_y^2\Phi $ estimate in \eqref{boot;phi}$:$
\begin{equation*}
	\left\|\partial_y^2\Phi(\cdot,s)\right\|_{L^\infty}\leq \frac{M^\frac{1}{2}}{2}e^{-3s}.
\end{equation*}
\end{remark}
\begin{proof}\textbf{(Proof of Proposition \ref{lem;bound phixx})}
 We begin by selecting $\epsilon>0$ to be sufficiently small such that $64M\epsilon\leq 1$. 	 We express a solution $\phi$  in the system  \eqref{subrie3} as 
\begin{equation*}
	\begin{split}
	 \phi-\bar\phi=e^{-\frac{2}{1+\alpha} t}H_{\frac{2}{1+\alpha} t}*(\phi_0-\bar\phi)+\frac{2}{1+\alpha} \int^t_0e^{-\frac{2}{1+\alpha} (t-t')}H_{\frac{2}{1+\alpha} (t-t')}*(\rho-\bar\rho)(t') dt'.
	\end{split}
\end{equation*}
Combining Lemma \ref{decay heat} and Proposition \ref{lem;rhou}, together with the initial data $\phi-\bar\phi$ in \eqref{ini final} and the bound $T^*\leq 2\epsilon$ in \eqref{mod spped}, leads to
\begin{equation}\label{weight phi}
	|\phi(x,t)-\bar\phi|\leq \frac{C_1}{(1+|x|^2)^\frac{1}{3} }\left(1+ 2\int^t_0 M^\frac{1}{4} dt'\right)\leq \frac{2C_1}{(1+|x|^2)^\frac{1}{3}}.
\end{equation}
Differentiating the system \eqref{subrie3} for a variable $x$,  we express $\phi_x$ as 
\begin{equation*}
	\begin{split}
		 \phi_x=e^{-\frac{2}{1+\alpha} t}H_{\frac{2}{1+\alpha} t}*(\phi_x)_0+\frac{2}{1+\alpha} \int^t_0e^{-\frac{2}{1+\alpha} (t-t')}H_{\frac{2}{1+\alpha} (t-t')}*\rho_x(t') dt'.
	\end{split}
\end{equation*}
In a similar previous argument with the initial data for $\phi_x$ in \eqref{ini final}, we obtain    
\begin{equation*}
	|\phi_x(x,t)|\leq \frac{C_1}{(1+x^2)^\frac{1}{3}}\left(1+2 \int^t_0 \frac{M^\frac{1}{4} }{\sqrt{ t-t'}}dt'\right)\leq \frac{2C_1}{(1+x^2)^\frac{1}{3} }.
\end{equation*}
Differentiating the system \eqref{subrie3} for time $t$ and using the equation $\rho_t+\frac{2}{1+\alpha} (\rho u)_x=0$,  we express $\phi_t$ as 
\begin{equation*}
	\begin{split}
		 \phi_t=e^{-\frac{2}{1+\alpha} t}H_{\frac{2}{1+\alpha} t}*(\phi_t)_0-\left(\frac{2}{1+\alpha}\right)^2 \int^t_0e^{-\frac{2}{1+\alpha} (t-t')}H_{\frac{2}{1+\alpha} (t-t')}*\partial_x (\rho u)(t') dt'.
	\end{split}
\end{equation*}
In a similar way with the initial data for $(\phi_t)_0\left[=\frac{2}{1+\alpha}(\phi_{xx}-\phi+\rho)_0\right]$ in \eqref{ini final},  we get
\begin{equation}\label{weight phi t}
	|\phi_t(x,t)|\leq \frac{C_1}{(1+|x|^2)^\frac{1}{3} }\left(6+ 2^{\frac{3}{2}}  \int^t_0 \frac{M^\frac{1}{2} }{\sqrt{ t-t'}}dt'\right)\leq \frac{7C_1}{(1+|x|^2)^\frac{1}{3}}.
\end{equation}
Applying  Proposition \ref{lem;rhou}, together with the bounds \eqref{weight phi} and \eqref{weight phi t}, to the equation \eqref{subrie3},
\begin{equation*}
	 \phi_{xx}=\frac{1+\alpha}{2} \phi_t +\phi-\rho,
\end{equation*}
we obtain
\begin{equation*}
	\begin{split}
		|\phi_{xx}(x,t)|\ls \frac{1}{(1+|x|^2)^\frac{1}{3}}(7C_1+2C_1+M^\frac{1}{4})\leq \frac{M^\frac{1}{2}}{2(1+|x|^2)^\frac{1}{3}},
	\end{split}
\end{equation*}
for sufficiently large $M>0$. This completes the proof.
\end{proof}

\subsection{$(1+y^2)^{\frac{1}{3}}\partial_yZ$ estimate.}\label{sec;weightZ}
We consider the following transport-type equation:
\begin{equation*}
	\begin{split}
		\partial_s R+D_RR+\mathcal{V}\partial_y R=F_R.
	\end{split}
\end{equation*}
We also define 
\begin{equation*}
	\begin{split}
		\widetilde{R}:=(1+y^2)^\mu R, \quad \hbox{for} \quad \mu\in\mathbb{R}.
	\end{split}
\end{equation*}
Using the following identity
\begin{equation*}
\partial_y\widetilde{R}=\frac{2\mu y}{(1+y^2)}\widetilde{R}+(1+y^2)^{\mu}\partial_yR,
\end{equation*}
we obtain
\begin{equation}\label{eq;variation}
	\begin{split}
		\partial_s\widetilde{ R}+D_{\widetilde{R}}\widetilde{ R}+\mathcal{V}\partial_y \widetilde{R}=F_{\widetilde{R}},
	\end{split}
\end{equation}
where 
\begin{equation*}
D_{\widetilde{R}}=D_R-\frac{2\mu y}{(1+y^2)}\mathcal{V}\qquad \hbox{and} \qquad F_{\widetilde{R}}=(1+y^2)^\mu F_R.
\end{equation*}
\begin{proposition}
 We have
	\begin{equation*}\label{weight in modul}
	  \left|(1+y^2)^\frac{1}{3}\partial_y Z(y,s)\right|\leq \frac{3}{4}.
\end{equation*} 
\end{proposition}

\begin{proof}
The bound \eqref{boot;Z} implies
  \begin{equation}\label{weight Z}
\left| (1+y^2)^\frac{1}{3} \partial_y Z\right|\leq 2Me^{-\frac{3}{2}s}+2Me^{-\frac{s}{2}}\leq 4Me^{-\frac{s}{2}},  \qquad \text{for all} \ |y|\leq  e^{\frac{3}{2}s}.
\end{equation}
Thus, we need only consider the case  $|y|\geq e^{\frac{3}{2}s}$. Let $\widetilde{R}:=(1+y^2)^{\frac{1}{3}}\partial_yZ$ in \eqref{eq;variation}. Then from \eqref{higher order Z}, the damping term $D_{\widetilde{R}}$ is
 \begin{align}\label{weight dampZ}
   	D_{\widetilde{R}}=&\frac{\beta e^{-s}}{(1-\dot{\tau})(1+\alpha)}+\frac{3}{2}+\frac{1}{1-\dot{\tau}}e^\frac{s}{2}\partial_yZ\nonumber\\
	& +\frac{1-\alpha}{(1-\dot{\tau})(1+\alpha)}\partial_y W-\frac{2y}{3(1+y^2)}\left(G_Z+\frac{3}{2}y+\frac{e^{\frac{s}{2}}}{(1-\dot{\tau})}Z  \right),
   	\end{align}
and the forcing term  $F_{\widetilde{R}}$ is 
 \begin{equation*}
	\begin{split}
		 F_{\widetilde{R}} (y,s)= (1+y^2)^\frac{1}{3}\partial_yF_Z.
	\end{split}
\end{equation*}
Applying the bounds \eqref{mod spped}, \eqref{rmk;tau}, \eqref{kappa estimate}, and \eqref{wholeW} to the $G_Z$ identity  in \eqref{g and force} gives
\begin{align}\label{estimate GZ}
| G_Z|=\left|\frac{e^{\frac{s}{2}}}{1-\dot{\tau}}\left(\frac{1-\alpha}{1+\alpha}\kappa-\dot{\xi}-\frac{\kappa_0}{1+\alpha} \right)+\frac{1-\alpha}{(1-\dot\tau)(1+\alpha)}W\right|\leq&\frac{9}{8} (M+ 8M+M+M)e^{\frac{s}{2}}\nonumber\\
\leq &13Me^\frac{s}{2}.
\end{align}
For all $|y|\geq e^{\frac{3}{2}s}$, the bounds \eqref{rmk;tau}, \eqref{boot;Z}, \eqref{eq;Wy}, and the $G_Z$ estimate \eqref{estimate GZ}  imply
 \begin{align}\label{weight dampZ1}
  & \left|	\frac{1}{1-\dot{\tau}}e^\frac{s}{2}\partial_yZ +\frac{1-\alpha}{(1-\dot{\tau})(1+\alpha)}\partial_yW-\frac{2y^2}{3(y^2+1)}\left(\frac{G_Z}{y}+\frac{e^{\frac{s}{2}}}{(1-\dot{\tau})y} Z \right)\right|\nonumber\\
   \ls &  Me^{- s}+(1+y^2)^{-\frac{1}{3}}+\left(\frac{ M e^{\frac{s}{2}}}{|y|}+\frac{e^\frac{s}{2}}{|y|}\right) \ls  Me^{-s}+e^{-s}+Me^{-s}+e^{-s}\leq CMe^{-s}.
   	\end{align}
 For all $|y|\geq e^{\frac{3}{2}s}$, combining \eqref{weight dampZ} and \eqref{weight dampZ1}, we estimate the damping term $D_{\widetilde{R}}$ as
\begin{equation*}
	D_{\widetilde{R}}\geq \frac{1}{2}+\frac{1}{1+y^2}-CMe^{-s}\geq \frac{1}{4},
\end{equation*}
where $\epsilon>0$ is sufficiently small. From the bounds \eqref{rmk;tau}, \eqref{rmk;phi}, and \eqref{prop;mod},  we obtain 
   \begin{equation*}\label{weight forceZ}
   \begin{split}
   | F_{\widetilde{R}}|=&\left|(1+y^2)^\frac{1}{3}\frac{2e^{-\frac{s}{2}}}{(1-\dot{\tau})(1+\alpha)} \left(e^{\frac{3}{2}s} \partial^2_{y}\Phi-\frac{\beta e^{-\frac{s}{2}} \partial_y W}{2}\right)\right|  \ls M^\frac{1}{2}	e^{-s}+Me^{-s}\leq e^{-\frac{s}{4}},
   \end{split}
   \end{equation*}
  where $M$ is absorbed into $e^{-\frac{s}{4}}$. Using  a similar argument as in the proof of the $\widetilde{\sigma}$ estimate in Proposition \ref{lem;rhou} with $s_1$ defined for the trajectory  $\psi_Z$, together with  the two preceding estimates, the bound \eqref{weight Z} and the initial data  $(1+y^2)^\frac{1}{3}|\partial_y Z(y,-\ln\epsilon)|\leq \frac{1}{2} $ for $y\in\mathbb{R}$ from \eqref{ini 3'}, we have
 \begin{equation*}
	\begin{split}
		\left|\left((1+y^2)^{\frac{1}{3}} \partial_y Z\right)\circ\psi_Z(y_0,s)\right|\leq & \left|\left((1+y^2)^{\frac{1}{3}} \partial_y Z\right)\circ\psi(y_0,s_1)\right|+\int^s_{s_1}\left| F_{\widetilde{R}}(\psi_Z(y,s'),s')\right|ds'\\
		\leq &\max\left\{\frac{1}{2} , 4Me^{-s}\right\}+\int^s_{-\ln\epsilon} e^{-\frac{s'}{4}}ds' \leq \frac{3}{4},
	\end{split}
\end{equation*}
provided that
\[
|\psi_Z|(y_0,s')\geq e^{\frac{3}{2}s'}, \qquad s'\in [s_1,s],
\]
for sufficiently small $\epsilon>0$. This completes the proof.
\end{proof}

\subsection{Weighted estimates for $\widetilde{W}$ and $\partial_y\widetilde{W}$  for $l\leq |y|\leq e^{\frac{3}{2}s}$.}\label{sec;interW}
\begin{proposition}
	 We have
	\begin{equation*}
	  (1+y^2)^{-\frac{1}{6}} \left|\widetilde{W}(y,s)\right|\leq \frac{1}{2} \epsilon^\frac{1}{6}\qquad  \hbox{and} \qquad	(1+y^2)^\frac{1}{3} |\partial_y\widetilde{W}(y,s)|\leq \frac{1}{2} \epsilon^\frac{1}{7},
\end{equation*} 
for all $l\leq |y|\leq e^{\frac{3}{2}s}$.
\end{proposition}
\begin{proof} 
$\bullet$ \textbf{($(1+y^2)^{-\frac{1}{6}}\widetilde{W}$ estimate)} Let $\widetilde{R}:=(1+y^2)^{-\frac{1}{6}}\widetilde{W}$ in \eqref{eq;variation}. Then from \eqref{higher order pertW}, the damping term $D_{\widetilde{R}}$ is
\begin{equation}\label{eq;weightdamp}
	\begin{split}
		D_{\widetilde{R}}=&-\frac{1}{2}+\frac{\beta e^{-s}}{(1-\dot{\tau})(1+\alpha)}+\frac{\partial_y\overline{W}}{1-\dot{\tau}}+\frac{y}{3(1+y^2)}\left(G_W+\frac{3}{2}y+\frac{1}{1-\dot\tau}W \right),
	\end{split}
\end{equation}
and the forcing term $F_{\widetilde{R}}$ is 
\begin{equation}\label{eq;weightforce}
	\begin{split}
		F_{\widetilde{R}}:=(1+y^2)^{-\frac{1}{6}} F_W-(1+y^2)^{-\frac{1}{6}}\frac{\beta e^{-s}}{(1-\dot\tau)(1+\alpha)}\overline{W}-(1+y^2)^{-\frac{1}{6}}\left(G_W+\frac{\dot\tau}{1-\dot\tau} \overline{W}\right)\partial_y\overline{W}.
	\end{split}
\end{equation}
We obtain
\begin{equation*}
	\begin{split}
		\frac{1}{2}-\frac{\beta e^{-s}}{(1-\dot{\tau})(1+\alpha)}-\frac{y^2}{2(1+y^2)}\leq &\frac{1}{2(1+y^2)}.
	\end{split}
\end{equation*}
%\begin{equation*}
%\left|\frac{y}{3(1+y^2)(1-\dot\tau)}W\right|\leq \frac{2}{3}(1+y^2)^{-\frac{1}{2}}\left(\widetilde{W}+\overline{W}\right)\leq +\frac{2}{3}(1+y^2)^{-\frac{1}{2}}(\epsilon^\frac{1}{6}+1)(1+y^2)^\frac{1}{6} \ or\ Me^{\frac{s}{2}}
%\end{equation*}
From Lemma \ref{lem;sbur} (2), the bounds \eqref{rmk;tau}, \eqref{eq;finalrevision}, and the $G_W$ estimate \eqref{estimate g}, we get
\begin{equation*}\label{eq;weightdamp1}
	\begin{split}
		&\left|\frac{\partial_y\overline{W}}{1-\dot{\tau}}+\frac{y}{3(1+y^2)}\left(G_W+\frac{1}{1-\dot\tau}W \right)\right|\\
		\leq &\frac{9}{8} (1+y^2)^{-\frac{1}{3}}+ \frac{y^2}{3(1+y^2)}\left(\frac{2e^{-\frac{s}{2}}}{|y|} +2Me^{-\frac{3}{2} s}\right)+\frac{9}{8}(1+\epsilon^\frac{1}{7})\frac{|y|^{\frac{4}{3}}}{(1+y^2)} \\
		\leq &3(1+y^2)^{-\frac{1}{3}}+\frac{2}{3}M e^{-\frac{s}{2}}+\frac{2}{3}M e^{-\frac{3}{2} s} \\
		\leq & 3(1+y^2)^{-\frac{1}{3}}+2M  e^{-\frac{s}{2} },
	\end{split}
\end{equation*}
for all $ l\leq|y| $. Applying the above two estimates to the damping term \eqref{eq;weightdamp} yields
\begin{equation*}
	\begin{split}
		D_{\widetilde{R}}\geq & -4 (1+y^2)^{-\frac{1}{3}}-2M e^{-\frac{ s}{2}},\quad \hbox{for all} \ \  l\leq|y|.
	\end{split}
\end{equation*}
By Lemma \ref{lowertraj} with $|\psi_W|(y_0,s_0)\geq l$, we get
\begin{align}\label{eq;W;traject}
\int_{s_0}^s \big(1+|\psi_W|^2(y_0,s)\big)^{-\frac{1}{3}}ds'\leq &\int_{s_0}^s \big(1+y_0^2e^{\frac{1}{4} (s'-s_0)}\big)^{-\frac{1}{3}}ds'\leq \int_{s_0}^s \big(1+l^2e^{\frac{1}{4}(s'-s_0)}\big)^{-\frac{1}{3}}ds'\nonumber\\
\leq & \int_{s_0}^{s_0+8\ln{\frac{1}{l}}}1\,ds'+\int_{s_0+8\ln{\frac{1}{l}}}^s l^{-\frac{2}{3}}e^{-\frac{1}{12}(s'-s_0)}\,ds'\nonumber\\
\leq & 8\ln{\frac{1}{l}}+12\leq 9\ln\frac{1}{l}.
\end{align}
Integrating the damping term $D_{\widetilde{R}}$ in time  from $s_0$ to $s$ and using the two preceding estimates  gives 
\begin{equation}\label{eq;weightWdamp}
	\begin{split}
		-\int^s_{s_0}D_{\widetilde{R}}\circ\psi_W(y,s')ds'\leq 36\ln\frac{1}{l}+2M \int^s_{s_0}e^{-\frac{ s'}{2}}ds'\leq 36\ln\frac{1}{l}+4M\epsilon^\frac{1}{2} ,
	\end{split}
\end{equation}
as long as $ l\leq|\psi_W|$. Applying Lemma \ref{lem;sbur} (2), the identity $F_W$ \eqref{FW}, the bounds \eqref{mod accel},  \eqref{rmk;tau}, \eqref{boot;Z}, \eqref{boot;phi}, the $F^{(2),0}_W/\partial^3_yW^0$ estimate \eqref{sec force W},  and the $G_W$ estimate \eqref{estimate g} to the forcing term \eqref{eq;weightforce}, we obtain
\begin{equation*}
	\begin{split}
		|F_{\widetilde{R}}|\leq &(1+y^2)^{-\frac{1}{6}}\left|\frac{F^{(2),0}_W}{\partial^3_yW^0}\right|+\frac{2 e^{-\frac{s}{2}}}{(1-\dot{\tau})(1+\alpha)}(1+y^2)^{-\frac{1}{6}} \left|e^{\frac{3}{2}s}\partial_y \Phi-\frac{\beta}{2}Z-e^{\frac{3}{2}s}\partial_y \Phi^0+\frac{\beta}{2}Z^0\right|\\
		&+(1+y^2)^{-\frac{1}{6}}\frac{\beta e^{-s}}{(1-\dot\tau)(1+\alpha)}|\overline{W}| +(1+y^2)^{-\frac{1}{6}}\left|G_W+\frac{\dot\tau}{1-\dot\tau} \overline{W}\right|\left|\partial_y\overline{W}\right|\\
		\ls& e^{-\frac{s}{2}}+e^{-\frac{s}{2}}(M+M+M+M)+Me^{-s}\\
		&+(1+y^2)^{-\frac{1}{2}}\left(e^{-\frac{s}{2}}+Me^{-\frac{3}{2}s}|y|+Me^{-s}(1+y^2)^{\frac{1}{6}}\right)\\
		\ls &Me^{-\frac{s}{2}} ,
	\end{split}	
\end{equation*}
for all $l\leq |y|$. Thus for $ l\leq|\psi_W|(y,s')$, we get
\begin{equation}\label{eq;weightWforce}
	\begin{split}
		&\int^s_{s_0}\left| F_{\widetilde{R}}\circ\psi_W(y,s')\right|ds'\ls M\int^s_{s_0}  e^{-\frac{ s'}{2}} ds'\ls M\epsilon^\frac{1}{2}.
	\end{split}
\end{equation}
 From the initial data \eqref{ini;weight1}, we have
\begin{equation}\label{self ini}
\left|(1+y^2)^{-\frac{1}{6}}\widetilde{W}(y,-\ln\epsilon)\right|\leq \frac{1}{3}\epsilon^\frac{1}{5},\quad \hbox{for all}\  l\leq|y|\leq  \frac{5}{4}M\epsilon^{-\frac{1}{2}}+\epsilon^{-\frac{3}{2}}.
\end{equation}  Then, we set $s_1=-\ln\epsilon$. Applying Lemmas \ref{max prin1} and \ref{lowertraj} for
\[
l\leq |y_0|
\leq
\frac{5}{4}M\epsilon^{-\frac{1}{2}}
+\epsilon^{-\frac{3}{2}},
\]
together with the damping estimate \eqref{eq;weightWdamp}, the forcing estimate
\eqref{eq;weightWforce}, and the initial data  \eqref{self ini}, we obtain
\begin{equation*}
	\begin{split}
		&\left|\left((1+y^2)^{-\frac{1}{6}} \widetilde{W}\right)\circ\psi_W(y_0,s)\right|\\
		\leq & \left|\left((1+y^2)^{-\frac{1}{6}} \widetilde{W}\right)\circ\psi_W(y_0,s_1)\right|\exp\left(-\int^s_{s_1}D_{\widetilde{R}}\circ \psi_W(y_0,s')ds'\right)\\
		&+\int^s_{s_1}\left| F_{\widetilde{R}}\circ\psi_W(y_0,s')\right|\exp\left(-\int^s_{s'}D_{\widetilde{R}}\circ \psi_W(y_0,s{''})ds''\right)ds'\\
		\ls &\frac{1}{3}\epsilon^\frac{1}{5}l^{-36}e^{4M\epsilon^\frac{1}{2} }+M\epsilon^\frac{1}{2} l^{-36}e^{4M\epsilon^\frac{1}{2} }.
	\end{split}
\end{equation*}
%\begin{equation*}
%	\begin{split}
%		&\left|\left((1+y^2)^{-\frac{1}{6}} \widetilde{W}\right)\circ\psi_W(y_0,s)\right|\\
%		\leq & \left|\left((1+y^2)^{-\frac{1}{6}} \widetilde{W}\right)\circ\psi_W(y_0,-\ln\epsilon)\right|\exp\left(-\int^s_{-\ln\epsilon}D_{\widetilde{R}}\circ \psi_W(y_0,s')ds'\right)\\
%		&+\int^s_{-\ln\epsilon}\left| F_{\widetilde{R}}\circ\psi_W(y_0,s')\right|\exp\left(-\int^s_{s'}D_{\widetilde{R}}\circ \psi_W(y_0,s{''})ds''\right)ds'\\
%		\ls &\frac{1}{3}\epsilon^\frac{1}{5}l^{-36}e^{4M\epsilon^\frac{1}{2} }+M\epsilon^\frac{1}{2} l^{-36}e^{4M\epsilon^\frac{1}{2} }.
%	\end{split}
%\end{equation*}
Next, suppose that $|y_0|<l$ and that there exists
\[
s_1:=
\min\left\{
s'\in[-\ln\epsilon,\infty):
|\psi_W(y_0,s')|
\geq l
\right\}.
\]
Then, from the bound \eqref{eq;tildeW}, we have
%If there exists $ s_1=\min\{s'\in [-\ln\epsilon,\infty): |\psi_W(y_0,s')|\geq l \}$ for $|y_0|\leq l$, then from the bound \eqref{eq;tildeW}, we have 
\begin{equation*}
	\begin{split}
		&\left|\left((1+y^2)^{-\frac{1}{6}} \widetilde{W}\right)\circ\psi_W(y_0,s)\right|\\
		\leq & \left|\left((1+y^2)^{-\frac{1}{6}} \widetilde{W}\right)\circ\psi_W(y_0,s_1)\right|l^{-36}e^{4M\epsilon^\frac{1}{2} }
		+M\epsilon^\frac{1}{2} l^{-36}e^{4M\epsilon^\frac{1}{2} }\\
		\ls &3\epsilon^\frac{1}{5}l^4l^{-36}e^{4M\epsilon^\frac{1}{2} }+M\epsilon^\frac{1}{2} l^{-36}e^{4M\epsilon^\frac{1}{2} }.
	\end{split}
\end{equation*}
Since Lemma \ref{uplw} implies that the trajectories $\psi_W(y_0,s)$ arising from the above two cases cover the region $l\leq|y|\leq e^{\frac{3}{2}s}$, we conclude that for  $l\leq |y|\leq e^{\frac{3}{2}s}$,
\begin{equation*}
	\begin{split}
		&\left|\left((1+y^2)^{-\frac{1}{6}} \widetilde{W}\right)(y,s)\right|
		\ls \max\left\{\frac{1}{3}\epsilon^\frac{1}{5},3\epsilon^\frac{1}{5}l^4\right\}l^{-36}e^{4M\epsilon^\frac{1}{2} }+M\epsilon^\frac{1}{2} l^{-36}e^{4M\epsilon^\frac{1}{2} }\leq \frac{1}{2}\epsilon^{\frac{1}{6}},
	\end{split}
\end{equation*}
for sufficiently small $\epsilon>0$.

\noindent $\bullet$ \textbf{($(1+y^2)^\frac{1}{3}\partial_y\widetilde{W}$ estimate)} Let $\widetilde{R}:=(1+y^2)^{\frac{1}{3}} \partial_y\widetilde{W}$ in \eqref{eq;variation}. Then from \eqref{higher order pertW}, we obtain the damping term $D_{\widetilde{R}}$ and the forcing term $F_{\widetilde{R}}$ as follows:
\begin{equation}\label{weightderidamp}
\begin{split}
D_{\widetilde{R}}:=&1+\frac{\beta e^{-s}}{(1-\dot{\tau})(1+\alpha)}+\frac{\partial_y\widetilde{W}+2\partial_y\overline{W}}{1-\dot{\tau}} +\frac{e^{\frac{s}{2}}(1-\alpha)}{(1-\dot{\tau})(1+\alpha)}\partial_yZ \\
&-\frac{2y}{3(1+y^2)}\left( G_W+\frac{3}{2}y+\frac{W}{1-\dot{\tau}}\right) \\
\geq&\frac{1}{1+y^2}+\frac{\partial_y W+\partial_y\overline{W}}{1-\dot{\tau}} +\frac{e^{\frac{s}{2}}(1-\alpha)}{(1-\dot{\tau})(1+\alpha)}\partial_yZ -\frac{2y^2}{3(1+y^2)}\left(\frac{ G_W}{y}+\frac{W}{y(1-\dot{\tau})}\right),
\end{split}
\end{equation}
and
\begin{align}\label{weightderiforc}
	F_{\widetilde{R}}=&(1+y^2)^\frac{1}{3}\frac{2 e^{-\frac{s}{2}} }{(1-\dot\tau)(1+\alpha)}\left(e^{\frac{3}{2}s}\partial^2_y\Phi-\frac{\beta\partial_y Z}{2}\right)-\frac{\beta e^{-s}}{(1-\dot{\tau})(1+\alpha)}(1+y^2)^\frac{1}{3} \partial_y\overline{W}\nonumber \\
	&-(1+y^2)^\frac{1}{3}G_W\partial^2_y\overline{W}-(1+y^2 )^\frac{1}{3} \left(\frac{\dot{\tau}\partial_y\overline{W}}{1-\dot{\tau}}+\partial_y G_W\right)\partial_y\overline{W}\nonumber\\
	&-(1+y^2)^\frac{1}{3} \left(\frac{\dot{\tau}\overline{W}+\widetilde{ W}}{1-\dot{\tau}}\right)\partial^2_y\overline{W}.
\end{align}
For $l\leq |y|$, using Lemma \ref{lem;sbur} (2), the bounds \eqref{rmk;tau}, \eqref{boot;Z},  \eqref{eq;Wy}, \eqref{eq;finalrevision}, and the $G_W$ estimate  \eqref{estimate g} yields 
\begin{align*}
&\left|\frac{\partial_y W+\partial_y\overline{W}}{1-\dot{\tau}} +\frac{e^{\frac{s}{2}}(1-\alpha)}{(1-\dot{\tau})(1+\alpha)}\partial_yZ -\frac{2y^2}{3(1+y^2)}\left(\frac{ G_W}{y}+\frac{W}{y(1-\dot{\tau})}\right) \right|\\
\leq &\frac{9}{8}\left(\epsilon^{\frac{1}{7}}+2\right)(1+y^2)^{-\frac{1}{3}}+\frac{9}{4}  M e^{- s }+\frac{2}{3}\left(\frac{2e^{-\frac{s}{2}}}{|y|}+2Me^{-\frac{3}{2} s}\right)+\frac{9}{4}(1+\epsilon^\frac{1}{7})\frac{|y|^{\frac{4}{3}}}{(1+y^2)}    \\
\leq & 5(1+y^2)^{-\frac{1}{3}}+3Me^{-\frac{s}{2} s}.
\end{align*}
Applying the above estimate to the damping term \eqref{weightderidamp} for $ l\leq |y|$, we get
\begin{align*}
 D_{\widetilde{R}}\geq& -5(1+y^2)^{-\frac{1}{3}}-3Me^{-\frac{ s}{2}}.
  \end{align*}
Combining the above estimate and \eqref{eq;W;traject} for $ l\leq |\psi_W|$ leads to 
\begin{equation}\label{eq;weight;damping}
	\begin{split}
		-\int^s_{s_0}D_{\widetilde{R}}\circ\psi_W(y_0,s')ds'\leq 45\ln\frac{1}{l}+3M\int^s_{s_0}e^{-\frac{ s'}{2}}ds'\leq 45\ln\frac{1}{l}+6M\epsilon^\frac{1}{2}  .
	\end{split}
\end{equation}
%For $l\leq|y|\leq e^{\frac{3}{2}s}$
%\begin{equation*}
%	(1+y^2)^\frac{1}{3} \left(\frac{\dot{\tau}\overline{W}+\widetilde{W}}{1-\dot{\tau}}\right)\partial^2_y\overline{W}\lesssim (1+y^2)^{-\frac{1}{2}}\left(M e^{-s}(1+y^2)^\frac{1}{6}+\epsilon^\frac{1}{6}(1+y^2)^\frac{1}{6} \right)\lesssim \epsilon^\frac{1}{6}(1+y^2)^{-\frac{1}{3}}
%\end{equation*}
%For $|y|\geq e^{\frac{3}{2}s}$
%\begin{equation*}
%	(1+y^2)^\frac{1}{3} \left(\frac{W}{1-\dot{\tau}}-\overline{W}\right)\partial^2_y\overline{W}\lesssim (1+y^2)^{-\frac{1}{2}}\left(Me^{\frac{s}{2}}+(1+y^2)^\frac{1}{6} \right)\ls Me^{-s}
%\end{equation*}
For $l\leq|y|\leq e^{\frac{3}{2}s}$, from Lemma \ref{lem;sbur} (2), the bounds \eqref{mod accel}, \eqref{rmk;tau}, and \eqref{eq;tildeW;inter}, we obtain
\begin{equation*}
	(1+y^2)^{-\frac{1}{2}} \left|\frac{\dot{\tau}\overline{W}+\widetilde{W}}{1-\dot{\tau}}\right|\lesssim (1+y^2)^{-\frac{1}{2}}\left(M e^{-s}(1+y^2)^\frac{1}{6}+\epsilon^\frac{1}{6}(1+y^2)^\frac{1}{6} \right)\lesssim \epsilon^\frac{1}{6}(1+y^2)^{-\frac{1}{3}}.
\end{equation*}
For $|y|\geq e^{\frac{3}{2}s}$, from Lemma \ref{lem;sbur} (2), the bounds  \eqref{rmk;tau},  and \eqref{wholeW}, we have
\begin{equation*}
	(1+y^2)^{-\frac{1}{2}} \left|\frac{\dot{\tau}\overline{W}+\widetilde{W}}{1-\dot{\tau}}\right|=(1+y^2)^{-\frac{1}{2}} \left|\frac{W}{1-\dot{\tau}}-\overline{W}\right|\lesssim (1+y^2)^{-\frac{1}{2}}\left(Me^{\frac{s}{2}}+(1+y^2)^\frac{1}{6} \right)\ls Me^{-s}.
\end{equation*}
Applying Lemma \ref{lem;sbur} (2), the bounds \eqref{rmk;tau}, \eqref{boot;Z}, \eqref{rmk;phi}, \eqref{eq;weightboot;Z}, the above two estimates,   and the $G_W$ estimate \eqref{estimate g} to the forcing term \eqref{weightderiforc}, we obtain
\begin{equation*}
\begin{split}
	|F_{\widetilde{R}}|\ls& (1+y^2)^\frac{1}{3}e^{-\frac{s}{2}} \left|e^{\frac{3}{2}s}\partial^2_y\Phi-\frac{\beta\partial_y Z}{2}\right|+\beta e^{-s}+(1+y^2)^\frac{1}{3}|G_W||\partial^2_y\overline{W}|\\
	&+ (1+y^2 )^\frac{1}{3} \left(|\dot{\tau}|\left|\partial_y\overline{W}\right|+\left|\partial_yG_W\right|\right)|\partial_y\overline{W}|+ (1+y^2)^\frac{1}{3}\left|\frac{\dot{\tau}\overline{W}+\widetilde{W}}{1-\dot{\tau}}\right|\left|\partial^2_y\overline{W}\right|\\
	\ls &  Me^{-s}+ Me^{-\frac{s}{2}}+Me^{-s}+(1+y^2)^{-\frac{1}{2}}\left( e^{-\frac{s}{2}}+Me^{-\frac{3}{2} s}|y|\right) \\
	&+Me^{- s}+Me^{-s} +\left(\epsilon^\frac{1}{6}(1+y^2)^{\frac{1}{6}}+Me^{-s}\right)\\
	\ls & Me^{-\frac{ s}{2}}+\epsilon^\frac{1}{6} (1+y^2)^{-\frac{1}{3}}, 
\end{split}
\end{equation*}
for $ l\leq |y|$. Combining the above estimate and the bound \eqref{eq;W;traject} for $l\leq |\psi_W|$, we get
\begin{equation}\label{eq;weight;force}
	\begin{split}
		&\int^s_{s_0}\left| F_{\widetilde{R}}\circ\psi_W(y_0,s')\right|ds'\ls\int^s_{s_0}M  e^{-\frac{s'}{2}}+\epsilon^\frac{1}{6} \left(1+|\psi_W|^2(y_0,s')\right)^{-\frac{1}{3}}ds' \ls M \epsilon^\frac{1}{2}+\epsilon^\frac{1}{6} \ln\frac{1}{l}.
	\end{split}
\end{equation}
From the initial data \eqref{ini;weight1}, we have
\begin{equation}\label{self ini1}
\left|(1+y^2)^\frac{1}{3}\partial_y\widetilde{W}(y,-\ln\epsilon)\right|\leq \frac{1}{3}\epsilon^\frac{1}{5},\quad \hbox{for all}\  l\leq  |y|\leq  \frac{5}{4}M\epsilon^{-\frac{1}{2}}+\epsilon^{-\frac{3}{2}}.
\end{equation}
Applying the same argument as in the proof of the
$(1+y^2)^{-\frac16}\widetilde{W}$ estimate, together with
the damping estimate \eqref{eq;weight;damping},
the forcing estimate \eqref{eq;weight;force},
the initial data \eqref{self ini1},
and the bound \eqref{eq;tildeW},
we obtain
%
%
%
%Using the same argument as in the proof of the
%$(1+y^2)^{-\frac16}\widetilde{W}$ estimate, together with the damping estimate \eqref{eq;weight;damping}, the forcing estimate \eqref{eq;weight;force}, the initial data \eqref{self ini1}, and the bound \eqref{eq;tildeW}
%we obtain
\begin{equation*}
	\begin{split}
		&\left|\left((1+y^2)^\frac{1}{3} \partial_y\widetilde{W}\right)\circ\psi_W(y_0,s)\right|\\
		\leq & \left|\left((1+y^2)^\frac{1}{3} \partial_y\widetilde{W}\right)\circ\psi_W(y_0,s_1)\right|\exp\left(-\int^s_{s_1}D_{\widetilde{R}}\circ\psi_W(y_0,s')ds'\right)\\
		&+\int^s_{s_1}\left| F_{\widetilde{R}}\circ\psi(y_0,s')\right|\exp\left(-\int^s_{s'}D_{\widetilde{R}} \circ\psi(y_0,s{''})ds''\right)ds'\\
		\ls & \max\left\{\frac{1}{3}\epsilon^\frac{1}{5},3\epsilon^\frac{1}{5}l^3\right\}l^{-45}e^{6M\epsilon^\frac{1}{2} }+\left(M\epsilon^\frac{1}{2}+\epsilon^\frac{1}{6}\ln\frac{1}{l}\right) l^{-45}e^{6M\epsilon^\frac{1}{2} }\leq \frac{1}{2}\epsilon^\frac{1}{7},
	\end{split}
\end{equation*}
for sufficiently small $\epsilon>0$. This completes the proof.
\end{proof}

\subsection{$(1+y^2)^\frac{1}{3}\partial_yW$ estimate for $|y|\geq e^{\frac{3}{2}s}$.}
\begin{proposition}\label{end prop}
	 We have
	\begin{equation*}
	  (1+y^2)^{\frac{1}{3}} \left|\partial_y W(y,s)\right|\leq \frac{3}{4}, 
\end{equation*} 
for all $ |y|\geq e^{\frac{3}{2}s} $.
\end{proposition}

\begin{proof}
$\bullet$ \textbf{($(1+y^2)^\frac{1}{3}\partial_y W$ estimate)} Let $\widetilde{R}:=(1+y^2)^{\frac{1}{3}}\partial_yW$ in \eqref{eq;variation}. Then from \eqref{higher order W}, we get the following damping term $D_{\widetilde{R}}$ and the forcing term $F_{\widetilde{R}}$:
\begin{align}\label{finaldampW}
D_{\widetilde{R}}=&1+\frac{\beta e^{-s}}{(1-\dot{\tau})(1+\alpha)}+\frac{1}{1-\dot{\tau}}\partial_yW+\frac{e^{\frac{s}{2}}(1-\alpha)}{(1-\dot{\tau})(1+\alpha)}\partial_yZ -\frac{2y}{3(1+y^2)}\left( G_W+\frac{3}{2}y+\frac{W}{1-\dot{\tau}}\right)  
\end{align}
and
\begin{align}\label{finalforceW}
	F_{\widetilde{R}}=&(1+y^2)^\frac{1}{3}\frac{2 e^{-\frac{s}{2}} }{(1-\dot\tau)(1+\alpha)}\left(e^{\frac{3}{2}s}\partial^2_y\Phi-\frac{\beta\partial_y Z}{2}\right).
\end{align}
For $e^{\frac{3}{2}s}\leq |y|$, using Lemma \ref{lem;sbur} (2), the bounds \eqref{rmk;tau}, \eqref{boot;W},  \eqref{boot;Z}, \eqref{eq;Wy}, \eqref{wholeW}, and the $G_W$ estimate \eqref{estimate g} yields
\begin{align*}
&\left|\frac{1}{1-\dot{\tau}} \partial_yW+\frac{e^{\frac{s}{2}}(1-\alpha)}{(1-\dot{\tau})(1+\alpha)}\partial_yZ -\frac{2y^2}{3(1+y^2)}\left(\frac{ G_W}{y}+\frac{W}{y(1-\dot{\tau})}\right) \right|\\
\leq & 3(1+y^2)^{-\frac{1}{3}}+ 2M e^{- s }+\frac{2}{3}\left(2\frac{e^{-\frac{s}{2}}}{|y|}+2Me^{-\frac{3}{2} s}\right)+\frac{Me^\frac{s}{2}}{(1+y^2)^\frac{1}{2}}\leq  4Me^{- s}.
\end{align*}
Applying the above estimate to the damping term \eqref{finaldampW} for $e^{\frac{3}{2}s}\leq|y|$, we obtain
\begin{align*}
 D_{\widetilde{R}}\geq&\frac{1}{1+y^2}-4Me^{-s}\geq -4Me^{-s}.
\end{align*}
Combining the above estimate and \eqref{eq;W;traject} for $e^{\frac{3}{2}s'}\leq|\psi_W|$ gives
\begin{equation}\label{final damp}
\begin{split}
		-\int^s_{s_0}D_{\widetilde{R}}\circ\psi_W(y_0,s')ds'\leq 4M\int^s_{s_0}e^{- s'}ds'\leq 4M \epsilon.
	\end{split}
\end{equation}
Using the bounds \eqref{rmk;phi} and \eqref{eq;weightboot;Z} to the forcing term \eqref{finalforceW}, we get 
\begin{equation*}
\begin{split}
	|F_{\widetilde{R}}|\leq& (1+y^2)^\frac{1}{3}2 e^{-\frac{s}{2}} \left|e^{\frac{3}{2}s}\partial^2_y\Phi-\frac{\beta\partial_y Z}{2}\right|\ls  M^\frac{1}{2} e^{-s}+Me^{-\frac{s}{2}}\ls Me^{-\frac{s}{2}}.
\end{split}
\end{equation*}
This leads to
\begin{equation}\label{final force}
	\begin{split}
		&\int^s_{s_0}\left| F_{\widetilde{R}}\circ\psi_W(y_0,s')\right|ds'\ls\int^s_{s_0}M  e^{-\frac{s'}{2}} \leq \epsilon^\frac{1}{4},
	\end{split}
\end{equation}
where $M$ is absorbed into $\epsilon^\frac{1}{4}$.  Applying a similar argument as in the proof of the $\widetilde{\sigma}$ estimate in Proposition
\ref{lem;rhou}, with $s_1$ defined for the trajectory $\psi_W$, together with
the damping estimate \eqref{final damp}, the forcing estimate
\eqref{final force}, and the bounds \eqref{ini;weight2} and
\eqref{eq;indfinalrevision}, we obtain
\begin{align*}
&\left|\left((1+y^2)^\frac13 \partial_yW\right)\circ\psi_W(y_0,s)\right|\\
\leq{}\,&\left|\left((1+y^2)^\frac13 \partial_yW\right)\circ\psi_W(y_0,s_1)\right|\exp\left(-\int_{s_1}^{s}D_{\widetilde{R}}\circ\psi_W(y_0,s')\,ds'\right)\\
&\quad+\int_{s_1}^{s}\left|F_{\widetilde{R}}\circ\psi_W(y_0,s')\right|\exp\left(-\int_{s'}^{s}D_{\widetilde{R}}\circ\psi_W(y_0,s'')\,ds''\right)\,ds'\\
\leq{}\,&\left(\frac{1}{2}+\epsilon^{\frac{1}{7}}\right)e^{4M\epsilon}+\epsilon^{\frac14}e^{4M\epsilon}\leq\frac34,
\end{align*}
provided that
\[
|\psi_W(y_0,s')|\ge e^{\frac{3}{2} s'}, \qquad \ s'\in [s_1,s],
\]
for sufficiently small $\epsilon>0$. This completes the proof.

%Applying the same argument as in the proof of the $\widetilde{\sigma}$ estimate in Proposition
%\ref{lem;rhou}, with $s_1$ defined for the trajectory $\psi_W$, and using \eqref{final damp} and \eqref{final force}, together with the bounds \eqref{ini;weight2} and \eqref{eq;tildeW;inter},

%Using \eqref{final damp} and \eqref{final force}, along with the bounds  \eqref{ini;weight2} and \eqref{eq;tildeW;inter}, a similar argument as in the $\widetilde{\sigma}$ estimate proof   in Proposition \ref{lem;rhou} {\color{cyan}with $s_1$ with respect to $\psi_W$ } yields
%\begin{align*}
%		&\left|\left((1+y^2)^\frac{1}{3} \partial_y W\right)\circ\psi_W(y_0,s)\right|\\
%		\leq & \left|\left((1+y^2)^\frac{1}{3} \partial_y W\right)\circ\psi_W(y_0,s_1)\right|\exp\left(-\int^s_{s_1}D_{\widetilde{R}}\circ \psi_W(y_0,s')ds'\right)\\
%		&+\int^s_{s_1}\left| F_{\widetilde{R}}\circ\psi_W(y_0,s')\right|\exp\left(-\int^s_{s'}D_{\widetilde{R}} \circ\psi(y,s{''})ds''\right)ds'\\
%		\leq &\max\left\{\frac{1}{2} ,\left(\frac{7}{20} +\epsilon^\frac{1}{7}\right)\right\}e^{4M\epsilon} +\epsilon^\frac{1}{4}e^{4M\epsilon}\leq \frac{3}{4},
%\end{align*}
%where $\color{cyan}|\psi_W|(y_0,s_2)\geq e^{\frac{3}{2}s_2}$, with   $\epsilon>0$ sufficiently small.
\end{proof}

\section{Proof of main result}\label{sec;finalmain}

We show Theorem \ref{delicate main} with Proposition \ref{prop;mod} -\ref{end prop} in the bootstrap argument. 

\subsection{Pointwise estimates of $w,z$ and $\phi$ at the blow-up time $T^*$.}\label{final section}
First, we address $z_x(x,T^*)$ for all $x$ and $w_x(x,T^*)$ for all $x\neq x^*(=\lim_{t\to T^*} \xi(t))$.
 Differentiating the equations \eqref{subrie1} and \eqref{subrie2} with respect to $x$, we have
\begin{subequations}
		\begin{align}
		\partial_{t} w_x+\left(\frac{\beta}{1+\alpha}+w_x+\frac{1-\alpha}{1+\alpha}z_x \right) w_x+\left(w+\frac{1-\alpha}{1+\alpha}z-\frac{\kappa_0}{1+\alpha} \right) w_{xx} =&\frac{2}{1+\alpha}\left(\phi_{xx}-\frac{\beta}{2}z_x\right),\\
		\partial_tz_x+\left(\frac{\beta}{1+\alpha}+\frac{1-\alpha}{1+\alpha}w_x+ z_x\right)z_x+\left(\frac{1-\alpha}{1+\alpha}w+z-\frac{\kappa_0}{1+\alpha} \right)z_{xx} =&\frac{2}{1+\alpha}\left( \phi_{xx}-\frac{\beta}{2}w_x\right).\label{def z}
		\end{align}
\end{subequations}
Since similar arguments hold for $w_x(x,T^*)$ when $x\neq x^*$, we focus solely on $z_x(\cdot,T^*)$. Define the characteristic curve $\psi$ as follows:
\begin{align}\label{traj z}
\frac{d}{dt} \psi(x,t) =\left(\frac{1-\alpha}{1+\alpha} w+z-\frac{\kappa_0}{1+\alpha} \right)(\psi(x,t),t) \, , \qquad \psi(x,0) = x \, .
\end{align}
For simplicity, we denote $\psi(x,t)$ as $\psi(t)$. Then, we prove the pointwise convergence; that is, $z_x(\psi(T^*),T^*)=\lim_{t\to T^*}z_x(\psi(t),t)$. Applying Lemma \ref{max prin1} with this flow $\psi$ to \eqref{def z} implies that 
it suffices to show that 
\begin{equation}\label{l1 estimate}
\int^{T^*}_0 | z_x(\psi(t),t)| dt<\infty \quad \hbox{and} \quad \int^{T^*}_0 | w_x(\psi(t),t)| dt<\infty.
\end{equation}
Since $\sup_{0\leq t<T^*}\| z_x\|_{L^\infty_x} \leq 2M$ from the bound \eqref{boot;Z}, we consider only the $L^1_t$ estimate of $w_x(\psi(t),t)$.
From the bound \eqref{eq;Wy}, we get
\begin{align}\label{bound WW}
| w_x (\psi(t),t)| = e^s \left|\partial_y W\left( e^{\frac{3}{2}s}(\psi(t) - \xi(t)) , s\right)\right|  \leq   \frac{2}{|\psi(t) - \xi(t)|^{\frac{2}{3}}  }\, .
\end{align}
First, we analyze the case  $\psi(T^*)\neq\xi(T^*)$.  From the time continuity of trajectories, we can choose $c,c'>0$ such that  $|\psi(t)-\xi(t)|\geq c'$  for $c\leq t\leq T^*$. Thus combining \eqref{bound WW} with the previous statement, $|w_x(\psi(t),t)|$ is uniformly bounded for $0\leq t<T^*$. What is left to show the case $\psi(T^*)=\xi(T^*)$.  From the definitions of $\psi$ and $\xi$ in \eqref{traj z} and  \eqref{form xi}, respectively, we obtain
 \begin{align*}
&\psi(t) - \xi(t)\\
 =& \int_{t}^{T^*}  \dot{\xi}(t')- \frac{1-\alpha}{1+\alpha} w(\psi(t'),t')-z(\psi(t'),t') +\frac{\kappa_0}{1+\alpha} dt' \notag\\
=& \int_{t}^{T^*}  -(1-\dot \tau) e^{-\frac{s'}{2}} \frac{F_W^{0,(2)}}{\partial^3_yW^0}(s')+\kappa(t')+\frac{1-\alpha}{1+\alpha} Z^0(s') -z(\psi(t'),t') \\
&\quad- \frac{1-\alpha}{1+\alpha} w(\psi(t'),t')  dt' \\
=& \int_{t}^{T^*}\frac{2\alpha}{1+\alpha} \kappa(t')  dt'+    \int_{t}^{T^*} \frac{1-\alpha}{1+\alpha} Z^0(s') - Z\left(e^{\frac{3}{2}s'}(\psi(t')-\xi(t')) ,s'\right)  dt' \notag\\
& - \int_{t}^{T^*}  \frac{1-\alpha}{1+\alpha} e^{-\frac{s'}{2}} W\left(e^{\frac{3}{2}s'}(\psi(t')-\xi(t')) ,s'\right) dt'-\int^{T^*}_t (1-\dot \tau) e^{-\frac{s'}{2}} \frac{F_W^{0,(2)}}{\partial^3_y W^0}(s') dt' \\
=& I_1(t) + I_2(t) - I_3(t)-I_4(t).
\end{align*}
Using the bound \eqref{kappa estimate} yields
 \begin{equation*}
 	I_1(t) \geq \frac{15}{2} (T^*-t).
 \end{equation*}
From the bound \eqref{boot;Z} for the $Z$ estimate, we have
 \begin{equation*}
| I_2(t)|\leq (2+16M\epsilon)(T^*-t).
 \end{equation*}
Since  $|W|(y,s)\leq (1+\epsilon^\frac{1}{6})(1+y^2)^\frac{1}{6}$ for $|y|\leq e^{\frac{3}{2}s}$ from  \eqref{eq;tildeW} and \eqref{eq;tildeW;inter}, we obtain 
 \begin{align*}
|I_3(t)|\leq & \int_{t}^{T_*}  \left| \frac{1-\alpha}{1+\alpha} e^{-\frac{s'}{2}} W\left(e^{\frac{3}{2}s'}(\psi(t')-\xi(t')) ,s'\right)\right|dt'\\
\leq &  2\int^{T^*}_t e^{-\frac{s'}{2}}(1+ e^{\frac{s'}{2}} |\psi(t')-\xi(t')|^\frac{1}{3}) dt'\\
\leq & 2\epsilon^\frac{1}{2}(T^*-t)+2\int^{T^*}_t  |\psi(t')-\xi(t')|^\frac{1}{3}    \leq3\epsilon^\frac{1}{2} (T^*-t),
\end{align*}
for $t$ sufficiently near $T^*$. From the bounds \eqref{mod accel} and \eqref{sec force W}, we get
 \begin{align*}
&|I_4(t)|\leq  \int_{t}^{T^*}\left|  (1-\dot \tau) e^{-\frac{s'}{2}} \frac{F_W^{0,(2)}(s')}{\partial^3_yW^0(s')} \right|dt'\leq 3\int^{T^*}_te^{-s'}dt'\leq 3\epsilon(T^*-t).
\end{align*}
Summing these estimates yields
\begin{equation*}
	\psi(t) - \xi(t)\geq 2(T^*-t),
\end{equation*}
for $t$ sufficiently near $T^*$. Combining the above estimate with \eqref{bound WW} implies \eqref{l1 estimate}.

Next, we prove the continuity of $w_x$ except for $x^*=\lim_{t\to T^*} \xi(t)$ and $z_x$ at $T^*$. Since a similar argument holds for $w_x$, we focus only on $z_x$. Fix $x\in\mathbb{R}$ and  $\epsilon'>0$. Since $\lim_{t\to T^*}z_x(\psi(t),t)=z_x(\psi(T^*),T^*)$  and smoothness of $z$ before $T^*$, we can choose $\delta>0$ such that if $0<|\psi(T^*)-\psi(x',T^*)|$, $T^*-t<\delta$, then
\begin{align*}
 		&\left|z_x(\psi(T^*),T^*)-z_x(\psi(x',T^*),T^*)\right|\\
		\leq &\left|z_x(\psi(T^*),T^*)-z_x(\psi_z(t),t)\right|\\
		&+\left|z_x(\psi(t),t)-z_x(\psi(x',t),t)\right|+\left|z_x(\psi(x',t),t)-z_x(\psi( x',T^*),T^*)\right|\\
		\leq &\epsilon'.
 \end{align*}
Thus continuity is established. In addition, the pointwise convergence of $w$ is obtained from the proof of  Proposition \ref{prop;Westi}. Similarly, the convergence of $z$ follows from \eqref{subrie2}. 

To prove that  the derivative of $z$ at $T^*$ coincides the pointwise limit of $z_x$,  let $f:\mathbb{R}\backslash\{0\}\times[0,T^*]\to \mathbb{R}$ be defined by
\begin{equation*}
f(h,t)=\frac{z(\psi(t)+h,t)-z(\psi(t),t)}{h}.
\end{equation*}
For fixed $|h|>0$, we  split $z_x(\psi(T^*),T^*)-f(h,T^*)$ as
\begin{align*}
z_x(\psi(T^*),T^*)-f(h,T^*)=&z_x(\psi(T^*),T^*)-z_x(\psi(t),t)+z_x(\psi(t),t)-z_x(\psi(t)+h,t)\\
&+z_x(\psi(t)+h,t)-f(h,t)+f(h,t)-f(h,T^*).
\end{align*}
Since $f$ and $z_x(\cdot,T^*)$ are continuous on their respective domains, it follows from the previous arguments that:\begin{equation*}
\begin{split}
z_x(\psi(T^*),T^*)=&\lim_{h\to0}\frac{z(\psi(T^*)+h,T^*)-z(\psi(T^*),T^*)}{h}.
\end{split}
\end{equation*}
An analogous result for $w_x$, except for $x^*$, holds. Applying these results and the basic property of the heat kernel to \eqref{subrie3}, we deduce $\phi\in C([0,T^*];C^2(\mathbb{R}))$.
\subsection{Proof of Theorem \ref{delicate main}.}
\begin{enumerate}[(1)]
\item From the continuation criterion \eqref{continucri} in Theorem \ref{local well1} with \eqref{boot;Z} and \eqref{eq;Wy}, the $H^m$-stability prior to singularity formation is established.
\item By Proposition \ref{prop;mod}, we obtain the blow-up time $T^*(\leq\frac{3}{2}\epsilon)$ and the blow-up position $x^*=\lim_{t\to T^*}\xi(t)$ with $|
x^*|\leq 6M\epsilon$.
\item By the constraint $\partial_yW(0,s)=-1$ in \eqref{W constraint} and $w(x,t)=e^{-\frac{s}{2}}W+\kappa(t)$ in \eqref{self similar}, we get
\begin{equation*}
\lim_{t\to T^*}w_x(\xi(t),t)=\lim_{s\to \infty}e^{s}\partial_yW\left(0,s\right)=-\infty.
\end{equation*}
Since $\tau(t)-t=\int^{T^*}_t(1-\dot{\tau}(t'))dt'$ and by \eqref{mod accel}, we have $ \tau(t)-t\thicksim T^*-t$. Then, from \eqref{eq;Wy}, we obtain
\begin{align*}
\| w_x(\cdot,t)\|_{L^\infty}=\frac{1}{\tau(t)-t}\left\|\partial_y W(\cdot,s)\right\|_{L^\infty}\thicksim\frac{1}{T^*-t},
\end{align*}
for all  $0\leq t<T^*$. Combining the bound $\|\partial_yZ(\cdot,s)\|_{L^\infty}\leq 2M$ in \eqref{boot;Z} and the definition of Riemann-type variables in \eqref{riemann;type} yields 
\begin{equation*}
 		\lim_{t\to T^*}  q_x(\xi(t),t)=\lim_{t\to T^*} u_x(\xi(t),t)=-\infty.
 	\end{equation*}
\item If $0<|x-x^*|< 1$, then there exists $t_1\in[0,T^*)$ such that $e^{-\frac{3}{2}s} l\leq\frac{1}{2}|x-x^*|\leq |x-\xi(t)|\leq 1$ for all $t\in[t_1,T^*)$. i.e. $\frac{1}{2}|x-x^*|e^{\frac{3}{2}s}\leq |y|\leq e^{\frac{3}{2}s}$. Therefore, from Lemma \ref{lem;sbur} $(4)$ and the bound \eqref{eq;tildeW;inter}, we get 
\begin{align*}
-\frac{2}{5} |x-\xi(t)|^{-\frac{2}{3}}\leq e^{s}\left(-\frac{7}{20}-\epsilon^\frac{1}{7}\right)|y|^{-\frac{2}{3}}\leq&  w_x(x,t)\\
 \leq& e^{s}\left(-\frac{1}{4}+\epsilon^\frac{1}{7}\right)|y|^{-\frac{2}{3}}\leq -\frac{1}{5} |x-\xi(t)|^{-\frac{2}{3}},
\end{align*}
 for $t_1\leq t<T^*$. So we obtain $w_x(x,T^*)\thicksim -|x-x^*|^{-\frac{2}{3}}$  for  $0<|x-x^*|< 1$. With \eqref{eq;W;final}, this leads to uniqueness of the blow-up point.
\item  Due to the structure of $w_x$ in (4),  $w$ has $C^\frac{1}{3}$ regularity in time up to $T^*$.
Then, applying  the bounds \eqref{boot;Z} and \eqref{eq;weightboot;phi}, we obtain  
\begin{equation*}
\sup_{t\in[0, T^*]} \left( \| w(t)\|_{ C^\frac{1}{3}}+ \| z(t)\|_{ C^1}+ \| \phi(t)\|_{ C^{2 }}  \right) <  +\infty.
\end{equation*}
In Section \ref{final section}, we already showed that $w$ is $C^1$ differentiable at $x\neq x^*$.

\end{enumerate}

\section*{Acknowledgments}
 The author  thanks the anonymous referees for their helpful and considerate comments on the construction and readability of the paper. The work of W. Lee is supported by NRF grants (Nos.  2022R1A2C1002820, and RS-2024-00453801).

%\begin{figure}[htbp]
%  \centering
%  \includegraphics[width=0.5\textwidth]{fig2}
%  \caption{}\label{fig2}
%\end{figure}

%\begin{table}
%  \centering
%\begin{tabular}{||||}
% \hline
%   &  &  \\
%   &  &  \\
%   &  &  \\
%  \hline
%\end{tabular}
 %\caption{}\label{}
%\end{table}


%          \begin{thebibliography}{}

%          % and use \bibitem to create references.

%           \bibitem{mnemonic}Author, {\em Article's title}, Journal Name,
%         Volume Number(Issue Number):page numbers, year.

          % Format for Journal Reference. For example

 %         \bibitem{taubes1} C. Taubes, {\em The Seiberg-Witten invariants
  %        and
   %       symplectic forms}, Math. Res. Letters, 1:809--822, 1994.
    %      \end{thebibliography}

          \end{document}